\documentclass{article}
\usepackage{amssymb,amsbsy,amsmath,amsfonts,amscd}
\usepackage{latexsym}

\usepackage{tocvsec2}
\usepackage{xy} 
\xyoption{all} 

\usepackage{latexsym}

\newcommand{\lie}[1] {\mathfrak{#1}}  
\newcommand{\bb}[1]{{\mathbb #1}}    

\newcommand{\bbR}{{\bb R}}
\newcommand{\bbA}{{\bb A}}
\newcommand{\bbZ}{{\bb Z}}

\newcommand{\bbC}{{\bb C}}

\newcommand{\GL}{{\rm GL}}
\newcommand{\PGL}{{\rm PGL}}
\newcommand{\SL}{{\rm SL}}

\newcommand{\SP}{{\rm Sp}}

\newcommand{\Aff}{{\rm Aff}}
\newcommand{\Aut}{{\rm Aut}}
 
\newcommand{\Mat}{{\rm Mat}}

\newcommand{\Hom}{{\rm Hom}}

\newcommand{\Diff}{{\rm Diff}}
\newcommand{\Ad}{{\rm Ad}}
\newcommand{\ad}{{\rm ad}}
\newcommand{\N}{{\rm N}}  
\newcommand{\Z}{{\rm Z}} 

\newcommand{\E}{{\rm E}}
\newcommand{\A}{{\rm A}}

\newcommand{\bU}{{\mathbf{U}}}
\newcommand{\bG}{{\mathbf{G}}}
\newcommand{\bH}{{\mathbf{H}}}
\newcommand{\bL}{{\mathbf{L}}}

\newcommand{\bT}{{\mathbf{T}}}

\newcommand{\bX}{{\mathbf{X}}}

\newcommand{\cA}{{\mathcal{A}}}
\newcommand{\cB}{{\mathcal{B}}}
\newcommand{\cC}{{\mathcal{C}}}
\newcommand{\cD}{{\mathcal{D}}}

\newcommand{\cU}{{\mathcal{U}}}

\newcommand{\lu}{{\lie{u}}}
\renewcommand{\lg}{{\lie{g}}}

\newcommand{\lz}{{\lie{z}}}

\newcommand{\tensor}[1]{\otimes_{#1}}

\newcommand{\id}{{\rm id}}

\newcommand{\rank}{{\rm rank }\, }
\renewenvironment{matrix}[1] {\left( \begin{array}{#1}}{\end{array}\right)}
\newcommand{\ac}[1]{\overline{#1}}

{ \par  \medskip  \noindent  {\bf Definition} \hspace{0.5em} }%
{\par \medskip } 

\newtheorem{example}{Example}[section]

\newenvironment {proof}%
{ \noindent {\em Proof. }}%
{\hspace*{\fill}$\Box$\par \medskip } 
\newenvironment{prf}[1]%
{ \noindent {\em #1 \hspace{0.5em}} }%
{\hspace*{\fill}$\Box$\par \medskip } 
\newenvironment {remark}%
{{\em Remark \hspace{0.5em}} }%
{\par \medskip }

\newtheorem{proposition}{Proposition}[section]
\newtheorem{definition1}[proposition]{Definition}
\newtheorem{theorem}[proposition]{Theorem}
\newtheorem{lemma}[proposition]{Lemma}
\newtheorem{corollary}[proposition]{Corollary}
\newenvironment{definition}%
{ \begin{definition1} \rm }%
{ \end{definition1}}

\numberwithin{equation}{section}

\newcommand{\SAff}{{\rm SAff}}

\newtheorem{question}{Question}

\newcommand{\ra}{\rightarrow}
\newcommand{\bbE}{{\mathbb E}}
\newcommand{\cK}{{\mathcal K}}

\newcommand{\Mid}{{\, \Big\vert \;}}
\newcommand{\Isom}{{\rm Isom}}
\newcommand{\Vect}{{\rm Vect}}
\renewcommand{\div}{{\rm div}}
\newcommand{\trace}{{\mathop{\rm trace\,}}}

\newcommand{\aut}{\lie{aut}}
\newcommand{\isom}{\lie{o}}
\renewcommand{\o}{\lie{o}}
\newcommand{\sym}{\lie{s}}
\newcommand{\g}{{<  ,   >}}
\newcommand{\aff}{{\lie{aff}}}
\newcommand{\lh}{{\lie{h}}}

\newcommand{\lra}{{\, \longrightarrow \, }}
\newcommand{\bell}{\bar{\ell}}

\title{Prehomogeneous Affine Representations and\\ Flat Pseudo-Riemannian Manifolds}

\author{Oliver Baues \thanks{e-mail: baues@math.uni-karlsruhe.de} \\
Institut f\"ur Algebra und Geometrie\\ 
Universit\"at Karlsruhe\\ 
D-76128 Karlsruhe
}

\bigskip
\bigskip
\bigskip
\bigskip

\begin{document}

\maketitle

\begin{abstract} 
The theory of flat Pseudo-Riemannian manifolds and
flat affine manifolds is closely connected to the topic
of prehomogeneous affine representations of Lie groups.
In this article, we exhibit several aspects of this correspondence.
At the heart of our presentation is a development of the theory
of characteristic classes and characters of prehomogeneous affine 
representations. We give applications concerning 
flat affine, as well as Pseudo-Riemannian
and symplectic affine flat manifolds. \\


\hrule 

\hspace{1cm} \\

\end{abstract}

\medskip

\tableofcontents

\newpage

\section{Introduction}  \label{sect:intro}

\subsection{Flat Pseudo-Riemannian and flat affine manifolds} \label{sect:overview}
A Pseudo-Riemannian manifold $(M,\g)$ is a smooth manifold
$M$ which is endowed with a possibly indefinite metric $\g$ 
on its tangent bundle $TM$.  
We let the expression $s(\g) = (n_{+}, n_{-})$, 
$n_{+} + n_{-} = \dim M$, denote the signature of $\g$. 
A positive definite metric $\g$ has signature $(n,0)$ and
is traditionally called a \emph{Riemannian metric}. If the signature of $\g$
is $(n-1,1)$, $\g$ is called a \emph{Lorentzian metric}. 

Every Pseudo-Riemannian manifold $(M,\g)$ has a unique
torsion-free connection $\nabla_{\g}$ on its tangent bundle
$TM$ which has the property that the metric tensor $\g$ 
is parallel for $\nabla_{\g}$.  
The connection $\nabla_{\g}$ is called the Levi-Civita connection for $\g$. 
Given vector fields $X,Y, Z$ on $M$, the tensorial expression, 
$$ R_{\nabla}{}(X,Y) Z = \nabla_{X} \nabla_{Y} Z -  \nabla_{Y} \nabla_{X} Z -\nabla_{[X,Y]} Z $$ is called the \emph{curvature tensor} for a connection $\nabla$.

The 
metric $\g$ is called \emph{flat} if 
the curvature $R_{\nabla_{\g}}$ vanishes everywhere on $M$. 
Thus, in particular, \emph{flat Pseudo-Riemannian manifolds} carry 
a torsion-free and flat connection on their tangent bundle.  A manifold 
$M$ together with a torsion-free and flat connection $\nabla$ is 
called a \emph{flat affine manifold}. Naturally, flat affine manifolds 
$(M, \nabla)$ share many of their geometric properties with the more restricted
class of flat Pseudo-Riemannian manifolds.   

\subsubsection{Global models and their quotients}
There does exist an abundance of examples of simply connected flat Pseudo-Riemannian and flat affine manifolds. Indeed, every open subset $U$ of Euclidean space $\bbR^n$, defines a flat Pseudo-Riemannian manifold (of any signature), 
and, every simply connected flat Pseudo-Riemannian manifold $(\tilde{M},\g)$ or flat affine manifold $(\tilde{M}, \nabla)$ is obtained by pulling back the flat structures along
a local diffeomorphism $\Phi: \tilde{M} \ra U$ onto an open subset 
of $\bbR^n$. 
 
A simply connected Pseudo-Riemannian manifold $(\tilde{M},\g)$ will be called a 
\emph{global model space} for flat Pseudo-Riemannian manifolds. By elementary covering theory,  every Pseudo-Riemannian manifold $(M,\g)$ is obtained as a quotient space of a global model $(\tilde{M},\g)$ by a \emph{properly discontinuous} group 
$\Gamma$ of isometries of $(\tilde{M},\g)$.  Similarly, any flat affine manifold $(M,\nabla)$ is a quotient of a simply connected manifold $(\tilde{M}, \nabla)$ by a group
of connection preserving diffeomorphisms. 

In principle, the study of flat Pseudo-Riemannian manifolds 
breaks into two parts, namely the determination of \emph{interesting}
global models $(\tilde{M},\g)$ (apart from the standard complete
model $\bbE^s$), and the study of the quotient 
spaces $(M,\g)$ which are modelled on $(\tilde{M},\g)$. 

A particular interesting class of model spaces will be furnished
by \emph{homogeneous} flat Pseudo-Riemannian manifolds. These are
flat spaces which admit a transitive group of isometries. More
generally,  \emph{homogeneous flat affine manifolds}, that is,
flat affine manifolds with a transitive group of affine
transformations constitute a natural class of models.

\subsubsection{Completeness and Pseudo-Euclidean space forms}
A Pseudo-Riemannian manifold $(M,\g)$ is called \emph{complete} if \emph{every} 
geodesic curve for the connection $\nabla_{\g}$ can be extended to infinity. 
If the flat manifold $(M,\g)$ is
complete, the \emph{Killing-Hopf theorem} asserts, that the universal Pseudo-\-Rie\-mann\-ian covering space $(\tilde{M},\g)$ for $(M,\g)$ is \emph{isometric} to the 
Pseudo-Euclidean space $\bbE^s= (\bbR^n, < , >_{s})$ of signature $s=s(\g)$, where
$\g_{s}$ denotes the standard representative for a scalar product
of signature $s$.  

In particular, for fixed signature $s$, there exists,
up to isometry, a unique simply connected  \emph{complete model
space} $\bbE^s$ for flat Pseudo-Riemannian manifolds of signature $s$. 

As another consequence, every complete flat Pseudo-Riemannian manifold 
$(M,\g)$ of signature $s$ is isometric to a quotient of $\bbE^s$
by a \emph{properly discontinuous subgroup}  $\Gamma$ of isometries, acting 
without fixed points on $\bbE^s$. Such quotient manifolds are 
also called \emph{Pseudo-Euclidean space forms}. 

\subsection{Flat Riemannian manifolds}
By the Hopf-Rinow theorem  (which holds solely in Riemannian geometry), 
the compactness of a Riemannian manifold implies its completeness. Similarly, 
the equivalence of metric and geodesic completeness for Riemannian manifolds
implies that every homogeneous Riemannian manifold is complete.  In fact, every homogeneous flat Riemannian manifold is a quotient of $\bbR^n$ by a
group of translations. 
In particular, every compact flat Riemannian manifold or homogeneous 
flat Riemannian manifold is a quotient of Euclidean space by
a \emph{discontinuous} group of isometries.
Thus, the theory of 
flat Riemannian manifolds concerns mostly the study of complete space forms,
and  it is roughly equivalent to the study of discontinuous subgroups of the 
isometry group $E(n)$ of Euclidean space.


The structure of discontinuous subgroups of $E(n)$ is rather well understood, 
by the famous three \emph{theorems of  Bieberbach} \cite{Bieberbach_1, Bieberbach_2},
dating from around 1910. 
According to Bieberbach, every discrete subgroup $\Gamma$ of
$E(n)$ is finitely generated and contains an abelian 
subgroup of finite index. If $\Gamma$ acts with compact quotient
space, then the  subgroup of translations is of finite index in $\Gamma$. 
In particular, every compact Euclidean space form is finitely and isometrically 
covered by a flat torus $\bbR^n /\Lambda$,  where
$\Lambda$ is a lattice of translations. Moreover, Bieberbach proved, that 
in every dimension $n$, there exist only finitely many compact flat 
Riemannian manifolds up to affine equivalence. Thus, the class of compact flat Riemannian manifolds is rather restricted from a topological point of view. 

The determination of 
compact Riemannian space forms and their geometric properties has a
long tradition, and remains a subject of geometric and algebraic (it is related
to the study of integral representations of finite groups) interest. 
Recent contributions concern, for example, 
isospectrality phenomena \cite{DR, MiPo, Schie}, and spin 
structures \cite{DSS,Pf} on flat Riemannian manifolds. 

The theory of \emph{almost flat} Riemannian manifolds as developed by Gromov
\cite{Grom_1}, in a sense, extends the theory of Bieberbach to a much more
general context. See \cite{Grom_1, BuKa}.  Here, the role of 
flat tori is taken over by compact \emph{nilmanifolds} (these are quotient spaces
of nilpotent Lie groups). Incidentally,  nilmanifolds  also appear (at least conjecturally)
as the fundamental building blocks for compact flat Pseudo-Riemannian manifolds.

\subsection{Flat manifolds of indefinite signature}

Much of the theory of flat Pseudo-Riemannian manifolds aims
to construct an analogy to the theory of Euclidean space forms.
But, as is it turned out, many new phenomena and principle difficulties
arise.  Many of them  constitute still open and difficult research 
questions. 

For example, it is widely believed, that, as in the Riemannian case, only the Pseudo-Euclidean spaces $\bbE^s$ admit compact quotient manifolds. This conjecture is so far
verified only for Lorentzian manifolds, see \cite{Carriere}. 
 
Likewise, the determination of simply connected  \emph{homogeneous model spaces}
 $(\tilde{M},\g)$, and more generally the determination of all \emph{homogeneous 
flat Pseudo-Riemannian manifolds} of a given signature $s$, is an unsolved
problem. 

\subsubsection{Structure theory of compact flat Pseudo-Riemannian manifolds}

\paragraph{Compactness and completeness} 
\label{sect:Markus}
A first stumbling stone for a structure theory is created by the general lack of understanding
about the relationship of compactness and completeness for
flat Pseudo-Riemannian and affine manifolds. In fact, although there are
known examples of compact (non flat) Lorentzian manifolds,
which are incomplete, it is expected (``\emph{Conjecture of Markus}'') that,
as in the Riemannian case, a flat compact Pseudo-Riemannian manifold is indeed also complete. This conjecture was proved by Carri\'ere \cite{Carriere} for the case of
flat Lorentzian manifolds. But it remains unverified for flat 
Pseudo-Riemannian manifolds of arbitrary signature. 

More generally, the \emph{Markus conjecture} (attributed to \cite{Markus}) 
suggests that an orientable
compact affinely flat manifold is complete if and
only if it admits a \emph{parallel volume}. Beside the result of Carri\'ere, 
this  conjecture is verified only in a few special cases. 
For example, \emph{compact  homogeneous} affine manifolds satisfy 
Markus' conjecture. 

Note that, by dividing out a cyclic group of linear dilatations on the affine
space with the point ${0}$ removed, one constructs
a simple example of a compact affinely flat and incomplete 
manifold. In general, the theory of compact flat
affine manifolds is considered as wild and possibly 
untractable, although some results on the topology 
of such manifolds have been obtained in low dimensions.
See \cite{SuThu} for wild examples. 

\paragraph{Group theoretic structure of the fundamental group} 
Assuming completeness,  new osbtacles appear to 
generalise Bieberbach's theory of discontinuous Euclidean 
groups and Euclidean crystallographic groups to 
a theory of affine crystallographic groups. These
concern the group theoretic structure of
the fundamental group of a complete
Pseudo-Riemannian manifold. Conjecturally,
(``\emph{Auslander's conjecture}''), an affine crystallographic
group has a solvable subgroup of finite index. 
The conjecture arose from a paper \cite{Auslanderc} 
of Louis Auslander, which contained a flawed
proof of the even stronger claim that every 
finitely generated discontinuous subgroup of affine transformations
has a solvable subgroup of finite index. 
  
The solvable group replaces the finite index abelian subgroup in the Euclidean case.
In this sense, the \emph{Auslander conjecture} serves as a weak analogue to Bieberbach's
first theorem. Auslander's conjecture is verified in low 
dimensions,  in the Lorentzian case,
and in some other cases (see \cite{Abels} for a survey). 
But it remains one of the main open questions of the subject.

Also, if the assumption of compactness is dropped, 
examples of flat complete three dimensional Lorentzian manifolds 
with a free non-abelian fundamental
group (cf.\ \cite{DG, Mar}) give a counter example to the 
original claim of Auslander.  In this case,
the analogy with the Euclidean theory breaks down.

\paragraph{The classification theory of Pseudo-Euclidean space forms}
Assuming completeness and solvability of
the fundamental group, new difficulties
and phenomena arise for a possible classification
program. For instance, the finitess part of the 
Bieberbach theory breaks down, as well. 
But it admits a weak and rather subtly 
defined replacement, as is described in \cite{GS}. 

The main achievement of the theory, so far,
gives a rather precise and strong link of the theory of 
Pseudo-Euclidean and affine crystallographic
groups with the theory of \emph{left-invariant} flat
Pseudo-Riemannian metrics on Lie groups. (See 
\cite{FriedGoldman,GS} for an exposition.)
For example, a classification theory of compact
flat Lorentzian manifolds is developed in 
\cite{AusM,Fried, GoldmanKamishima_2, GM, Gued}. 

However, despite the achievements
in the Lorentzian case, the structure theory 
for compact complete flat Pseudo-Riemannian
manifolds remains widely open. 

\subsubsection{Homogeneous flat Pseudo-Riemannian manifolds}
Unlike the Riemannian case, a homogeneous
Pseudo-Riemannian manifold need \emph{not} be complete. 
The simplest example of a \emph{non-complete} homogeneous flat
Pseudo-Riemannian manifold is an open orbit of the two-dimensional
non-abelian solvable simply connected Lie group in $\bbE^{1,1}$
(see  \cite[\S 11]{Wolf_4}). In general, the classification of non-complete 
homogeneous flat Pseudo-Riemannian manifolds is not fully understood.

Only if the additional assumption of compactness is made then 
homogeneity implies completeness for Pseudo-Riemannian manifolds 
(in fact, even without the assumption of flatness), see \cite{Hermann, Marsden}. 
Furthermore, in some cases, the group theoretical structure of a homogeneous 
space is directly linked to its completeness properties.  
For example, a flat Pseudo-Riemannian Lie group, more generally a volume
preserving flat affine Lie group, is complete if and only if it is unimodular.  
In this article, we introduce the related new result that a flat affine homogeneous space of 
a nilpotent Lie group is complete if and only if the action preserves
a parallel volume form. In particular, a flat homogeneous Pseudo-Riemannian 
manifold or a flat symplectically homogeneous affine
manifold of a nilpotent group is complete. See section \ref{sect:HAM} and 
section \ref{sect:crittrans} 
for further discussion of these results.

The structure and classification of \emph{complete} homogeneous flat 
Pseu\-do\--Rie\-mannian manifolds is more accessible than the general case. 
Every Riemannian
homogeneous flat manifold is obtained by dividing out a group of translations in
$\bbR^n$. In the Pseudo-Riemannian case  interesting phenomena do occur.

For example, contrasting the Riemannian case, there exists a large class 
of (non-compact) complete homogeneous flat Pseu\-do\--Rie\-mann\-ian manifolds 
with abelian but \emph{non-translational} holonomy groups.  See \cite{Wolf_1, Wolf_2, Wolf_3} for a detailed  study of this examples. 

In general, the holonomy group of a homogeneous
flat Pseudo-Riemannian manifold must be a two-step nilpotent group. 
(See section \ref{sect:homhol_psr} for details on this.)
Non-abelian fundamental groups occur, in particular, as
fundamental groups of \emph{compact} homogeneous flat 
Pseudo-Riemannian manifolds. In fact, there exists a large
class of compact two-step nilmanifolds, which admit an (essentially unique)
homogeneous flat Pseudo-Riemannian metric. The first examples, which
are not homotopy equivalent to a torus arise in dimension six. See 
section \ref{sect:comphPSR} for more details and proofs.   

For general  homogeneous affinely flat manifolds similar results hold. 
By \cite{GH_2}, the affine holonomy group
of a compact affine manifold with parallel volume does not
preserve any proper algebraic subsets in affine space $\bbA^n$. 
As a corollary 
a compact homogeneous affine manifold with parallel volume is complete.   
(An independent proof is given in section \ref{sect:comphom}.)

More generally, if a volume preserving homogeneous affine 
manifold admits a compact Clifford-Klein form then it must be
complete (see section \ref{sect:GHvol} of this article for further discussion).
In particular, if the universal covering of a compact 
flat Pseudo-Riemannian manifold is homogeneous then it must
be complete. It follows that non-complete homogeneous flat Pseudo-Riemannian 
manifolds do not occur as models of compact flat Pseudo-Riemannian 
manifolds.  

%




\subsection{Overview of the article}

This article aims to exhibit  several aspects, which link the theory 
of flat affine and Pseudo-Riemannian manifolds with the topic of
representations of Lie groups and algebraic groups on affine 
space $\bbA^n$. 

In section \ref{sect:intro} we started with an overview on 
some of the main achievements and open problems in the topic of flat 
affine and  Pseudo-Riemannian manifolds. 

In the following section \ref{sect:foundations} 
we discuss flat manifolds 
from the point of view of Thurston's theory of locally homogeneous 
$(X,G)$-manifolds. 
Here we introduce our notation, as well as basic definitions
and methods. 

In section \ref{sect:AffM}, we describe the structure 
of the group of affinities and of the
Lie algebra of affine vector fields on flat affine manifolds.
By the development process, these groups are represented 
as subgroups of the affine group, and relate to certain associative matrix 
algebras.  Compactness  poses strong
restrictions on the symmetries of volume preserving, 
in particular, of Pseudo-Riemannian and symplectic 
affine flat manifolds, which we describe in detail. 
Similarly, the holonomy groups of flat homogeneous manifolds
are determined by the centralisers of prehomogeneous
representations. They turn out to be a nilpotent, if
the homogeneous manifold is complete.

The following section \ref{sect:HAM} is devoted to the basic properties
of homogeneous affine manifolds. Their automorphism groups
develop to prehomogeneous subgroups of the affine group.
We show that the Markus 
conjecture is satisfied for compact homogeneous affine manifolds, and
we describe the  structure of \emph{compact} homogeneous 
flat Pseudo-Riemannian,  and also of \emph{compact} symplectically 
homogeneous affine flat manifolds.

We also introduce a result which states that a homogeneous 
affine manifold of a nilpotent group is complete if and
only if the group is volume preserving. In particular,  a
homogeneous flat Pseudo-Riemannian or symplectic 
affine flat manifold of a nilpotent group is always complete.
The proof of these results depends on methods which are
developed in section \ref{sect:crittrans}. 

In section \ref{sect:FALGs}, we review the relationship 
between the geometry of a flat affine Lie group and the behaviour of its
left and right Haar measures. The main result shows that the 
completeness of a flat affine Lie group is determined by
the interaction of its unimodular character with the volume
character which is defined by the affine structure. 
These results build on the study of \'etale affine 
representations of Lie groups.  

In section \ref{sect:AHDs},  we discuss the basic properties of affinely homogeneous
domains and of prehomogeneous affine representations of Lie groups. 

In the following section \ref{sect:crittrans}, we develop a criterion 
for the transitivity of prehomogeneous affine representations, 
which extends corresponding results for \'etale affine representations. 
The main application shows that every volume preserving 
nilpotent prehomogeneous group of affine transformations 
is transitive. 

In section \ref{sect:Ccc}, we explain how the geometry of invariant measures 
on an affine  homogeneous space and the transitivity properties 
of its associated prehomogeneous representation are linked 
by certain naturally defined \emph{characteristic classes} of 
the affine representation. 

In  section \ref{sect:parvol}, we study  properties of the Zariski closure 
$A(\Gamma)$, where $\Gamma$ is the holonomy group of a \emph{compact} affine 
manifold $M$.   One of the main objectives of the subject is to characterise the groups $
A(\Gamma)$ in relationship with
the geometric properties of $M$.
An important  result of Goldman and
Hirsch states that $A(\Gamma)$ acts transitively on affine
space if $M$ is a volume preserving compact flat affine manifold.
This  has strong consequences for the geometry of
compact volume preserving flat affine manifolds. 
We explain some of the applications
of this result,  and we also explain how its proof relates to the methods 
developed in the previous sections. 

%

\subsubsection*{Acknowledgement} 
I thank Joseph A.\ Wolf, Vicente Cort\'es and Wolfgang Globke for their 
detailed comments, numerous helpful suggestions, and careful reading 
of a first draft of this article. 

\section{Foundations} \label{sect:foundations}

Flat Pseudo-Riemannian and affine manifolds constitute actually 
a particular class of locally homogenous manifolds. This 
point of view allows to express many geometric properties
of flat manifolds in an elegant and transparent way.  

\subsection{The development map and holonomy}

We start by briefly recalling the fundamental notion of 
an $(X,G)$-manifold. For further and more detailed 
reference on $(X,G)$-structures,  see, for example, 
\cite{Epstein,Ratcliffe,Thurston}.

\subsection{$(X,G)$-manifolds}

Let $X$ be a homogeneous space for the
Lie group $G$. 
A manifold $M$ is said to be 
\emph{locally modelled on $(X,G)$}
if $M$ admits an atlas of charts with range in $X$ such that the
coordinate changes are restrictions of elements of $G$. 
A maximal atlas with this property
is then called a {\em $(X,G)$-structure on $M$}, and $M$ is called a
{\em $(X,G)$-manifold}, or \emph{locally homogeneous space} modelled on 
$(X,G)$. 

A map $\Phi$ between two $(X,G)$-manifolds is called an
\emph{$(X,G)$-map} if it looks like the action of an element of $G$ in the local charts. If the $(X,G)$-map
$\Phi$ is a diffeomorphism it is called an {\em $(X,G)$-equivalence}. 

Every $(X,G)$-manifold comes equipped with some extra structure, called the development
and the holonomy. 
Let $\pi:\tilde{M} \rightarrow M$ denote the universal covering space of 
the $(X,G)$-manifold $M$. We fix $x_0 \in M$. 
The 
{\em development map\/} of the $(X,G)$-structure on $M$ is the 
local diffeomorphism 
$$ D:  \tilde{M}  \rightarrow   X$$ which is obtained by 
analytic continuation of a local $(X,G)$-chart of $M$ in $x_0$. 
The development map is an $(X,G)$-map, and, \emph{for 
any $(X,G)$-equivalence $\Phi$ of $\tilde{M}$,
there exists an element $h(\Phi) \in G$ such that }
\begin{equation}  D \circ \Phi= h(\Phi) \circ D \; .  \label{eq:devel1}
\end{equation} 

The fundamental group $\pi_1(M) = \pi_1(M,x_0)$ acts 
on  $\tilde{M}$ via covering transformations, which are
$(X,G)$-equivalences of $\tilde{M}$. 
This induces the {\em holonomy homomorphism\/} 
$h: \pi_1(M,x_0) \rightarrow G$ which satisfies 
\begin{equation}  \label{eq:hol} D \circ \gamma =   h(\gamma)  \circ D  \; ,   
\end{equation}  
for all $\gamma \in \pi_1(M,x_0)$.  After the choice of the development map (which corresponds to a choice of a germ of a $(X,G)$-chart in $x_0$),  
the holonomy homomorphism $h$ is well defined.  
We note that the $(X,G)$-structure on $M$ determines
the conjugacy class of $h$ under the action of $G$. 

Clearly, the development map already determines the $(X,G)$-structure
on $\tilde{M}$, and specifying a development pair for the action of $\pi_1(M,x_0)$ on
$\tilde{M}$ is equivalent to constructing an $(X,G)$-structure on $M$:

\begin{proposition}
Every local diffeomorphism $ D:  \tilde{M}  \rightarrow  X$ which satisfies \eqref{eq:hol},
for some $h: \pi_1(M,x_0) \rightarrow G$ defines a unique $(X,G)$-structure on $M$,
and every $(X,G)$-structure on $M$ arises in this way.
\end{proposition}

\paragraph{Properly discontinuous actions}
Let $\Gamma$ be a group of diffeomorphisms
of a manifold $M$. Then $\Gamma$ is said to act \emph{properly 
discontinuously} on $M$ if, for all compact subsets $\cK \subset M$,
the set $$ \Gamma_{\cK} = \{ \gamma \in \Gamma \mid \gamma \cK \cap \cK \neq
\emptyset \} $$ is finite. 
If $\Gamma$ acts properly 
discontinuously and freely on $M$ then 
the quotient space $M/\, \Gamma$ 
is a smooth manifold, and the projection
map $\pi: M \ra M/\, \Gamma$ is a smooth
covering map. Now if $\Gamma$ acts by
$(X,G)$-equivalences on an $(X,G)$-manifold
$M$, then the quotient space $M/\, \Gamma$ inherits 
a natural  $(X,G)$-manifold structure from $M$.
In fact, the development maps of $M$ and $M/\, \Gamma$
coincide (as well, as their universal coverings.)

\begin{example}[$(X,G)$-space forms] 
\rm 
Assume that $X$ is simply connected, and $\Gamma$
is a group of $(X,G)$-equivalences of $X$ (that is, $\Gamma$ is a subgroup 
of $G$) acting properly discontinuously and freely on $X$. Then 
$X/\, \Gamma$ is an $(X,G)$-manifold, and the identity map
of $X$ is a development map for $X/\, \Gamma$. 
\end{example}

See \cite{Epstein,Thurston} for further discussion of $(X,G)$-geometries
and the properties of the development process. 

\subsubsection{Affine and projectively flat manifolds}
The geometry of the Pseudo-Euclidean space  
$\bbE^s = (\bbR^n, E(s))$  is determined by the transitive action of the 
isometry group $E(s)$ of the standard scalar product of signature $s$.
Similarly the geometry of affine space $\bbA^n = (\bbR^n, \Aff(n))$ is determined 
by the action of the full affine group $\Aff(n)$.

\paragraph{Coordinate representation of the affine group}

We view affine space $\bbA^n$ as a hyperplane 
$$ \bbA^n = \{( x,1) \mid x \in \bbR^n \}$$  embedded in $\bbR^{n+1}$. 
In this setting, the group of affine transformations
$\Aff(n)$ identifies naturally with a group of linear transformations
of $\bbR^{n+1}$. Namely, 
$$  \Aff(n)   =   \left\{ A =  
\begin{matrix}{cc} g & t \\ 0 & 1 \end{matrix} \Mid g \in \GL(n,\bbR)\right\}  \; , $$ 
which is a subgroup of $\GL(n+1, \bbR)$. 
Note that the affine 
group decomposes as a semi-direct product 
$$ \Aff(n) = T(n) \rtimes GL_{n}(\bbR) \; , $$
where  $$ T(n)   =   \left\{ A =  
\begin{matrix}{cc} 1 & t \\ 0 & 1 \end{matrix} \Mid t \in \bbR^n \right\}  \;  $$
is the group of translations of $\bbR^n$. 
The natural quotient 
homomorphism $$ \ell: \Aff(n) \ra GL_{n}(\bbR) \; , \; 
A \mapsto \ell(A)= g   \;  $$ associates to the affine transformation $A \in \Aff(n)$ its \emph{linear part}.
The vector $t(A)= t$, is called the  \emph{translational part} of $A$.

The Lie algebra $\lie{aff}(n)$ of $\Aff(n)$ is 
$$ \lie{aff}(n) = \left\{  X=   \begin{matrix}{cc} \varphi  & v \\ 0 & 0 \end{matrix} 
\Mid \varphi \in \lie{gl}(n,\bbR)  \right\}  \; , $$ 
which is a Lie subalgebra of the matrix algebra $\lie{gl}(n+1,\bbR)$. 

Note that the \emph{evaluation map} of the affine action at $x \in \bbA^n$ 
$$ o_{x}: \Aff(n) \ra  \bbA^n \; , \;  \; A \mapsto A \cdot x \, = \,g(x) +t$$ is expressed by matrix multiplication,
and so is its derivative at the identity, which is  the map 
\begin{equation} \label{eq:tx} 
t_{x}: \lie{aff}(n) \ra \bbR^n  \; , \;  \; X \mapsto X \cdot x \, =  \, 
\varphi(X) +v \; . 
\end{equation}

\paragraph{Transitive subgroups of $\Aff(n)$}
We shall also consider various subgroups of $\Aff(n)$. 
We let $O(s) = O( < , >_{s})$ denote the group of linear isometries of the standard scalar product
 $< , >_{s}$ of signature 
$s$. The Pseudo-Euclidean isometry group $\E(s)$ is a semidirect product 
$\E(s)= T(n) \rtimes O(s)$.  Thus 
$\E(s)$ embeds into the affine group as 
$$  \E(s)  =   \left\{ A =  \begin{matrix}{cc} g & t \\ 0 & 1 \end{matrix} \Mid g \in O(s)\right\}  \; . $$
The group of volume preserving affine
transformations is $\SAff(n) = \{ A \in \Aff(n) \mid \ell(A) \in  \SL_{n}(\bbR) \}$.
The group of symplectic affine transformations is 
$\Aff(\omega_{n}) = \{ A \in \Aff(2n) \mid \ell(A) \in \SP(\omega_{n}) \} $, where 
$\omega_{n}$ is a non-degenerate skew bilinear form on $\bbR^{2n}$. 


\paragraph{Fixed points for reductive affine actions} 
We call a subgroup $G$ of  a linear group \emph{reductive}
if every $G$-invariant linear subspace admits an invariant 
complementary subspace. The following Lemma is a
 particular useful fact:

\begin{lemma} \label{lemma:fixedpoint}
Let $G \leq \Aff(n)$ be a reductive subgroup. Then the affine
action of $G$ has a fixed point on $\bbA^n$.  
\end{lemma}
\begin{proof} Take a complementary line to the $G$-invariant subspace 
$\{ (v,0) \mid v \in \bbR^n \} \subset \bbR^{n+1}$, and intersect with 
the hyperplane $\bbA^n$. 
\end{proof}

\paragraph{Projective geometry} 
The space of lines through the origin in $\bbR^{n+1}$ is called
real projective space, and it is denoted by $\mathbb P^n \, \bbR$. 
The projective linear group is $$ \PGL(n, \bbR) = \GL(n+1 , \bbR)/ \{ \pm E_{n+1} \}$$
acting on the projective space $\mathbb P^n \, \bbR$.  
Via the above coordinate representation, affine space $\bbA^n$ embeds as an 
open subset in $\mathbb P^n \, \bbR$, and $\Aff(n)$ embeds as a subgroup 
of $\PGL(n, \bbR)$. In particular, affine geometry (that is, geometry modelled on
$(X,G)= (\bbA^n, \Aff(n))$), is a subclass of projective $(\mathbb P^n \, \bbR, \PGL(n, \bbR))$
geometry. Note also that a cone over any $(\mathbb P^n \, \bbR, \PGL(n, \bbR))$-manifold
becomes a $(\bbR^{n+1},  \GL(n+1 , \bbR))$ manifold. (See \cite{Goldman1},
for an exposition about projectively flat manifolds, and their relation with 
affine flat manifolds.)
%



%

\subsection{Flat manifolds are $(X,G)$-manifolds}

As already remarked, the geometry of the Pseudo-Euclidean space  $\bbE^s$ is determined by the transitive action of its isometry group $E(s)$, and similarly the geometry of affine space $\bbA^n$ is determined by the action of the affine group $\Aff(n)$. 

\begin{theorem} A flat Pseudo-Riemannian manifold $(M,\g)$ is a
locally homogeneous space modelled on the standard pseudo-Euclidean 
space $\bbE^s$, $s = s (\g)$.  A flat affine manifold $(M, \nabla)$ is a
locally homogeneous space modelled on affine space $\bbA^n$. 
\end{theorem}
\begin{proof} In fact, by flatness of $(M,\g)$, the exponential map 
$\exp_{p}: T_{p} M \ra M$ allows to define an isometry of a neighbourhood of $0$ in 
$\bbE^s$ to a normal neighbourhood in $(M,\g)$, for every point $p \in M$.  This defines a 
\emph{compatible atlas} for $M$, where all charts are local isometries, and with coordinate changes in the Pseudo-Euclidean group $E(s)$.  Evidently, 
the analoguous argument
also works for a flat affine manifold $(M, \nabla)$.
\end{proof}

Note that the locally homogeneous $\bbE^s$-structure on $(M,\g)$ is uniquely determined by the \emph{compatibility condition} that its charts are local isometries. 
The \emph{compatible $\bbA^n$-structure} for $\nabla$ is determined by
the condition that its charts are affine maps. In particular, we may speak of
the development map $D: \tilde{M} \ra \bbA^n$ and holonomy of a flat affine 
manifold, where $D$ is then a local locally affine diffeomorphism 
(respectively, $D: \tilde{M} \ra \bbE^s$ is a local isometry, 
for a flat Pseudo-Riemannian manifold). 
 
The traditional point of view (Pseudo-Riemannian metric and
flat connection) and the $\bbE^s$-manifold point of view are
completely equivalent. Namely, every locally homogeneous space modeled        
on $\bbE^s$ has a unique flat Pseudo-Riemannian metric $\g$
which turns the local charts of the $\bbE^s$-structure into local isometries for $\g$.
A map between two flat manifolds $(M,\g)$ and $(M',\g')$
is a local isometry if and only if it is a map of $\bbE^s$-structures, and so on.

\paragraph{Holonomy and parallel transport}
Let $h: \pi_{1}(M) \ra \E(s)$ denote the holonomy homomorphism
of an $\bbE^s$-structure on $M$. 

\begin{definition} The homomorphism $h$ is called the 
\emph{affine holonomy homomorphism} of the flat manifold $(M,\g)$.
The composition $hol= \ell \circ  h: \pi_{1}(M) \ra O(s)$ is called the
\emph{linear holonomy} homomorphism.
\end{definition}

Recall that the \emph{parallel transport} of a flat manifold $(M,\g)$
defines a homomorphism $p_{x}: \pi_{1}(M,x) \ra O(T_{x} M, \g_{x})$
of the fundamental group of $M$ into the group of linear
isometries of the tangent space. The image $p_{x}(\pi_{1}(M,x))$
is called the holonomy group of $(M,\g)$ (at $x$). 
The following is easy to see:

\begin{proposition} In a local chart for the induced $\bbE^s$-structure
on $M$, based at $x \in M$,  the parallel transport homomorphism 
$p_{x}: \pi_{1}(M,x) \ra O(T_{x} M, \g_{x})$
corresponds to (that is, it is conjugate to) 
the linear holonomy homomorphism of the  $\bbE^s$-structure.
\end{proposition}

Of course, the analogous result holds for a flat affine manifold 
$(M,\nabla)$ and its parallel transport. Note that also the 
affine holonomy may be interpreted as  a parallel transport in 
a suitable associated bundle over $(M,\nabla)$ (cf.\ \cite{KN}). 
For yet another interpretation of the affine parallel transport,
see \cite{BuKa}.

\paragraph{Affine structures of type $\cA$}


We do not need to restrict our attention to the transitive groups 
$E(s)$ or $\Aff(n)$. More generally, we may consider any 
subgroup $\cA$ of the affine group $\Aff(n)$  to
define a flat model geometry. This gives rise to the following notion. 

\begin{definition} A locally homogeneous space modelled on 
 $(\bbA^n, \cA)$ is called a flat affine
manifold of type $\cA$. 
\end{definition}

The holonomy determines if the affine structure group 
for a flat affine manifold $(M,\nabla)$ can be reduced
to a subgroup $\cA$ of $\Aff(n)$.

\begin{definition}
A flat affine manifold $(M,\nabla)$ is called 
of type $\cA$ if (a conjugate) of its affine holonomy homomorphism $
h: \pi_{1}(M) \ra \Aff(n)$ takes image in $\cA$. 
\end{definition}
If $(M,\nabla)$ is of type $\cA$ then it admits a compatible
structure of an $(\bbA^n, \cA)$ manifold. 
The definition allows to consider various kinds of
geometric flatness conditions for $(M, \nabla)$. For example, 
the groups $\cA = S\Aff(n)$, 
$\Aff(\omega)$, $\Aff(2n,\bbC)$, $UE(2n)$
decribe the concepts of volume preserving, symplectic, complex 
or K\"ahler flat affine manifold, respectively.  

\subsubsection{Completeness and crystallographic groups}

Geodesic completeness of a flat affine manilfold is interpreted
naturally by its development map.

\begin{theorem}  A flat affine manifold $(M, \nabla)$ is complete 
if and only if its development map $D: \tilde{M} \ra \bbA^n$ is 
a diffeomorphism. 
\end{theorem}
\begin{proof} The development map $D$ is a local affine diffeomorphism.
By completeness of $(\tilde{M}, \nabla)$ it is onto $\bbA^n$. Moreover,
geodesics of $\bbA^n$ admit a lift along any preimage for $D$.
Therefore, $D$ is a covering  map. Hence, it must be a diffeomorphism. 
\end{proof}

\begin{question} Does there exists a flat affine manifold with development
map $D$, which is onto $\bbA^n$, but not a diffeomorphism? 
\end{question}

The characterisation of completeness via the development implies the following: 

\begin{corollary}[Killing-Hopf theorem] \label{cor:Killing-Hopf}
Let $(M, \nabla)$ be a 
complete flat affine manifold of type $\cA$.  
Let $\Gamma = h(\pi_{1}(M)) \leq \cA$ be its affine holonomy group. Then $\Gamma$ acts properly discontinuously 
and freely on $\bbA^n$, and  $M$ is $\cA$-equivalent to the $\cA$-manifold
$\bbA^n/ \, {\mbox{$\Gamma$}}$. 
 \end{corollary}

It follows that every complete Pseudo-Riemannian manifold 
$(M,\g)$ of signature $s$, is isometric to a quotient of $\bbE^s$ 
by a properly discontinuous group $\Gamma \leq \E(s)$. 
Such a manifold $$ M= \bbE^s / \, \Gamma$$ is called 
a \emph{Pseudo-Euclidean space form}. 

A complete flat affine manifold $(M,\nabla)$ is called an  \emph{affine space form}. By Corollary  \ref{cor:Killing-Hopf}, the study of affine space forms
reduces to the study of properly discontinuous subgroups
of $\Aff(n)$.

\begin{example} Every discrete subgroup $\Gamma$ of the Euclidean group 
$\E(n)$ acts properly discontinuously on $\bbR^n$.  If $\Gamma$ is in
addition torsion-free then $\bbE^n/\Gamma$ is a complete flat Riemannian
manifold. 
\end{example}

In general, if $\Gamma \leq  \cA$ is properly discontinuous then it
is discrete, but the converse may not hold.

\paragraph{Crystallographic groups}
A uniform discrete subgroup of $\E(n)$ is traditionally
called a (Euclidean-) crystallographic group, cf.\ \cite{Wolf}. 
This motivates the following

\begin{definition}
A properly discontinuous subgroup $\Gamma  \leq \cA$, $\cA  \leq \Aff(n)$, is called 
an \emph{affine crystallographic group of type $\cA$} if the quotient space
$\bbA^n/\Gamma$ is compact.  
\end{definition}

\begin{question} 
Is an affine crystallographic group (up to finite index)
contained in the group of volume preserving affine transformations
of $\bbA^n$? Or equivalently, has an orientable compact complete affine manifold
always a parallel volume form? (This is one direction of Markus'
conjecture, see section \ref{sect:Markus}). 
\end{question}

Note that a \emph{finite volume} flat complete Riemannian manifold
is necessarily compact. This is a consequence of the classification
of  discrete subgroups of the Euclidean group $\E(n)$, as given
by Bieberbach, see \cite{Wolf}. 

\begin{question}  Does there exist a non-compact, finite volume
complete affine or Pseudo-Riemannian manifold? For example,
does a complete  flat Lorentzian manifold admit finite volume,
without being compact? 
\end{question}

\subsubsection{The development image of a compact affine manifold}
Let $M$ be a compact affinely flat manifold, $D: \tilde{M} \ra \bbA^n$ its
development map. The development image of $M$ is an open domain 
$$ U_{M}= D(\tilde{M}) \subseteq \bbA^n  \; . $$ 
It is conjectured that $M$ is complete,  if $M$ admits a parallel volume
(Markus' conjecture). In particular, in this case, the conjecture claims that
$D(\tilde{M}) = \bbA^n$.

\begin{example} The development image of an affine two-torus $U_{T^2}$
is one of the four affinely homogeneous domains which admit an abelian 
simply transitive group of affine transformations, see Corollary \ref{cor:dev2t}.
See also  \cite{BauesG} for construction methods of affine structures on $T^2$,
and visualisation of the development process. 
\end{example}

In general, the development image of a compact affinely flat, or projectively
flat, manifold can have a very complicated development image. See \cite{SuThu}, 
for the construction of examples. 

\section{The group of affinities} \label{sect:AffM}

Let  $\nabla$ be  a torsion-free 
connection on $M$. We put $\Aut(M,\nabla)$ for the group of
connection preserving diffeomorphisms of $M$. 
Recall 
(cf.\  \cite{Kobayashi}) 
that  $\Aut(M,\nabla)$ is a Lie group.

\subsection{Affine vector fields}
A vector field $X$ on $(M,\nabla)$ is called an \emph{affine vector field} if its local
flow preserves the connection $\nabla$. 
The affine vector fields 
form a subalgebra $\lie{aut}(M,\nabla)$ of the Lie algebra
$\Vect(M)$ of all $C^\infty$-vector fields. 

A vector field on $M$ is called \emph{complete} 
if its local flow generates a one-parameter group of diffeomorphisms
of $M$. The \emph{complete} affine vector fields on $M$ form
a subalgebra $\lie{aut}_{c}(M,\nabla)$ of $\lie{aut}(M,\nabla)$.  
If $M$ is compact or $\nabla$ is complete
then  $\lie{aut}_{c}(M,\nabla) = \lie{aut}(M,\nabla)$ (see \cite[VI, \S 2]{KN}). 

Note that the Lie algebra of complete affine vector fields is  
(anti-) isomorphic to the tangent Lie algebra of (left-invariant)
vector fields on $\Aut(M,\nabla)$. 

%

\subsubsection{The associative algebra of affine vector fields on a flat manifold}

Let $(M,\nabla)$ be a flat affine manifold, and let $X$ be 
a vector field on $M$. 
The covariant derivative of $X$  defines an endomorphism field 
$$ A_{X} : \;  Y \, \mapsto  \;  - \nabla_{Y} X$$
on $M$. Since $\nabla$ is flat, the  vector field 
$X$ is affine if and only if $A_{X}$ is parallel, that is, $X$ is affine 
if and only if,   
for all vector fields $Y, Z$,
\begin{equation} \label{eq:AXpar}
\nabla_{Z}  A_{X} Y =  A_{X} \nabla_{Z} Y
 \; . \end{equation}
Moreover, for an affine vector field  $X$ on the flat manifold $(M,\nabla)$ 
we have, for all $Y \in \Vect(M)$, 
the  relation
 \begin{equation}
  A_{[X,Y]}  \, =  \, [A_{X}, A_{Y} ]  \;    .
\end{equation}

\paragraph{Induced algebra structure on $\Vect(M)$}
We declare the product of two vector fields $X,Y \in \Vect(M)$ 
by \begin{equation} \label{eq:rsap} X *_{\nabla} Y := \,  A_{X} Y =  - \nabla_{Y} X \; \; . 
\end{equation}
%
In terms of the product $*_{\nabla} $, condition  \eqref{eq:AXpar}, that 
$X \in \Vect(M)$ is affine, is equivalent to  \begin{equation}  
 X  *_{\nabla}  \, (Y  *_{\nabla}  Z) \,  = \, 
( X   *_{\nabla} Y)  *_{\nabla}  Z  \; \; ,
\;     \label{eq:akernel} 
\end{equation}  
for all $Y, Z \in  \Vect(M)$. \\ 

A short calculation, 
involving the associativity condition  \eqref{eq:akernel}
shows: if $X$ and $X'$ are affine vector fields 
then $A_{X} X' =  X *_{\nabla} X'$  also satisfies \eqref{eq:akernel},
for all $Y,Z \in \Vect(M)$. 
Hence, $ X *_{\nabla} X'$ is an affine vector field.
In particular,  \emph{the  Lie algebra 
of affine vector fields $ \lie{aut}(M,\nabla)$ is 
a subalgebra of  
$(\Vect(M),  *_{\nabla})$ which 
is associative:}


\begin{proposition} \label{prop:asu} 
The Lie algebra of affine  vector fields $(\lie{aut}(M,\nabla)$ forms 
an associative subalgebra of $(\Vect(M),  *_{\nabla})$.
\end{proposition}

Let us  furthermore remark that the centralisers of affine vector fields form a subalgebra in 
$(\Vect(M), *_{\nabla})$:

\begin{lemma}  \label{lemma:centraliserY}
Let $Y  \in \aut(M,\nabla)$ be an affine vector field. Let  
$X, X'  \in \Vect(M)$ be vector fields on $M$,
which centralise $Y$. Then also $\nabla_{X} X'$ 
centralises $Y$. 
\end{lemma}
\begin{proof} Using \eqref{eq:AXpar} for the affine field $Y$, 
as well as, $[X,Y] =  [X',Y] = 0$, we calculate 
\begin{eqnarray*}  [\nabla_{X} X', Y]  & =  &\nabla_{\nabla_{X} X'}  Y - \nabla_{Y} \nabla_{X} X' \\
&= &   \nabla_{X} \nabla_{X'}  Y -   \nabla_{Y} \nabla_{X} X' \\  &= &
 \nabla_{X} \nabla_{Y} X' -   \nabla_{Y} \nabla_{X} X' \\   &  = & \nabla_{[X,Y]} X'  \, = 0  \; \, . 
\end{eqnarray*} \end{proof}

\begin{remark} 
The fact that  $ \lie{aut}(M,\nabla)$ inherits the structure of
an associative matrix algebra, has strong consequences, 
especially  in a situation, where $ \lie{aut}(M,\nabla)$ describes
the tangent algebra for  the group $\Aut(M,\nabla)$.  
For example,  \emph{homogeneous} affine manifolds, which
are compact or complete are naturally related to associative 
matrix algebras.  See, for example,  \cite{Yagi} for applications. 
See also some of the results in section \ref{sect:AHDs} and section \ref{sect:crittrans}
for further exploitation of this principle. 
\end{remark}

\subsubsection{Right- and left-symmetric algebras} \label{sect:LRalgs}
A bilinear product $\cdot$ on a vector space which satisfies the identity
\begin{equation}
  X \cdot (Y \cdot Z)  -  ( X \cdot Y) \cdot Z \,   =  \;   Y \cdot (X \cdot Z)  -  ( Y \cdot X) \cdot Z   \; \;   \label{eq:LSA2}
\end{equation}  is called a \emph{left-symmetric algebra}. 
Complementary, a product is called \emph{right-symmetric} if 
 \begin{equation}
Z \cdot (Y \cdot X)  -  ( Z \cdot Y) \cdot X \,   =  \;   Z \cdot (X \cdot Y)  -  ( Z \cdot X) \cdot Y   \; \;   \label{eq:RSA2} \; . 
\end{equation}
Both identities  naturally generalise the associativity condition for algebras, and,
by exchanging factors, right- and left symmetric products are in one to one 
correspondence. Moreover, as in the associative case, there is an 
associated Lie product which is defined by 
\begin{equation*}
  [ X, Y]  \,: = \; X \cdot Y - Y \cdot X  \label{eq:LSA1} \; . 
\end{equation*}
The associative subalgebra of elements $X$ which is defined by equation \eqref{eq:akernel} is called the \emph{associative kernel}   
of the right-symmetric algebra. \\

\begin{example} Let $(M,\nabla)$ be a flat affine manifold. 
As declared
in \eqref{eq:rsap}, the flat torsion-free connection $\nabla$
induces an algebra product $*_{\nabla}$ on the
Lie algebra $\Vect(M)$ of $C^\infty$- vector fields on $M$.
Note then that $*_{\nabla}$  is a compatible \emph{right-symmetric 
algebra} product 
on  the Lie algebra $\Vect(M)$.  The associative 
kernel of $(\Vect(M),  *_{\nabla})$ is finite dimensional, and it is precisely the
subalgebra  $(\lie{aut}(M,\nabla), *_{\nabla})$ of affine vector fields. 
\end{example} 

Algebras satisfying  \eqref{eq:LSA2} appeared and a play mayor role in the study of convex homogeneous cones \cite{Vinberg}. Naturally finite dimensional LSAs over the real numbers play an important role in the study of left invariant flat affine structures on Lie groups, see section \ref{sect:FALGs}. 
Left symmetric algebras also appear in several other mathematical and physical contexts, see \cite{Burde_2} for a recent survey. 

\subsection{Development representation of affinities} \label{sect:holrepA}
Let $(\tilde{M},\nabla)$  be the universal
covering manifold of $(M, \nabla)$. We put $\tilde \Gamma \leq   \Aut(\tilde{M},\nabla)$
for the the group of covering transformations, and 
$\Gamma = h(\tilde \Gamma)$ for  the affine holonomy group of $(M,\nabla)$. 
The development process provides  a local 
representation of the group of affine transformations $\Aut(M,\nabla)$:
\\

We define 
$$ \hat \Aut(M, \nabla) = \Aut(\tilde{M}, \tilde \Gamma, \nabla)  = \,  \{ \Phi \in \Aut(\tilde{M}, \nabla) \mid
 \Phi \,   \tilde \Gamma \,  \Phi^{-1} =  \tilde \Gamma \} \; . $$ 
The group $\hat \Aut(M, \nabla)$ is a covering group 
of $\Aut(M,\nabla)$, and  the
development  homomorphism \eqref{eq:devel1} induces a
homomorphism of Lie groups
\begin{equation} \label{eq:devautM} 
h: \hat \Aut(M, \nabla) \ra \Aff(n) 
\end{equation} 
into the affine group $\Aff(n)$. Note that $h$ has discrete kernel. Moreover, 
the image $h( \hat \Aut(M, \nabla))$ normalises $\Gamma$, and  
$h( \hat \Aut(M, \nabla)^0)$ centralises $\Gamma$. 

\paragraph{Representation of affine vector fields}

We study now the tangent representation of the development 
homomorphism. 
\begin{proposition} \label{prop:linrep} 
The development defines a 
a natural  faithful associative algebra representation 
\begin{equation} \label{eq:linrep} 
 \bar h: \, 
(\lie{aut}(M,\nabla),  *_{\nabla})  \, \ra \, \lie{aff}(n)  \; .  
\end{equation}
\end{proposition}
\begin{proof} 
We already remarked that $ \lie{aut}(M,\nabla)$ forms an associative algebra
with respect to $ *_{\nabla})$.
Via local affine coordinates, we can identify
the tangent space $T_{p} M$ with $\bbR^n$. By the formula
$ A_{ X *_{\nabla} Y}= A_{ A_{X} Y} =   A_{X} A_{Y } $, which is deduced from 
\eqref{eq:akernel}, the map 
\begin{equation} \label{eq:linrep1} 
\bar{h}:  X \mapsto \;  \phi_{{X}_{p}} =  \begin{matrix}{cc} A_{X_{p}} & X_{p}  \\ 0 & 0
\end{matrix} 
\end{equation}
is easily seen to be a faithful representation of
$(\lie{aut}(M,\nabla), *_{\nabla})$ into the associative algebra of 
$\lie{aff}(n)$ (where the associative algebra 
structure on $\lie{aff}(n)$ is given by the usual product 
of matrices). 
\end{proof}

The development representation  $$ \bar h: \, 
\lie{aut}(M,\nabla)  \, \ra \, \lie{aff}(n)  \;  , $$ 
as defined in \eqref{eq:linrep1},
is associated to the local representation $h$ of $\Aut(M,\nabla)$. 
In fact, on the subalgebra  $\lie{aut}_{c}(M,\nabla)$, the development 
representation $\bar h$ corresponds to the 
derivative of $h$. We further remark: 

\begin{proposition} \label{prop:linrep2} 
The representation $\bar h$
identifies $\lie{aut}(M,\nabla)$ with the subalgebra  
$\lie{aff}(n)^\Gamma$ of\/ $\Gamma$-invariant affine vector fields 
on $\bbA^n$. 
\end{proposition}
\begin{proof} 
If $M$ is simply connected the map $\bar{h}: \lie{aut}(M,\nabla)  \, \ra \, \lie{aff}(n)$
is an isomorphism, since affine vector fields may be extended uniquely to all of $M$ 
from any coordinate patch. For the general case,  note that $\lie{aut}(M,\nabla)$
is isomorphic to the subalgebra $\lie{aut}(\tilde M,\nabla)^{\tilde \Gamma}$ 
of $\Gamma$-invariant affine vector fields on $\tilde M$.
Since $\bar{h}: \lie{aut}(\tilde M,\nabla)  \, \ra \, \lie{aff}(n)$ is equivariant with respect to $h: \tilde \Gamma \ra \Gamma$, $\lie{aut}(M,\nabla)$ is mapped 
onto $\lie{aff}(n)^\Gamma$.
\end{proof}


\subsection{Affinities of compact volume preserving affine manifolds} 
\label{sect:compvp}

Recall that a compact volume preserving affine manifold is conjectured 
to be complete. Now if  $(M,\nabla)$ is a compact complete affine manifold then
the centraliser of the affine crystallographic group $\Gamma = hol(\pi_{1}(M))$
in $\Aff(n)$ is unipotent (see for example  \cite{GS}, or section \ref{sect:parvol}). 
Since the centraliser of  $hol(\pi_{1}(M))$ is  the development of 
$\hat{\Aut}(M, \nabla)^0$, the latter group is a simply connected nilpotent Lie 
group, which is faithfully represented by unipotent matrices. In this section 
we show that the analogous facts do indeed hold  for the affinities of 
a compact volume preserving affine manifold, without assuming its
completeness.  \\

Let  $(\lie{aut}(M,\nabla), *_{\nabla})$ denote the Lie algebra of 
affine vector fields on $M$ with the associative algebra structure
induced by $\nabla$. We note: 

\begin{proposition} \label{prop:compvp}
Let $(M, \nabla)$ be a flat affine manifold
with parallel volume, which is compact (or with finite volume).
Then $(\lie{aut}(M,\nabla), *_{\nabla})$ is isomorphic
to an associative algebra of nilpotent matrices. In particular,
$\lie{aut}(M,\nabla)$ is a nilpotent Lie algebra. 
\end{proposition}
\begin{proof} Define the divergence of $X \in \Vect(M)$ relative to the
parallel volume $\mu$ on $M$ by the formula
$$ \div X \, \mu =  L_{X} \mu \; . $$ 
Since $M$ is compact, by Green's theorem (cf. \cite[Appendix 6]{KN})
we have, 
\begin{equation} \int_{M} \div X d\mu  =  0 \; , \label{eq:green} 
\end{equation}
and, moreover, since the volume form is
parallel, we also have  (again according to \cite[Appendix 6]{KN}), 
$$  \div X =  - \trace A_{X} \;  (= \trace \nabla_{\bullet} X)  \; . $$
In particular,  if $X \in \lie{aut}(M,\nabla)$
then $\div X$ is a constant function, which implies $\div X =0$.
Thus, via the linear representation $X \mapsto \phi_{X_{p}}$ 
of $(\lie{aut}(M,\nabla), *_{\nabla})$
constructed in Proposition \ref{prop:linrep}, $\lie{aut}(M,\nabla)$ 
is isomorphic to an associative algebra of linear operators of trace zero. 
By Lemma \ref{lemma:nilopsalg}, $\phi_{X_{p}}$ is a nilpotent operator,
for all $X \in  \lie{aut}(M,\nabla)$. 

In the finite volume case, we use an analogous
argument, noting that \eqref{eq:green} holds for an affine 
vector field $X$ on finite volume $M$ as well. 
\end{proof}

\begin{corollary} \label{cor:compvp}
Let $(M, \nabla)$ be a flat affine manifold
with parallel volume which is compact (or with finite volume).
Then the centraliser of the affine holonomy group $\Gamma$ of $M$
is a connected unipotent group, which is isomorphic
to $\hat{\Aut}(M,\nabla)^0$ under the development. 
\end{corollary}
\begin{proof} In appropriately chosen
coordinates, the homomorphism $\bar h$, defined in \eqref{eq:linrep},
corresponds to the differential of the holonomy representation 
$h: \hat{\Aut}(M,\nabla) \ra \Aff(n)$.
By Proposition \ref{prop:linrep},  the tangent Lie algebra of
the holonomy image is the associative algebra $\lie{aff}(n)^\Gamma$, which
is a subalgebra in $\lie{aff}(n)$. By Proposition \ref{prop:compvp},
the elements of $\lie{aff}(n)^\Gamma$ are nilpotent. 
If $A\in Z_{\Aff(n)}(\Gamma)$ is an affine transformation which centralises
$\Gamma$ then $A - E_{n+1} \in \lie{aff}(n)^\Gamma$. 
Thus, $A$ is a unipotent affine transformation. 
It follows that $Z_{\Aff(n)}(\Gamma)$ is
an algebraic subgroup of $\Aff(n)$ consisting of unipotent
elements. In particular, $Z_{\Aff(n)}(\Gamma)$ must be connected,
and simply connected.  It follows that the development homomorphism \eqref{eq:devautM},
$$ h: \hat{\Aut}(M,\nabla)^0 \ra  Z_{\Aff(n)}(\Gamma) $$ is a covering map onto, 
and, in fact, it is an isomorphism. 
\end{proof}

The following is actually an application of Theorem \ref{thm:GHvol}. 

\begin{corollary} \label{cor:compvp2}
Let $(M, \nabla)$ be a flat affine manifold
with parallel volume, which is compact.  
If an affine vector field on $M$ has a zero, it must be trivial. 
In particular, $\dim \Aut(M,\nabla) \leq n= \dim M$.
\end{corollary}
\begin{proof} In fact, by Theorem \ref{thm:GHvol}, the Zariski
closure of $\Gamma$ acts transitively on $\bbA^n$. This implies 
that the unipotent group $Z_{\Aff(n)}(\Gamma)$ acts without
fixed points on $\bbA^n$. 
\end{proof}
In particular, any locally faithful connection preserving  action on a compact volume preserving affine manifold is locally free (the stabiliser of any point is discrete).

\subsection{Lie groups of isometries} 
 The group of isometries $\Isom(M,\g)$ of a Pseudo-Riemannian manifold 
 is a finite-dimensional Lie group, and it is a subgroup of  the group of affine (connection preserving) transformations for the Levi-Civita connection of $\g$.  
 If $M$ is compact and $\g$ is Riemannian then $\Isom(M,\g)$ is
 compact.  The identity components of the automorphism groups of 
 compact Loren\-tzian manifolds are known to be of rather restricted type. 
 (Compare \cite{AS, Zeghib, Zimmer}). 
In fact, by \cite{Zimmer}, non-compact
semisimple factors in the isometry group
of a compact Lorentzian manifold are locally isomorphic to $\SL(2,\bbR)$.
Moreover, a connected solvable subgroup of isometries of a compact Lorentzian manifold
is a product of an abelian group with a 2-step nilpotent group of Heisenberg type 
(see \cite{Zeghib}).

In the presence of a compatible \emph{flat} torsion-free connection,
particular strong restrictions hold for the group of isometries of a compact 
Pseudo-Riemannian manifold $(M, \g)$. 
Proposition \ref{prop:compvp} already
implies that the identity component of $\Isom(M,\g)$  is a nilpotent Lie group.  
Here we show, more specifically, that, in case $(M, \g)$ is a compact flat  Pseudo-Riemannian manifold, $\Isom(M,\g)$ is a \emph{two-step nilpotent group} of 
rather restricted type. If $(M, \g)$ is a compact flat  Lorentzian manifold then
the identity component $\Isom(M,\g)^0$ is abelian and consists of translations only. 

\subsubsection{Killing vector fields} 
Let  $(M,\g)$ be a Pseudo-Riemannian manifold and $\nabla$ the Levi-Civita 
connection for $\g$. A vector field $X$ on $M$ is called a \emph{Killing
 vector field} if its flow preserves the metric $\g$. Equivalently, we
have $L_{X} \, \g = 0$.  Thus, $X$ is Killing if and only if 
 \begin{equation} \label{eq:LXg}
 L_{X}  <U,V> \; =  \; < [ X,U] \, , \, V> +  <U \, , \, [X,V]> \; , 
\end{equation}  for all
 vector fields $U,V$ on $M$.  Since
 $\g$ is parallel, \eqref{eq:LXg} holds if and only if $A_{X}$ is skew with
 respect to $\g$.  Namely,  \eqref{eq:LXg} is equivalent to 
   \begin{equation}  \label{eq:AXskew} 
  <\nabla _{U} X\, ,\, V>  \, + \,  <U \, , \, \nabla _{V} X> \; =  \; 0 \; . 
  \end{equation} 
 
Let $\isom(M, \g)$ denote 
 the Lie algebra of Killing vector fields for $(M,\g)$. Since
Killing vector fields preserve the Levi-Civita connection, every
Killing vector field is affine. 
We have the inclusion 
 $$ \isom(M, \g) \subseteq \aut(M, \nabla) \; . $$

\begin{lemma}  \label{lemma:ocentraliserY}
Let $X,  Y  \in \isom(M, \g )$ be Killing vector fields,
where $Y$ centralises $X$. Then \begin{enumerate}
\item  $   <\nabla _{Y} X,  X >  \; =  \; 0$. 
\item  $  < \nabla_{X} X, Y >  \; =  \; 0$. 
\end{enumerate}
\end{lemma}
\begin{proof}  Since $Y$ commutes with $X$, $\nabla_{X} Y = \nabla_{Y} X$.
Because $Y$ is Killing, $A_{Y}$ is skew. Hence, 
 $$   <\nabla _{Y} X ,  X > \, =  \, <\nabla _{X} Y,  X > \, = -  < X,  \nabla _{X} Y > \; .  $$ 
Thus 1.\  follows. 
We note next that  $  L_{X} \!  < Y, X > \, = 0$, since $[X,Y] =0$.
Since $\g$ is parallel, this implies 
$$ 0 =  \;  < \nabla_{X}Y, X > +   < Y, \nabla_{X }X > \, = \, < \nabla_{Y} X, X > + 
  < Y, \nabla_{X }X > \; . $$
By 1., $< \nabla_{Y} X, X >\,  = 0$, and, hence, 2.\ holds. 
\end{proof}


\begin{proposition} \label{prop:Gcentrvf}
Let $X$ be a Killing vector field whose flow centralises a 
group $G$ of isometries of $(M, \g)$, which has an open orbit on $M$. 
Then \begin{enumerate}
\item  $  \nabla_{X} X  \; =  \; 0$. 
\end{enumerate}
If $(M ,\g)$ is also flat then: \begin{enumerate}
\item[2.]  $A_{X} A_{X} = 0$.
\end{enumerate}
\end{proposition}
\begin{proof} The Killling fields $Y$  corresponding to the 
action of $G$ span the tangent spaces on an open subset of 
$M$. 
Thus,  $  \nabla_{X} X  \; =  \; 0$
is a consequence of Lemma \ref{lemma:ocentraliserY}.

If $(M,\g)$ is flat then $A_{X}$ is parallel. Therefore, 
$$ A_{X} A_{X} Y   =  - A_{X} \nabla_{Y} X =  - \nabla_{Y} \A_{X} X= 
 \nabla_{Y}  \nabla_{X} X = 0 \; . $$
 Hence, $A_{X} A_{X} = 0$. \end{proof} 

We deduce a few consequences: 

\begin{proposition} \label{prop:Killingcent}
Under the assumptions of Proposition \ref{prop:Gcentrvf}, 
the following hold, for all Killing vector fields $X$, $X'$, $X''$, 
which are commuting with $G$:  \begin{enumerate} 
\item $[X', X] =  - 2 A_{X'} X = 2 A_{X} X'$.
\end{enumerate}
If $(M ,\g)$ is also flat then: \begin{enumerate}
\item[2.] $A_{X'} A_{X'' } X = A_{X''} A_{X' } X$.
\item[3.] $[[X', X'' ], X] = 0$. 
\item[4.] $A_{[X,X']} =  [ A_{X},  A_{X'} ] = 2 A_{X} A_{X'} $.
\item[5.]  $A_{X} A_{X'} A_{X''} = 0$.
\end{enumerate}
\end{proposition}

\begin{proof} Note that $X + X'$ is Killing and commutes 
with $G$. Therefore, by 1.\ of Proposition \ref{prop:Gcentrvf} , 
$ \nabla_{X + X'} X + X' =  \nabla _{X'} X + \nabla_{X}  X' = 0$. 
It follows that $[X, X' ] =  \nabla_{X} X' - \nabla_{X'} X = - 2 A_{X} X'$.
Therefore 1.\ holds. 

Next  (using 1.) we note that 
$A_{X'} A_{X'' } X = - A_{X'} A_{X} X''$.
Since  $A_{X'}$ is parallel, $ - A_{X'} A_{X} X''  = - A_{A_{X'} {X}} X''$.
We remark that the Killing vector fields centralising
$G$ form a Lie algebra with respect to the bracket of vector fields. 
Thus, by 1.\,  $A_{X'} X$ is  a Killing vector field, and also
centralises $G$. Therefore, $ A_{X'} A_{X'' } X = A_{X''} A_{X'} {X}$.
Thus, 2.\ holds. 

Now  $[[X', X'' ], X] =  - 2 A_{[X', X'' ]} X  = - 2 ( A_{X'} A_{X'' } X -  A_{X''} A_{X' } X) = 0$
follows.   Thus, 3.\ holds. 

Using polarisation, 2.\ of Proposition \ref{prop:Gcentrvf} 
implies that $A_{X} A_{X'} = - A_{X'} A_{X}$.
Therefore,  $A_{[X,X']} = [A_{X}, A_{X'}] = 2 A_{X} A_{X'} $. 
Using these facts, we can compute,
$ A_{X} A_{X'} A_{X''}  =  (A_{X} A_{X'}) A_{X''}  = -  A_{X''} A_{X} A_{X'} =   A_{X}  A_{X''}  A_{X'} = - A_{X} A_{X'} A_{X''}$. 
\end{proof}

\begin{corollary} \label{cor:Gcentralg}
Let $(M, \g)$ be a flat Pseudo-Riemannian manifold which admits 
a group $G$ of isometries, which has an open orbit on $M$. 
Then \begin{enumerate}
\item  The Lie algebra  $\isom(M,\g)^G$
 of Killing vector fields on $(M, \g)$ which centralise $G$
 forms a subalgebra of the associative algebra
 of affine vector fields  $(\aut(M,\nabla), *_{\nabla})$.  
 \item  The Lie algebra $\isom(M,\g)^G$ is (at most) two step nilpotent. 
\end{enumerate}
\end{corollary}

The result implies that the centraliser of a prehomogeneous group
of isometries on $\bbE^s$ is a unipotent group, see Corollary 
\ref{cor:psrcisunip}.  The main applications concern 
the automorphism groups of compact  flat Pseudo-Riemannian 
manifolds (see section \ref{sect:isomcomp} below) and the holonomy groups of
homogeneous  flat Pseudo-Riemannian manifolds (see section
\ref{sect:homhol}).

\subsubsection{Isometries of compact flat Pseudo-Riemannian manifolds}
\label{sect:isomcomp}

The following result gives a rough description of the possible 
connected groups of isometries of compact flat Pseudo-Riemannian
manifolds. It is another consequence of Theorem \ref{thm:GHvol}:  

\begin{theorem} \label{thm:Killing_vf}
 Let $(M,\g)$ be a compact flat Pseudo-Riemannian 
 manifold, and $\nabla$ the Levi-Civita connection. 
 Let $X \in  \isom(M,\g)$ be a Killing vector field. 
 Then the following hold:
 \begin{enumerate}
 \item If $X$ has a zero then $X=0$. 
\item  $A_{X} A_{X} = 0$. 
\item  The Lie algebra  $\isom(M,\g)$
 of Killing vector fields on $(M, \g)$
 forms a  subalgebra of the associative algebra
 of affine vector fields  $(\aut(M,\nabla), *_{\nabla})$.  
 \item  The Lie algebra $\isom(M,\g)$ is (at most) two step nilpotent. 
 \end{enumerate}
\end{theorem}
\begin{proof} Let $\Gamma$ be the holonomy group of $(M,\g)$. 
By Theorem \ref{thm:GHvol}, the Zariski closure $A(\Gamma)$ of the
holonomy group acts transitively on $\bbA^n$. This transitive group
of isometries commutes with the development of the 
Killing vector fields on $M$. Therefore, $X$ cannot have a zero. 
Moreover, Proposition \ref{prop:Gcentrvf} and Proposition \ref{prop:Killingcent}
imply the next three claims. 
%
\end{proof} 
 
 \begin{remark}
 The development homomorphism  \eqref{eq:linrep},
  $ \bar h: \aut(M,\nabla)  \ra \aff(n)$, maps the subalgebra
  $\isom(M,\g)$ to an associative subalgebra of $\lie{aff}(n)$,
  which is contained in $\isom(\g_{s})$. By the theorem, 
  the elements $\bar h(X) = \phi_{X_{p}}$,
  $X \in \isom(M,\g)$,  satisfy the conditions $X_{p} \notin {\rm Im}\,  A_{X_{p}}$ and $A_{X_{p}}^2  = 0$. Subalgebras of linear maps satisfying 
  both conditions appear in the context  of \emph{complete} homogeneous 
  Pseudo-Riemannian manifolds, as well. (See section \ref{sect:homhol}.)  
  \emph{Abelian} algebras of this type have been further investigated in \cite{Wolf_2,Wolf_3}. 
  For the construction of non-abelian examples, see section \ref{sect:biinvariantmetricLAs}.
 \end{remark}

On a  compact Riemannian manifold with non-positive Ricci curvature,
every Killing vector field must be parallel, see \cite{KN}. 
This holds also in the flat Lorentzian case:  

\begin{corollary} Let $(M,\g)$ be a compact flat Pseudo-Riemannian 
 manifold. If $\g$ is Riemannian or $\g$ has Lorentzian
 signature  then every Killing vector field on $M$ is parallel. 
\end{corollary} 
\begin{proof} Let $X$ be Killing. By Proposition \ref{thm:Killing_vf}, the 
linear operator $A_{X}$ develops to a two-step nilpotent
element contained in $\o(\g)$.  In the Riemannian case 
 $\o(\g)$ has no non-zero nilpotent elements. 
For the Lorentzian case, note that 
every two-step nilpotent element of $\o(\g)$ is zero. 
This implies $A_{X} = 0$, in both cases. 
\end{proof} 

In section \ref{sect:biinvariantmetricLAs}, we construct a
compact flat manifold $(M,\g)$ of dimension six, with signature $s(\g)= (3,3)$,
and six-dimensional non-abelian algebra $\isom(M,\g)$. 

 
\subsection{Lie groups of symplectic transformations}
Let $\omega$ be a non-degenerate $2$-form on $M$
which is closed. 
Then $(M,\omega)$ is called a \emph{symplectic manifold}. 
A diffeomorphism of $M$ which preserves $\omega$
is called a \emph{symplectomorphism} of $(M,\omega)$. 
 In general, the group of symplectomorphisms $\Diff_{\omega}(M)$ of 
 a symplectic manifold is very large, and it is not a Lie group,
 even if $M$ is compact.    
 However, if $M$ is compact, there do exist strong restrictions on the finite-dimensional
 Lie subgroups of $\Diff_{\omega}(M)$. For example, 
 Zwart and Boothby \cite{ZB} proved that a \emph{compact} symplectically
 homogeneous manifold of a solvable Lie group $S$ is diffeomorphic 
 to a torus $T^{2n}$. More generally, by Guan's work \cite{Guan}, 
 for any compact symplectic manifold $(M,\omega)$, every 
 connected solvable Lie subgroup of  $\Diff_{\omega}(M)$  
 is $2$-step solvable, with a compact adjoint image. Moreover, any finite 
 dimensional connected subgroup of $\Diff_{\omega}(M)$  is a semi-direct 
 product of a compact group and a group $S$ of the above type. 
 
 Let $(M,\omega)$ be compact symplectic manifold, which admits
 a compatible flat affine connection $\nabla$.
 We show below (cf.\ Theorem \ref{thm:symplectic_autos})
 that every Lie subgroup of $\Diff_{\omega}(M)$,  which 
 also preserves $\nabla$,  is abelian.  \\
 
 \begin{remark}
 The above mentioned restrictions do \emph{not} apply if $M$ is non-compact. 
 For example, there do exist plenty solvable Lie subgroups of $\Diff(\bbR^n)$, 
 not $2$-step  solvable, which preserve the standard 
 symplectic structure on $\bbR^{2n}$. Such examples may be constructed using 
 solvable Lie groups with left-invariant symplectic structure  (symplectic Lie groups).
 See section \ref{sect:dtc}, for a particular construction method
 for such groups, which also produces non-solvable examples.  
 Further examples are discussed, for example, in \cite{Tralle}. 
\end{remark}
 
 \subsubsection{Symplectic vector fields}
 A vector field $X$ on $(M,\omega)$ is called a \emph{symplectic
 vector field} if its flow preserves $\omega$. Equivalently, $X$
 satisfies $L_{X} \omega = 0$, which means that
 $$ L _{X}  \, \omega(U,V) = \omega([X,U],V) +  \omega(U, [X,V]) \, , $$ for all
 vector fields $U,V$ on $M$. 
  
 \begin{lemma} Let $X,Y,Z$ be symplectic vector fields on $(M,\omega)$.
 Then \begin{equation}   \label{eq:omegaclosed1} 
 \omega(X, [Y,Z]) +   \omega(Z, [X,Y]) +   \omega (Y, [Z,X]) = 0 \; . 
 \end{equation}
 \end{lemma}
 \begin{proof} 
 Since $ d \omega = 0$,  for all vector fields $X,Y,Z$ on $M$, the relation
 \begin{equation}   \label{eq:omegaclosed}
\begin{split} 
L_{X}  \, \omega(Y,Z)  + L_{Z}   \,   \omega(X,Y)  + L_{Y}   \,   \omega(Z,X)  + \;  \\ 
   \omega(X, [Y,Z]) +   \omega(Z, [X,Y]) +  \omega (Y, [Z,X])    \; =  \;0  \; 
\end{split} 
 \end{equation} is satisfied. 
  If  $X,Y,Z$ are symplectic,  $L_{X} \omega = 
  L_{Y} \omega =  L_{Z} \omega = 0$. It follows that
 \begin{multline*}   L_{X}  \, \omega(Y,Z)  + L_{Y}   \,   \omega(Z,X) + L_{Z}   \,   \omega(X,Y) =  \\
 - 2 \left(   \,  \omega(X, [Y,Z]) +   \omega(Z, [X,Y]) +   \omega (Y, [Z,X]) \,   \right) \; .
 \end{multline*} 
 Substituting the right hand side in \eqref{eq:omegaclosed}, the lemma follows. 
 \end{proof} 
 
 \begin{proposition} \label{prop:symplecticentraliser}
 Let $(M, \omega)$ be a symplectic manifold which
 admits a Lie group $G \leq \Diff_{\omega}(M)$, which
 acts transitively on $M$. 
 Then every connected Lie subgroup of $ \Diff_{\omega}(M)$, 
 which centralizes $G$, is abelian.  
 \end{proposition}
 \begin{proof} Let $X,Y$ be vector fields on $M$, whose flows
 preserve $\omega$, and which centralise $G$. Let $Z$ be
 a vector field tangent to $G$.  Since $X,Y,Z$  are symplectic,
 equation \eqref{eq:omegaclosed1} holds. Since $[X,Z] = [Y,Z] = 0$,
 it follows that $$ \omega(Z, [X,Y]) = 0 \; . $$ 
 Since $G$ acts transitively on $M$, the tangent spaces at 
 every point $x \in M$, are spanned by vectors $Z_{x}$, where
 $Z$ is tangent to $G$. Hence, $[X,Y] = 0$. This proves the 
 proposition.
 \end{proof}
 
 The following special case of Proposition \ref{prop:symplecticentraliser}
 is well known:
 
 \begin{corollary} \label{cor:bisymplectic} 
 Let $(G, \omega)$ be a symplectic Lie group, where 
 $\omega$ is biinvariant. Then $G$ is abelian. 
 \end{corollary}
 
\subsubsection*{Symplectic affine vector fields}
 Let $(M,\nabla,\omega)$ be a flat affine manifold with 
 parallel symplectic form $\omega$. As in the metric case, a vector field 
 $X$ is symplectic if and only if $A_{X}$ is skew with respect to $\omega$. 
  An affine vector field $X$ on $(M,\nabla)$ which is also symplectic is 
 called a \emph{symplectic affine vector field}. The symplectic affine vector fields 
 form a Lie subalgebra $\sym(M,\nabla,\omega)$ of $\aut(M,\nabla)$.\\

 The following gives a useful analogue of Lemma \ref{lemma:ocentraliserY} 
 for the symplectic case:
 
 \begin{lemma} \label{lemma:scentraliserY}
 Let $X$ be a symplectic affine vector field on $(M,\nabla, \omega)$,
 and let $Y,Z$ be symplectic affine vector fields which commute with $X$. Then 
  \begin{equation}   \label{eq:omegaclosed2} 
 2\,  \omega(A_{X} A_{X} Y, Z) \,  =  \, \omega(\nabla_{X} X,  [Y,Z ] ) \; . 
 \end{equation}
 \end{lemma}
 \begin{proof} The proof is a straightforward computation. See
  \cite[Proof of Theorem 9]{BC_1} for a special case.
 \end{proof} 
 

\subsubsection{Automorphisms of compact symplectic affine manifolds}
 On a symplectic affine manifold $M$ the group of affine symplectic 
 transformations forms a Lie subgroup of $\Diff_{\omega}(M)$. 
 If $M$ is compact, we show that this group must be abelian. 
 In fact, the following holds: 
  
 \begin{theorem} \label{thm:symplectic_autos}
 Let $(M,\nabla)$ be a compact flat affine manifold with
 parallel symplectic structure $\omega$. Then the identity component of
 the group of symplectic affine transformations is an abelian
 group which develops to a unipotent subgroup of symplectic
 affine transformations of $\bbR^n$.   
 \end{theorem} 
  \begin{proof} 
 Let $\Gamma \leq \Aff(\omega_{n})$ be the holonomy group of $M$. 
 The group of symplectic 
 affine transformations $\hat \Aut(M,\nabla,\omega)^0$ develops
 onto a subgroup $h(\hat \Aut(M,\nabla,\omega)^0)$ of $\Aff(\omega_{n})$
 (compare section \ref{sect:holrepA}). 
 This group centralises $\Gamma$, and therefore also the Zariski closure $A(\Gamma) \leq  \Aff(n)$.
 The Zariski closure $A(\Gamma)$ is a group of symplectic affine
 transformation which acts 
 transitively on affine space $\bbA^{2n}$ (by Theorem \ref{thm:GHvol}). 
 Thus, by Proposition \ref{prop:symplecticentraliser}, 
$h(\hat \Aut(M,\nabla,\omega)^0)$ must be abelian. Hence, so
is $\Aut(M,\nabla,\omega)^0$.
 \end{proof}

 It follows that a compact affine manifold which is homogeneous under the 
 group of symplectic  affine transformations  
 must be diffeomorphic to a torus: 
 
 \begin{corollary}  \label{cor:symplectic_hom}
 Let $(M,\nabla)$ be a compact flat affine manifold with
 parallel symplectic structure which is homogeneous for 
 the group of symplectic affine transformations. Then 
 $M$ is diffeomorphic to a torus. 
 \end{corollary}

Let $\sym(M,\nabla,\omega)$ denote the Lie algebra of symplectic affine vector fields on $M$. As shown above, if $M$ is compact then $\sym(M,\nabla,\omega)$ is abelian. 
We further note: 

\begin{proposition}   \label{prop:symplectic_autos}
Let $(M,\nabla)$ be a compact flat affine manifold with
parallel symplectic structure $\omega$.
Then the following hold:
 \begin{enumerate}
\item The Lie algebra 
 of symplectic affine vector fields $\sym(M,\nabla,\omega)$ is abelian,
 and $\dim \sym(M,\nabla,\omega) \leq  \dim M$. 
 \item For all $X, Y \in  \sym(M,\nabla,\omega)$, $A_{X} A_{X} Y = 0$.  
\end{enumerate}
\end{proposition}
\begin{proof} The Lie algebra $\sym(M,\nabla,\omega)$ is abelian,  
by Theorem \ref{thm:symplectic_autos}. 
Let  $\bar X$ be the
development image of $X$ and $Y$. Let $Z \in \sym(\bbA^n, \omega_{n})$ be a symplectic affine vector field on $\bbA^n$
centralising $\bar X$ and  $\bar Y$. By Lemma  \ref{lemma:scentraliserY}, we obtain 
$\omega_{n}(A_{\bar X}  A_{\bar X}  \bar Y, Z) = 0$.
Since $A(\Gamma)$ acts transitively, the vector fields $Z$ 
commuting with $\bar Y$ and $\bar X$ span all tangent spaces. Therefore, 
$A_{X} A_{X} Y = 0$. 
\end{proof}

If $(M,\nabla,\omega)$ is homogeneous then  $\sym(M,\nabla,\omega)$  forms a 
 subalgebra of the algebra of affine vector fields: 

\begin{proposition}   \label{prop:symplectic_autoshom}
If $M$ is homogeneous under the group of symplectic affine transformations then 
\begin{enumerate}
\item For $X, Y \in  \sym(M,\nabla,\omega)$, $A_{X} A_{Y}  = 0$. 
\item  $\sym(M,\nabla,\omega)$  forms a subalgebra of the
associative algebra $(\aut(M,\nabla, *_{\nabla})$.
\end{enumerate}
\end{proposition}
\begin{proof} Since $M$ is homogeneous, the vector
fields $Y \in  \sym(M,\nabla,\omega)$ span all tangent spaces.
Therefore,  $A_{X} A_{X} = 0$, for all $X \in \sym(M,\nabla,\omega)$.
Since $\sym(M,\nabla,\omega)$ is abelian, also  $A_{X} A_{Y}  = 0$. 
\end{proof}

\begin{remark} The structure and classification of algebras which are
appearing in Proposition \ref{prop:symplectic_autoshom} have been 
investigated in \cite{BC_1}. 
\end{remark}

%
%
%
%

\section{Homogeneous model spaces} \label{sect:HAM}

Let $G$ be a Lie group, and $H \leq G$ a closed subgroup. 
Then $X = G/H$ is a homogeneous space for $G$. 
If $X$ admits a flat affine connection $\nabla$, which is $G$-invariant,
then $(X,\nabla)$ is called a
\emph{flat affine homogeneous space}. 
Let  $\Gamma \leq G$ be a discrete subgroup, such that 
the quotient space $$  M =  \Gamma \,  \backslash G/H \; $$
is  a manifold. Since $\nabla$ is invariant by $\Gamma$, 
$M$ inherits a flat connection $\nabla$, such that
the natural map $(X,\nabla) \ra (M,\nabla)$ is an affine covering map.
Then $M$ is  called a  \emph{Clifford-Klein form} for $X= G/H$, and  
the  (simply connected) covering  space $(X, \nabla)$ is called a 
global model for $(M,\nabla)$. \\
 
The Auslander and Markus conjectures suggest a strong relationship 
between the geometry of a flat affine homogeneous space $(X,\nabla)$, 
and the existence of compact Clifford-Klein forms for $X$.  \\

The following result a consequence of Theorem \ref{thm:GHvol}:

\begin{theorem}  \label{thm:CKf}
Let $(X=G/H,\nabla)$ be a flat affine homogeneous space.
If $X$ has a compact Clifford-Klein form $(M,\nabla)$ with parallel volume 
then $(X,\nabla)$ is complete. 
\end{theorem}
In particular, if the universal covering $(X,\nabla)$ of a compact volume 
preserving flat affine manifold $(M,\nabla)$ is homogeneous under
its group of affine transformations, then $(M,\nabla)$ is complete. 
\\

In general, it seems difficult to understand the precise conditions
on an arbitrary flat affine homogeneous space $(X,\nabla)$, which ensure 
the completeness of $(X, \nabla)$. For flat affine Lie groups $(G,\nabla)$, 
a simple characterisation is known (see section \ref{sect:FALGs}).
An approach to the general problem will be introduced in  
section \ref{sect:crittrans} and section \ref{sect:Ccc}. 
In particular,  we will prove the following result: 

\begin{theorem} \label{thm:nilvol}
Let $(M, \nabla)$ be a homogeneous flat affine 
manifold for a nilpotent Lie group $G$. Then $M$ is complete
if and only if $G$ preserves a parallel volume form on $M$. 
\end{theorem}
The theorem is a direct consequence of a corresponding 
result for affine  homogeneous domains of nilpotent groups
(see section \ref{sect:niltransitive}). \\


In what follows, we will then mainly investigate \emph{compact homogeneous} 
affine manifolds. Before doing so,  we discuss the
development map of a homogeneous space and indicate the proofs
of Theorem \ref{thm:CKf} and Theorem \ref{thm:nilvol} in section
 \ref{sect:homdeviscovering}. Next we provide a structure theorem for volume 
preserving compact homogeneous affine manifolds in section \ref{sect:structurech}.
This result implies a direct proof (independent of 
Theorem \ref{thm:CKf}) that every compact volume preserving 
homogeneous affine manifold is complete. 
As another consequence of the structure theorem, 
we derive a classification theorem for compact 
flat Pseudo-\-Rie\-mannian homogeneous manifolds. 
We also show that a compact symplectic 
homogeneous affine manifold is \emph{diffeomorphic} 
to a torus. In section \ref{sect:homhol} we further discuss 
the properties of the holonomy groups of homogeneous
manifolds. 

\subsection{The development map of a homogeneous space}  
\label{sect:homdeviscovering}
Recall that a group action on $\bbA^n$ is called \emph{prehomogeneous} if it has a
Zariski-dense open orbit. A subset of $\bbA^n$ is \emph{semi-algebraic} if it is
defined by  polynomial equations and inequalities. The following is a basic
observation: 

\begin{proposition} \label{prop:hdevmap}
Let $(M,\nabla)$ be a homogeneous flat
affine manifold. Then the development map is a covering 
map, and its development image is a semi-algebraic
subset in $\bbA^n$. Let $G$ be a group which acts
transitively on $(M, \nabla)$.  Then its universal
covering group acts prehomogeneously on $\bbA^n$
by the development homomorphism. 
\end{proposition}
\begin{proof} Let $G$ be a Lie group which acts transitively 
on $(M,\nabla)$, and $\tilde{G} \leq \Aff(\tilde{M})$ be 
a covering group,  which lifts the action of $G$. 
Then the development $h(\tilde{G}) \leq \Aff(n)$ acts transitively 
on the development image $U= D(\tilde M)$. In fact, it follows that
$U \subseteq \bbA^n$ is an affinely homogeneous domain.
By Proposition \ref{prop:basic}, every affinely
homogeneous domain $U$ is a semi-algebraic subset of $\bbR^n$.
The development map $D: \tilde M \ra U$ is a covering map, because it identifies with 
a covering of homogeneous spaces 
$\tilde M = \tilde{G}/ \tilde H \ra  h(\tilde{G})/ L = U$,
which is induced by the locally faithful homomorphism 
$h : \tilde{G} \ra \Aff(n)$.
\end{proof} \hspace{1cm} 

\begin{prf}{Proof of Theorem \ref{thm:CKf}}
Since $X$ has a compact Clifford-Klein form with parallel volume, its 
development image is a homogeneous affine domain $U$, which is divisible
by volume preserving affine transformations. Hence, by Corollary \ref{cor:GHvol},
$U = \bbA^n$. Since $D$ is a covering, $(X,\nabla)$ must be complete.
\end{prf}   \hspace{1cm}  

\begin{prf}{Proof of Theorem \ref{thm:nilvol}} Assume $M$ has a transitive volume 
preserving action of some nilpotent Lie group. Then the universal covering $X$ of $M$
is a homogeneous space for a nilpotent Lie group $N$ of affine transformations, which preserves a parallel volume form. The development of $N$ preserves the parallel
volume on $\bbA^n$, and acts prehomogeneously. By Corollary \ref{cor:niltransitive},
the development of $N$ acts transitively on $\bbA^n$, and 
therefore $D(X) = \bbA^n$. Since $D$ is
a covering, it must be a diffeomorphism. Therefore, $M$ is complete. 

The converse is a direct consequence of  Corollary \ref{cor:niltransitive}.
\end{prf} 

\subsection{Compact homogeneous affine manifolds} \label{sect:comphom}
Robert Hermann \cite{Hermann_2} observed that a compact 
Pseudo-Riemannian manifold which admits a transitive group of
isometries must be complete (see \cite{Marsden} for a complete proof). 

\begin{example}  \label{ex:gln}
Consider $\GL(n,\bbR)$ as an open subset
in $\bbR^{n^2}$. Then the natural affine strucure on $\GL(n,\bbR)$
is invariant under left- and right-translations on the group. 
As a consequence of a result of Borel \cite{Borel}, 
there exist cocompact lattices $\Gamma \leq \GL(n,\bbR)$.
The compact manifold $M = \GL(n,\bbR) / \, \Gamma$ inherits
an affine flat structure, and $\GL(n,\bbR)$ acts 
as a transitive group of affine transformations on $(M,\nabla)$.
Thus $(M,\nabla)$ is a compact affinely homogeneous 
manifold, which is not complete. 
\end{example}

Markus' conjecture asserts that an (orientable) compact  flat affine manifold
is complete if and only if it admits a parallel volume form. The conjecture
is known to hold for compact affinely homogeneous affine manifolds. (This is proved in 
\cite{GH_2}.  See section \ref{sect:GHvol} for further discussion.) We 
obtain an independent  proof in Corollary 
\ref{cor:compact_hom} below. 


\subsubsection{Compact homogeneous manifolds with parallel volume} 
\label{sect:structurech}
We have the following structure result for 
compact homogeneous flat affine manifold with 
parallel volume:

\begin{theorem} \label{thm:compact_hom}
Let $(M, \nabla)$ be a compact homogeneous flat affine manifold with volume preserving affine structure. Then
the following hold: 
\begin{enumerate}
\item  $(M,\nabla)$ is affinely diffeomorphic to a quotient space of
a simply connected nilpotent affine Lie group $(N,\nabla)$ with biinvariant affine structure $\nabla$ by a discrete subgroup 
$\Gamma \leq N$.  
\item $(M, \nabla)$ is geodesically complete.
\end{enumerate}
Furthermore, every parallel tensor field on $(M,\nabla)$, which is 
preserved by the action of $N$, pulls back to a biinvariant (and parallel) tensor
field on $N$. 
\end{theorem}
\begin{proof} By Proposition \ref{prop:compvp}, the group of affinities 
$\hat{\Aut}(M, \nabla)^0$ develops to a connected unipotent Lie group. Moreover, 
$M$ is a compact homogeneous space of 
the simply connected nilpotent Lie group $N= \hat{\Aut}(M, \nabla)^0$.
Since $N$ acts almost effectively on $M$, if follows that 
$\dim N = \dim M$ (see  \cite{Malcev}).  
Therefore, $M = N / \Gamma$ is a quotient of $N$ 
by a uniform discrete subgroup $\Gamma \leq N$. Pull back $\nabla$  to obtain a left-invariant affine connection $\nabla$ on $N$, which is right-invariant under the lattice 
$\Gamma$. Since $\Gamma$ is Zariski-dense in the adjoint representation 
(see \cite{Raghunathan}), $\nabla$ is a biinvariant connection on $N$. (In fact,  
$\nabla$ is determined by the associative product structure $*_{\nabla}$ on
the Lie algebra $\aut(M,\nabla) = \lie{n}$.) The volume preserving 
left-invariant connection $\nabla$ on the nilpotent Lie group $N$
is complete. (See section \ref{sect:compuni}.) Therefore, $M$ is complete.  
\end{proof}

\begin{remark}
By a  result of Mostow \cite{Mostow} (see also, \cite{Raghunathan, Tol}), 
a finite volume
homogeneous space of a solvable Lie group is compact. 
Thus, the above proof also shows that  a finite volume homogeneous affine manifold is compact and geodesically complete.
\end{remark}

The Markus conjecture holds  for compact homogeneous affine manifolds:

\begin{corollary} \label{cor:compact_hom}
Let $(M, \nabla)$ be a compact homogeneous flat affine manifold. 
Then $(M, \nabla)$ is complete if and only if $(M, \nabla)$ has
parallel volume. 
\end{corollary}
\begin{proof}  If $(M,\nabla$ has parallel volume Theorem \ref{thm:compact_hom} 
applies. 

Assume now that $(M, \nabla)$ is complete. Then, as already remarked in 
the beginning of section \ref{sect:compvp}, $\Aut(M,\nabla)^0$ is a nilpotent Lie group,
which develops to a unipotent subgroup of $\Aff(n)$. 
Therefore, as in the proof of Theorem \ref{thm:compact_hom}, $M = N/\Gamma$, where
$N$ is a nilpotent Lie group with complete left invariant affine connection $\nabla$, 
such that $(N,\nabla) \ra (M,\nabla)$ is an affine covering. 
Since $N$ is nilpotent and $\nabla$ is complete,  the left-multiplications
on $N$ are volume preserving affine transformations (see  section \ref{sect:compuni}). Since $N$ is nilpotent, also the right-multiplications preserve the left-invariant 
parallel volume on $N$.  Since $hol(\pi_{1}(M)) = \Gamma$ is 
contained in the right-multiplications on $N$, $\Gamma$ 
preserves the parallel volume on $N$. 
This shows that $(M,\nabla)$ has parallel volume.  
\end{proof}

Note also that Theorem \ref{thm:compact_hom} strengthens 
the completeness result in \cite{Hermann_2,Marsden}
considerably, since it is only required that the group of affinities $\Aut(M,\nabla_{\g})$ acts
transitively on $M$, to ensure completeness. 

\subsubsection{Pseudo-Riemannian examples}  \label{sect:comphPSR}

Let $(M,\g)$ be a homogeneous flat
Pseudo-Riemannian manifold, which is compact.  
Then, by the above,  $(M, \nabla_{\g})$ is complete.
Moreover, the following structure theorem holds: 

\begin{theorem}  \label{thm:structure_ch}
Let $M$ be a compact (or finite volume)
homogeneous flat Pseudo-Riemannian manifold. Then
$M$ is isometric to a quotient of a flat Pseudo-Riemannian 
Lie group $N$ with biinvariant metric. 
\end{theorem}
\begin{proof} In fact, by Theorem \ref{thm:compact_hom}, $(M, \g)$
is isometric to a quotient of a nilpotent flat Pseudo-Riemannian Lie group
$(N, \g)$ with biinvariant metric $\g$, such that the natural
(orbit-) map $N  \ra M = N/\Gamma$ is a Pseudo-Riemannian
covering. 
\end{proof}

Theorem \ref{thm:structure_ch} implies the following strong geometric \emph{rigidity property} for compact homogeneous flat Pseudo-Riemannian manifolds:

\begin{corollary}   \label{cor:rigidity_ch}
Let $M$ and $M'$ be compact 
homogeneous flat Pseudo-Riemannian manifolds with isomorphic
fundamental groups.  Then $M$  and $M'$ are affinely diffeomorphic.
\end{corollary}
\begin{proof}  As above, we write $M = N/\Gamma$ and $M' = N' / \Gamma'$, where
$\Gamma$, $\Gamma'$ are lattices in the simply connected nilpotent 
Lie groups $N$ and $N'$, respectively. 
By Malcev-rigidity, cf.\ \cite{Malcev}, every isomorphism $\phi: \Gamma \ra \Gamma'$
extends to an isomorphism of Lie groups $\Phi: N \ra N'$.  We contend that 
$\Phi$ is an affine isomorphism of metric Lie groups $(N,\g) \ra (N',\g')$. 
In fact, both $\g$ and $\g'$ are biinvariant metrics, and therefore (see 
\cite[Proposition 11.9]{ONeill}), the Levi-Civita connection $\nabla$ (respectively $\nabla'$)
is the canonical torsion-free connection on the Lie group $N$ (respectively $N'$). 
That is, for left-invariant vector fields $X$ and $Y$, 
$\nabla_{X} Y =  {1 \over 2} [X,Y]$.  In particular, it follows that the isomorphism of Lie groups $\Phi$ is 
affine with respect to $\nabla$ and $\nabla'$.  Therefore, the induced 
diffeomorphism $\bar{\Phi}: M \ra M'$ is affine. 
\end{proof}

In particular, the fundamental group $\pi_{1}(M)$ determines the Levi-Civita 
connection on a compact homogeneous flat Pseudo-Riemannian
manifold $M$ up to affine equivalence. 

\begin{example} If\/ $\pi_{1}(M)$ is abelian then $\pi_{1}(M) \cong \bbZ^n$, and $M = T^n$
is diffeomorphic to a torus. Moreover,  by 
Corollary \ref{cor:rigidity_ch}, every homogeneous flat manifold
$(T^n, \g)$ is isometric to a quotient of  $\bbE^{s}$ by a lattice of translations. 
The affine structure of\/ $(T^n, \g)$ 
is the uniquely determined translational structure on $T^{n}$. 
\end{example}

\begin{remark} In fact, \emph{every} flat metric on the torus $T^{n}$ \emph{is} 
homogeneous, and thus isometric to a quotient of  $\bbE^{s}$ by a lattice of translations.
To the contrary, the set of affine equivalence classes of volume preserving  
complete affine structures on $T^n$ is rather large. See \cite{BauesG, BG}
for further discussion.
%
%
\end{remark}

A nilpotent Lie group $N$, which acts transitively by isometries 
on a compact flat Pseudo-Riemannian manifold $M$ is of rather restricted type:

\begin{lemma} Let  $M = N/\Gamma$  be a compact 
homogeneous flat Pseudo-Riemannian manifold. 
Then $N$ is a two-step nilpotent Lie group.  
\end{lemma}
\begin{proof}  
Remark  (see \cite[loc.\ cit.]{ONeill}), that the curvature of
the Levi-Civita connection $\nabla$ of the biinvariant 
flat metric $\g$ on $N$ is computed as, $
R_{\nabla}(X,Y) Z = {1 \over 4} [X, [Y,Z]]$, for all left-invariant vector fields
$X,Y,Z$. Since $\g$ is flat, this implies $[X, [Y,Z]] = 0$. It follows that $N$ is two-step
nilpotent. 
\end{proof}

Note also that every quotient space $M = N/\Gamma$ of a flat Pseudo-Riemannian 
Lie group $N$ with biinvariant metric by a discrete subgroup $\Gamma$ is a
homogeneous flat Pseudo-Riemannian manifold. In fact, left-multiplication in $N$
induces a transitive isometric action of $N$ on $M = N / \Gamma$. 

In section \ref{sect:biinvariantmetrics}, we will discuss the structure of flat Pseudo-Riemannian Lie groups with biinvariant metric in detail. 
In particular, see Corollary \ref{cor:semprod},
we explicitly construct a  family of such Lie groups, which gives 
rise to an interesting class of examples of compact (complete) homogeneous flat 
Pseudo-Riemannian manifolds with non-abelian fundamental group. 
This is summarised in Corollary \ref{cor:familypsrh} below. \\

First, we state a few general consequences concerning the classification
of \emph{compact} homogeneous flat Pseudo-Riemannian manifolds:

\begin{corollary} Let  $M$  be a compact 
homogeneous flat Pseudo-Riemannian manifold. Then the
following hold:
\begin{enumerate}
\item If $\dim M \leq 5$ then $M$ is isometric to a translation torus (that is,
$M$ is a quotient of\/  $\bbE^s$, by a lattice of translations).
\item The fundamental group $\Gamma$ of $M$ is two-step nilpotent. 
\end{enumerate}
\end{corollary}

The first examples of compact homogeneous flat Pseudo-Riemannian manifolds
with \emph{non-abelian} fundamental group arise in dimension six.
More generally, there is an interesting class of examples, as follows: 

\begin{corollary} \label{cor:familypsrh}
Let  $N / \Gamma$ be a compact two-step nilmanifold (where 
$\Gamma$ is a lattice in the two-step nilpotent Lie group $N$, $n  = \dim N$).  
Then there exists a flat $\rank n$  torus torus bundle $$ T^n \ra M \ra N / \Gamma $$ 
such that  $M$ admits a homogeneous flat Pseudo-Riemannian metric
of signature $s= (n,n)$ with totally geodesic and isotropic fibers $T^n$. 
\end{corollary}
\begin{proof} Let $A: N \ra \GL(\lie{n}^*)$ denote the coadjoint 
representation of $N$ on the dual $\lie{n}^*$ of its  Lie algebra $\lie{n}$. 
By Corollary \ref{sect:biinvariantmetrics}, the corresponding 
semidirect product $N \rtimes_{A} \lie{n}^*$ admits a
flat biinvariant Pseudo-Riemannian metric of signature $(n.n)$,
$n =\dim M$. Now $\Gamma$ is a lattice in $N$, and under the
representation $A$, maps to a finitely
generated subgroup of the group of unipotent upper
triangular matrices with rational coefficients relative
a basis of $\lie{n}$. In particular, $A(\Gamma)$ preserves
a lattice $\Lambda \cong \bbZ^n$ in $\lie{n}^*$. The
semi-direct product $\Gamma' = \Gamma \rtimes_{A} \Lambda$
is a lattice in $N \rtimes_{A} \lie{n}^*$, and the manifold $M = N \rtimes_{A} \lie{n}^* /\, \Gamma'$
has the required properties.
\end{proof}

In particular,  let $\Gamma$ be a torsion-free finitely generated two-step nilpotent
group.  Then $\Gamma$ arises as a quotient $\Gamma = \Gamma'/ \Lambda$, 
where $\Gamma'$ is the holonomy of a compact homogeneous flat 
Pseudo-Riemannian manifold, and $\Lambda$ is an abelian normal
subgroup of $\Gamma'$, with $\rank \Lambda  = \rank \Gamma$.

\begin{example} In dimension $n=6$, there exist up to affine equivalence
two classes of compact homogeneous flat Pseudo-Riemannian manifolds. 
The six-dimensional translation torus $T^6$, and torus bundles 
$T^3 \ra M_{6, \Gamma} \ra H_{3}/ \Gamma$, where $\Gamma$ is a lattice
in the three-dimensional Heisenberg group $H_{3}$. 
\end{example} 

\subsubsection{Symplectic examples}
As noted before, and contrasting the Pseudo-Riemannian situation,
only tori admit the structure of a homogeneous symplectic flat
affine manifold. However, such a structure is, in general, not affinely
equivalent to a translation torus:
 
\begin{corollary}  Let $M$ be a compact (or finite volume)
homogeneous symplectic flat affine manifold. Then $M$
is a compact abelian flat affine symplectic Lie group. 
In particular,  $M$ is diffeomorphic to a  torus. 
\end{corollary} 
\begin{proof}  By Theorem 
\ref{thm:compact_hom}, $(M,\nabla)$ 
is a quotient of an affine
symplectic Lie group $N$ with biinvariant symplectic
structure. By Proposition \ref{cor:bisymplectic}, $N$ is abelian.  
(To show that $M$ is a torus,  we might have applied the
classification results for connected groups of symmetries
of compact symplectic manifolds given in \cite{ZB,Guan}, as well.)
\end{proof}

Note that the space of affine equivalence classes
of flat symplectic connections on the torus $T^{2n}$ is very large. 
In fact, as is shown in \cite{BC_1}, the set of 
simply connected abelian flat affine symplectic Lie groups
forms an algebraic variety of (real) dimension $n^3$.  
For the precise classification result and 
further geometric properties of abelian flat affine symplectic Lie groups, 
see \cite{BC_1}. 

\subsection{Holonomy of homogeneous affine manifolds} \label{sect:homhol}

The fundamental group of a homogeneous  flat Riemannian 
manifold is abelian, and its holonomy acts by translations. 
Starting in dimension five, 
there do exist complete flat homogeneous Pseudo-Riemannian manifolds
with abelian but non-translational holonomy group, see  \cite{Wolf_1}.
As Example \ref{ex:gln} shows, the holonomy of a  (compact)
homogeneous  affine manifold can be far away from
being abelian or nilpotent. This example is however a 
non-complete manifold.

\begin{theorem} \label{thm:homhol}
Let $\Gamma \leq \Aff(n)$ be the holonomy group of 
a complete affinely homogeneous affine manifold.  Then $\Gamma$
is a unipotent subgroup of $\Aff(n)$. 
\end{theorem}
\begin{proof} Since $M$ is homogeneous, the development $G$ of $\Aff(M,\nabla)$
acts transitively on the development image of $M$, which is $\bbA^n$.
Moreover, $\Gamma$ centralises $G$. By Proposition \ref{prop:centisu},
the centraliser of the transitive group $G$ is a connected
unipotent group. In particular, $\Gamma$ is unipotent. 
\end{proof}

In particular: 
\begin{corollary} The fundamental group of a 
complete affinely homogeneous affine manifold is nilpotent. 
\end{corollary} 

\subsubsection{Holonomy of homogenous flat Pseudo-Riemannian manifolds}
\label{sect:homhol_psr}
The above result can be strengthened considerably for 
homogeneous flat Pseudo-Riemannian manifolds 
even without the assumption of completeness. 

\begin{proposition}  \label{prop:homhol_psr}
Let $\Gamma \leq \E(s)$ be the holonomy group of 
a  homogeneous flat Pseudo-Riemannian manifold.  Then $\Gamma$
is a unipotent group, and $\Gamma$ is nilpotent of nilpotency class
at most two.  
\end{proposition}
\begin{proof} Let $(M,\g)$ be a homogeneous flat Pseudo-Riemannian 
manifold. The development image of $M$ is a 
Pseudo-Riemannian domain $U$ in $\bbE^s$, which is homogeneous 
with respect to the development $G$ of $\Isom(M,\g)$.
The holonomy group $\Gamma$ preserves the domain $U$ and is
contained in the centraliser of $G$ in the full isometry group $\E(s)$. 
By Corollary \ref{cor:psrcisunip}, 
the centraliser is a connected unipotent group of
nilpotency class at most two. In particular, $\Gamma$
has this property.
\end{proof}

\begin{remark}
The proposition leads to a classification theory of  flat homogeneous Pseudo-Riemannian manifolds, which are complete. 
Wolf studied the classification of complete homogeneous Pseudo-Riemannian manifolds with \emph{abelian} holonomy $\Gamma$, discovering interesting phenomena, 
see  \cite{Wolf_1,Wolf_2,Wolf_3}. In section \ref{sect:comphPSR}, we proved the existence of a six-dimensional compact (complete) homogeneous Pseudo-Riemannian manifold with non-abelian holonomy $\Gamma$. Moreover,  all examples constructed in section \ref{sect:comphPSR} relate to certain two-step nilpotent Lie groups.  \emph{It remains to join these
results to a full structure theory for complete flat homogeneous 
Pseudo-Riemannian manifolds.} 
\end{remark}

The determination of all homogeneous flat 
Pseudo-Riemannian manifolds (that is, including non-complete examples) 
seems not at hand, because of the difficulties
associated with the classification of homogeneous Pseudo-Riemannian domains. 
(See section \ref{sect:psrdomains}.)\\

We mention furher that Proposition \ref{prop:homhol_psr} admits an analogous version for symplectically homogeneous affine manifolds. This is a consequence of
Proposition \ref{prop:scisunip} below: 

\begin{proposition}  Let $\Gamma \leq \Aff(\omega_{n})$ be the holonomy group of 
a symplectically homogeneous affine manifold with parallel symplectic
form.  Then $\Gamma$ is an abelian group of unipotent symplectic affine transformations.
\end{proposition}

\section{Flat affine Lie groups} \label{sect:FALGs} 
We consider a particular tractable class of homogeneous spaces,  
namely those which admit a \emph{simply transitive} group of 
equivalences. Many phenomena can be illustrated, 
and in many cases examples constructed from flat affine Lie groups 
serve as the basic building blocks of the theory. 

\subsection{Left-invariant geometry on Lie groups} 
Let $G$ be a Lie group. Recall that a geometric structure 
on $G$ is called {left-invariant if it is preserved by 
all left-multiplication maps $l_{g } : G \ra G$,  $g \in G$.
A connection $\nabla$ on $G$ is called
a \emph{left-invariant connection} if, for all $g \in G$,
$l_{g}$ is a connection preserving diffeomorphism.  
A Pseudo-Riemannian metric is
called left-invariant if all left-multiplications are isometries. 
More generally, if $G$ has the structure of an $(X,A)$-manifold
then the structure on $G$ is called left-invariant if all
left multiplications are $(X,A)$-equivalences. 

\begin{definition} Let $G$ be a Lie group and $\nabla$
a connection on $G$.  If $\nabla$ is flat affine and left-invariant, 
then $(G,\nabla)$ is called a \emph{flat affine Lie group}.
Moreover, $\nabla$ is called a flat left-invariant connection of 
type $\cA$, where $\cA \leq \Aff(n)$, if $(G,\nabla)$ has a compatible 
\emph{left-invariant} locally homogeneous structure of type $(\bbA^n, \cA)$.
\end{definition}
In the latter  situation we shall call the pair $(G,\nabla)$ a \emph{flat affine Lie group
of type $\cA$}. If $\g$ is a flat left-invariant metric, then $(G,\g)$ is called a \emph{flat Pseudo-Riemannian Lie group}.  

\begin{example} Let $\nabla_{\g}$ be the Levi-Civita connection
of a  flat Pseudo-\-Rie\-mannian Lie group. Then $(G, \nabla_{\g})$ is a
flat affine Lie group of type $E({s})$. 
\end{example}


\subsubsection{Left-symmetric algebras}
There is a one to one correspondence of left-invariant
flat affine connections $\nabla$ on the simply connected group $G$, 
with left-symmetric algebra structures on the Lie algebra $\lg$
of left-invariant vector fields on $G$. In fact, if $X,Y$ are 
left-invariant then $\nabla_{X} Y = - Y {*}_{\nabla} X$ is a left-invariant vector field.
Hence $\nabla$ defines a left-symmetric product on $\lg$. (Compare
section \ref{sect:LRalgs}.) The flat connection $\nabla$ is a {biinvariant
connection}  if and only if the induced product is associative.  

\begin{example} Let $N$ be a two-step nilpotent Lie group. The Lie algebra
of left-invariant vector fields $\lie{n}$ thus  satisfies $[\lie{n}, [\lie{n}, \lie{n}]] =  \{ 0 \}$.
It follows that the canonical torsion-free connection on $N$, defined by
$$ \nabla_{X} Y = {1 \over 2} [X,Y]  \; \;  , \,  X, Y \in \lie{n} \; , $$ is flat. Therefore, 
$N$ is naturally a flat affine Lie group with biinvariant flat connection $\nabla$.
The bilinear product $X \cdot Y := \nabla_{X} Y$ defines an associative algebra structure
on $\lie{n}$.  (See \cite{Scheune1,Scheune2} for more constructions of left-invariant
flat connections on nilpotent Lie groups.)  
\end{example}

\subsection{The development map of a flat affine Lie group}

Let $\rho: G \ra \Aff(n)$ be a homomorphism. Recall that 
$\rho$ is called an \emph{\'etale affine representation}  (see Definition \ref{def:etrep})
if there exists $x \in \bbR^n$
such that the orbit map at $x$ $$ o_{x}: G \ra \bbA^n,  \;  \; 
g  \mapsto \rho(g) \cdot x $$ is a local diffeomorphism onto 
an open subset $U$ of $\bbA^n$. 
The representation $\rho$ is
\emph{simply transitive} if $o_{x}$ is a diffeomorphism
for one (and, hence, for all) $x \in \bbA^n$. 

\begin{example}
Let $\rho:  G \ra \cA$ be an \'etale affine representation, and $x \in \bbA^n$
such that $o_{x}: G \ra \bbA^n $ is a local diffeomorphism. By pull back of the 
standard flat connection on $\bbA^n$ along $o_{x}$, we obtain a 
flat left invariant connection of type $\cA$ on $G$. 
\end{example}

Moreover, every 
flat left invariant connection of type $\cA$ arises in this way:

\begin{proposition} \label{prop:flatLiegroups}
Let $G$ be a simply connected Lie group, and $\nabla$ 
a  a torsion-free flat connection on $G$. Then the following are equivalent: 
\begin{itemize}
\item $\nabla$ is a left-invariant connection of type $\cA$.
\item There exists an \'etale affine representation $\rho: G \ra \cA$, and $x \in \bbR^n$, 
such that $\nabla$  and its compatible $(\bbA^n, \cA)$-structure are a 
pullback along the orbit map $o_{x}$. 
\end{itemize}
Moreover, the connection $\nabla$ is complete if and only if $\rho$ is
simply transitive. 
\end{proposition} 
\begin{proof}
Let $D: G \ra \bbA^n$ be any development map for the compatible 
left-invariant $(\bbA^n,\cA)$-structure. 
Let $h$ denote the holonomy homomorphism with respect to $D$.
Since the locally homogeneous structure is compatible with $\nabla$, 
the map $D$ is also an affine map $(G, \nabla) \ra \bbA^n$. 
Since the maps $l_{g}: G \ra G$ are equivalences for the locally homogeneous
$(\bbA^n,\cA)$-structure on $G$, 
$\rho(g) = h(l_{g})$ defines a representation $\rho: G \ra \cA$.

We put $x = D(1)$.
Then, by \eqref{eq:hol}, $$ D(g) = D(l_{g} 1) =  \rho(g) D(1) = o_{x}(g) \; . $$  
Thus, the orbit map for $\rho$ at $x= D(1)$ is the development map $D$.  
This shows that $\rho$ is an \'etale affine representation,
and the left-invariant $(\bbA^n,\cA)$-structure on $G$ is obtained by pullback along $o_{x}$.
\end{proof}

\subsubsection{The boundary of the development image} 

An  \'etale affine representation $\rho: G \ra \cA$, which  induces 
a flat left-invariant connection 
$\nabla$ of type $\cA$, as in Proposition \ref{prop:flatLiegroups}, 
will be called a \emph{compatible} \'etale  affine representation for $\nabla$. 

\begin{corollary} The development map of a flat affine Lie group $(G, \nabla)$ 
is a covering map onto its image, 
and it is an orbit-map for any compatible \'etale  affine representation
of its universal covering group $(\tilde G, \nabla)$. 
\end{corollary}
\begin{proof} Since $D: \tilde G \ra \bbA^n$ is a local diffeomorphism and an orbit-map, it is a covering map.
\end{proof}

 Let   $\Ad: G \ra \GL(\lg)$ denote the adjoint representation of $G$. Furthermore, 
 if $(G, \nabla)$ is affine we put  $\det_{\nabla} l: G \ra \bbR^{\neq 0}$ for the character which is obtained by the left-action of $G$ on parallel volume forms.

\begin{theorem}
Let $(G,  \nabla)$ be a  flat affine Lie group. 
Then the following hold: 
\begin{enumerate}
\item  Either $(G,\nabla)$ is complete or the development image of $(G,  \nabla)$ is a connected component of the complement of a non-empty hypersurface 
in $\bbR^n$, $n = \dim G$. In particular, 
the development image is a semialgebraic subset of $\bbR^n$.
\item  Assume that $G$ is simply connected. 
Then the  connection $\nabla$ is complete if and only if 
$\det \Ad(g)^{-1}  \det_{\nabla} l_{g} =  1$, for all $g \in G$. 
\end{enumerate}
\end{theorem} 
\begin{proof} 
In fact, the development image of $G$ is the orbit $\rho(\tilde{G}) x$, 
for a  compatible affine \'etale representation $\rho$ for the induced flat
structure on $\tilde{G}$. By section \ref{sect:etreps}, $D(\tilde{G}) = \rho(\tilde{G}) x$
is a connected component of the open semi-algebraic subset 
$  U_{\delta} = \{ x \in \bbR^n \mid  \delta(x) \neq 0 \} \subseteq  \bbR^n$, where
$\delta = \delta(\rho)$ is the relative invariant for $\rho$. This implies 1.\
Moreover, since $D$ is a covering map, $\nabla$ is complete if and
only if $D$ is surjective.  
If $G$ is simply connected, the volume character of $G$
satisfies $\det_{\nabla} l_{g} =  \det \rho_{\nabla} (g)$, for any compatible \'etale affine representation $\rho_{\nabla}: \tilde G \ra \Aff(n)$.  
Thus 2.\ follows from Theorem \ref{thm:relinv_strans}
below.
\end{proof}

\subsubsection{Completeness and unimodularity of flat affine Lie groups}
\label{sect:compuni}

Let  $(G, \nabla)$ be a flat affine Lie group. If $\nabla$ is complete 
then, by  Proposition \ref{prop:flatLiegroups}, the universal covering Lie group 
of $G$ acts simply transitively on affine space. This, more or less, 
reduces the study of complete affine Lie groups to 
a study of simply transitive representations on affine space.
(See \cite{Auslander}, 
for a seminal paper on simply transitive groups of 
affine transformations.)

\begin{corollary}
Let  $(G, \nabla)$ be a flat affine Lie group. If $\nabla$ is complete then $G$ is 
a solvable Lie group.   
\end{corollary}
\begin{proof} 
By Levi's theorem $\tilde{G}$ splits as a semi-direct product of
a semisimple group $H$ and a solvable Lie group $S$.  
Now $\tilde{G}$ can not have a reductive Levi-component $H$
since every affine representation of $H$ has a fixed point.
Hence, $\tilde{G}=S$ must be solvable. 
\end{proof}

Recall that a Lie group $G$ is called \emph{unimodular}
if $\det \Ad(g)= 1$, for all $g \in G$. If a flat Lie group
$(G, \nabla)$ has a left-invariant
parallel volume form, we call it  \emph{volume preserving}.

\begin{corollary} \label{cor:Gcomplete_vp}
Let $(G,\nabla)$ be a volume-preserving 
flat affine Lie group. Then $(G,\nabla)$ is complete if and only if $G$
is unimodular. 
\end{corollary}

\begin{corollary} \label{cor:Gcomplete_um}
Let $(G,\nabla)$ be a flat affine Lie group, where $G$
is unimodular. Then $(G,\nabla)$ is complete if and only if $(G,\nabla)$
is a volume-preserving flat affine Lie group. 
\end{corollary}

\begin{example}[Semi-simple Lie groups are not affine.]   
Let $S$ be a semisimple Lie group. Then $S$ does not 
admit a flat left-invariant connection. 
\end{example}

For further generalisation, see section \ref{sect:relisc}. 

\subsection{Flat Pseudo-Riemannian Lie groups}

\begin{corollary} Let $(G,\g)$ be a flat Pseudo-Riemannian Lie group. Then 
$(G,\nabla_{\g})$ is complete if and only if $G$ is unimodular. In particular,
a unimodular  flat Pseudo-Riemannian Lie group is solvable. 
\end{corollary}

\begin{example} A  flat Pseudo-Riemannian Lie group  $(G,\g)$, with
biinvariant metric $\g$ is complete and $G$ is a two-step nilpotent 
Lie group. (See section   \ref{sect:biinvariantmetrics} below.)
\end{example}

The above  corollary already excludes a large class of Lie groups
from carrying flat left-invariant metrics:

\begin{example} Let $G= \GL(n,\bbR)$. Since $G$ is an open subset of 
$\bbR^{n^2}$, it has a natural structure of a flat affine Lie group. 
However, $G$ does not admit a volume preserving flat left-invariant connection. 
In particular, $G$ does not admit a flat Pseudo-Riemannian 
metric. 
\end{example}

The latter argument holds for all \emph{reductive} Lie groups $G$.
Here we call Lie group  \emph{reductive} if it is an almost semi-direct 
product of its center and a semisimple Lie group. Such a group is
necessarily unimodular. Hence: 

\begin{corollary} Let $G$ be reductive. Then 
$G$ does not admit a volume preserving flat left-invariant connection. 
In particular, $G$ does not admit a left-invariant flat Pseudo-Riemannian 
metric. 
\end{corollary}

\begin{remark}
Similar restrictions, as stated here for the Pseudo-Riemannian
case,  also hold for symplectic affine Lie groups. 
\end{remark}

Concerning the existence of flat Pseudo-Riemannian Lie groups  $(G,\g)$,
where $G$ is not solvable, we observe: 

\begin{example} Every flat Lorentzian Lie group $G$ is solvable. This is 
a consequence of the classification of homogeneous domains in $\bbE^{n-1,1}$.
In fact, every such domain is diffeomorphic to $\bbR^n$. (See section
\ref{sect:psrdomains}, repectively \cite{DuIh_1}.) Therefore, the 
universal covering group $\tilde{G}$ of $G$ is diffeomorphic to $\bbR^n$. 
It can not have a reductive part, since $\widetilde{\SL(2,\bbR)}$ does not
admit a faithful representation. 
\end{example}

Not every flat Pseudo-Riemannian Lie group is solvable, as the following
example shows: 

\begin{example}  Let $U = \GL(2,\bbR) \subset \bbA^4$.
Then the dual tube domain $\check{T}(U) \leq \bbE^{4,4}$ 
is a homogeneous Pseudo-Riemannian domain, 
which has a simply transitive isometric action of the group 
$\GL(4,\bbR) \ltimes (\bbR^4)^*$. (Compare section \ref{sect:dtc}.)
\end{example}

\begin{question} Are flat Pseudo-Riemannian Lie groups of signature $s=(n-2,2), (n-3,3)$ 
always solvable? 
\end{question}

Milnor \cite{Milnor_2} gives the classification and structure of 
simply connected \emph{flat Riemannian Lie groups}. These
are automatically complete. In \cite{GM} (see also \cite{Aubert-Med}
for the case of nilpotent Lie groups) the classification
of all simply connected  \emph{complete} flat Lorentzian Lie groups is 
described.

\subsubsection{Flat biinvariant metrics on Lie groups}  \label{sect:biinvariantmetrics}
 A left-invariant metric $\g$ on $G$ is called
 \emph{biinvariant} if the right-multiplications of $G$ are isometries as well.
 If $G$ is connected,  $\g$ is biinvariant if 
and only if the induced scalar product $\g$ on the Lie algebra 
$\lg$ is skew with respect to
the Lie bracket.  That is, $\g$ is biinvariant iff, for all $X,Y,Z \in \lg$, 
\begin{equation} \label{eq:biinvmetric}
  < [X, Y ], Z   >  \;   =  \,  - <  Y , [X, Z ]  > \; .   
\end{equation}
The following fact was already used in the proof of Corollary \ref{cor:rigidity_ch}: 
\begin{lemma} Let $(G,\g)$ be a Pseudo-Riemannian Lie group 
with biinvariant metric. Then the Levi-Civita  connection $\nabla= \nabla_{\g}$
 is the canonical torsion-free connection on $G$, that is,
 $\nabla_{X} Y = {1 \over 2} [X,Y] $, for all $X,Y \in \lg$.
 Moreover,  $(G,\g)$ is flat if and only if $G$ is two-step nilpotent.
\end{lemma}

\begin{corollary} Every flat Pseudo-Riemannian Lie group with
biinvariant metric is complete.
\end{corollary}


\subsubsection{Construction of flat biinvariant metric Lie algebras}  
\label{sect:biinvariantmetricLAs}
We describe now a method, which allows to construct examples of flat biinvariant 
metric Lie algebras $(\lg, \g)$. This, of course, implies the construction of
flat Lie groups with biinvariant metric. 
\paragraph{Coadjoint extensions}
Let $\lie n$ denote a Lie algebra, and $\lie{n}^*$ 
its dual vector space.  By the coadjoint representation $\ad^*: \lie n \ra \lie {gl} (\lie{n}^{*})$, 
$\lie{n}^*$ is a module for $\lie{n}$. Let $\omega \in Z^2(\lie{n}, \lie{n}^*)$
be a two-cocycle (cf.\ Appendix \ref{sect:Lieco}). To $\omega$ there
belongs a \emph{Lie algebra extension} 
$$  0 \ra \lie{n}^* \ra  \lie{t}_{\lie{n},\omega} \ra \lie{n} \ra 0 \; , $$
where  $\lie{n}^*$ is an abelian ideal in $ \lie{t}_{\lie{n},\omega}$. 
Put $W = \lie{n}^* \oplus \lie{n}$ for  
the direct sum of vector spaces. Explicitly, the Lie algebra 
$ \lie{t}_{\lie{n},\omega}$ arises from the Lie product on $W$,
which is declared by 
$$  [ (\lambda, v), (\lambda', v') ] =  (  \ad^*(v) \lambda' - \ad^*(v') \lambda+ 
\omega(v,v'), [v, v']) \; . $$ 
(See, for example, \cite{MacLane} for the general theory of Lie-algebra extensions.)\\

Next we define a natural scalar product $\g_{o}$ of signature
$s = (n,n)$ on $W$, where $n = \dim \lie{n}$. With respect to 
$\g_{o}$,  $\lie{n}$ and $\lie{n}^*$ are dual
totally isotropic subspaces under $\g_{o}$.  Moreover, if 
$v \in \lie{n}$ and $\lambda \in \lie{n}^*$, we put 
$$ < \!v, \lambda \! >_{o} = \lambda(v) \; .$$  
We put $V = \lie{n}$, and
we define a three-form $F_{\omega} \in \bigwedge^2 V^* \tensor{}  V^*$
by putting $$ F_{\omega}(u,v,w) =  \; < \omega(u,v), w >_{o} \; .$$ 

The following lemma is easily verified by direct calculation: 

\begin{lemma} Let $\lie{n}$ be two-step nilpotent, and $\omega \in Z^2(\lie{n}, \lie{n}^*)$. Then the  following hold:
\begin{enumerate}
\item The canonical scalar product $\g_{o}$ on $ \lie{t}_{\lie{n},\omega}$ 
satisfies \eqref{eq:biinvmetric} if and only if  $F_{\omega}$ is alternating 
(that is, if $F_{\omega} \in  \bigwedge^3 V^*$).  
\item The Lie algebra $ \lie{t}_{\lie{n},\omega}$ is two-step nilpotent, iff, for 
all $u,v,w, z \in \lie{n}$, \begin{equation} \label{eq:tistwostep}
 F_{\omega}(u,[v,w], z)  =   F_{\omega}(v,w, [u,z]) \; . 
\end{equation}
\end{enumerate} 
\end{lemma}
\begin{proof} Note, since $\lie{n}$ is two-step nilpotent, triple commutators in $ \lie{t}_{\lie{n},\omega}$
satisfy $$ [ (\lambda, u), [ (\mu, v), (\tau, w)  ]] = (ad_{u}^* \omega(v,w) + \omega(u, [v,w]) 
, 0) \; .  $$
This proves 2.
\end{proof} 

Thus we have, in particular:

\begin{proposition}
Let $\lie{n}$ be two-step nilpotent, and $\omega \in Z^2(\lie{n}, \lie{n}^*)$.
If $F_{\omega}$ is alternating and satisfies
\eqref{eq:tistwostep}, then the metric Lie algebra 
$( \lie{t}_{\lie{n},\omega},\g_{o})$  is flat and $\g_{o}$ is biinvariant. 
\end{proposition}
\begin{proof}
Since $F_{\omega}$ is alternating the metric $\g_{o}$ is biinvariant.
Since , \eqref{eq:tistwostep} is satisfied, $ \lie{t}_{\lie{n},\omega}$ is two-step nilpotent, 
and therefore, the biinvariant metric $\g_{o}$ is flat. 
\end{proof}

\paragraph{Split extensions}
We may always choose
 $\omega = 0$.  In this case, the Lie algebra $\lie{t}_{\lie{n},\omega}$
is the semi-direct product  $\lie{n}^* \oplus _{\ad^*}\! \lie{n}$ of $\lie{n}$ 
with its representation space $\lie{n}^*$. 
Thus, a particular interesting and rich class of examples arises as
follows:  

\begin{corollary} \label{cor:semprod}
Let $\lie{n}$ be a two-step nilpotent Lie algebra.  
Then the metric Lie algebra $( \lie{n}^*  \oplus _{\ad^*}\! \lie{n}, \g_{o})$ is flat 
with biinvariant metric. 
\end{corollary}

The corollary shows in particular that the class of
Lie algebras with flat biinvariant metrics is as rich
as the class of two-step nilpotent Lie algebras. 

\paragraph{Twisted extensions with abelian base} 
If $\lie n = \lie a$ is abelian, examples may be constructed using alternating three-forms
on $\lie a$. 

\begin{example}
Let $\lie{n} = \lie{a}_{3}$ be the three-dimensional abelian
Lie algebra. Let $\det \in  \bigwedge^3  \lie{a}_{3}^*$ be 
a non-degenerate alternating three-form, and define 
$$ \omega_{\det}(X,Y) = \det(X,Y, \, \cdot \, ) \; . $$ Then 
$(\lie{t}_{\lie{a}_{3}, \omega_{\det}}, \g_{o})$ is a non-abelian 
biinvariant flat metric Lie algebra.
\end{example}

The following result states that every 
flat metric Lie algebra $(\lie{g}, \g)$  with
biinvariant metric arises in this way: 

\begin{theorem} \label{thm:biinv_struct}
Let $(\lie{g}, \g)$ be a flat metric Lie algebra with
biinvariant metric. Then there exists an abelian Lie algebra $\lie{a}$ 
and an alternating three-form  $F_{\omega} \in  \bigwedge^3  \lie{a}^*$,  
and an  abelian Lie algebra
$\lie z$ such that $(\lie{g}, \g)$ can be written as a direct
product of metric Lie algebras  
$$   (\lie{g}, \g)  \; =  \; (\lie{z}, \g) \,  \oplus \,  ( \lie{t}_{\lie{a},\omega},\g_{o})   \; . $$
\end{theorem}
\begin{proof}  We show that it is  possible to choose $\lie n$ as an abelian
Lie algebra, as in the previous example.  
Note that the commutator subalgebra  $[\lie g , \lie g ]$ of 
$(\lie{g}, \g)$ is an isotropic ideal in $(\lie{g}, \g)$, and its orthogonal
complement with respect to $\g$ is the center $ \lie z(\lg)$ of $\lie g$. This shows that 
we can choose an isotropic subspace $\lie a$ of $(\lg,\g)$, such that 
there is a direct decomposition of subspaces $\lie g = \lz \oplus \lie a \oplus [\lie g , \lie g ]$,
where  $\lie z \subset \lie z(\lg)$. 
\end{proof}

\begin{example} Let $(\lg, \g)$ be a flat biinvariant metric Lie algebra of
dimension six. Then $\lg$ is abelian or $(\lg, \g) = (\lie{t}_{\lie{a}_{3}, \omega_{\det}}, \g_{o})$.
In particular,  in the second case $\lg$, is a semidirect product of the three-dimensional
Heisenberg algebra $\lie h_{3}$ with its coadjoint representation. 
\end{example}

Remark that a version of Theorem \ref{thm:biinv_struct}  already appeared in  \cite[Proposition 7.5]{Wolf_5}. A similar result is stated recently also in \cite{CS}.

\subsection{\' Etale affine representations} \label{sect:etreps}

Let $\rho: G \ra \Aff(n)$ be an  \'etale affine representation, and
let $\rho_{\lg}: \lg \ra \lie{aff}(\bbR^n)$ be the
corresponding Lie algebra representation. 

\subsubsection{The relative invariant for an  \'etale affine representation}
Choose a basis $X_{1}, \ldots, X_{n}$ for $\lie{g}$, and $x \in \bbR^n$.
Let us consider the linear map $\tau_x(\rho): \bbR^n  \rightarrow \bbR^n$, which is 
defined by 
$$   \tau_x(\rho): \;  (\alpha_1, \ldots , \alpha_n) \longmapsto  
         t_x(\;  \sum_{i=1}^n \alpha_i \rho_{\lg}(X_i) \, ) \;.
$$
Here,  $t_x$ is the derivative of the evaluation map at $x$, see \eqref{eq:tx}. 
We easily see that $\tau_x(\rho)$ is an isomorphism 
if and only if the orbit map $o_x: {G} \rightarrow \bbR^n$ is
non-singular, that is, if $\rho$ is \'etale at $x$. \\

If $A \in \Aff(\bbR^n)$, we put $\rho^{A}(g) = A \,  \rho(g) A^{-1}$.
The following formulas are easily verified by direct calculation:

\begin{lemma} Let $A \in \Aff(\bbR^n)$, with $\ell(A) \in \GL(n,\bbR)$ its 
linear part. Let $g \in G$, and $\Ad(g): \lie{g} \ra \lie{g}$ its adjoint
action on the Lie algebra $\lie{g}$ of $G$. We put $B= \rho(g)$. 
Then the following hold:
$$  \det \tau_x(\rho^A) = (\det \ell(A)) \det \tau_{A^{-1}x} (\rho) $$ 
$$  \det \tau_{x}(\rho^B) = (\det \Ad(g)) \det \tau_{x}(\rho)$$
\end{lemma} 

We define the polynomial function 
$$ \delta = \delta(\rho): \bbR^n \ra \bbR  \; , \;  x \mapsto \delta(x) =  \det \tau_x(\rho)$$
which (by the above lemma) satisfies the transformation rule 
\begin{equation} \label{eq:character}
 \delta( \rho(g) x ) \; = \, \det \Ad(g)^{-1} \,   \det \ell(\rho(g)) \; \delta(x)  \; .
\end{equation}

\begin{definition}
We call $\delta = \delta(\rho)$ the \emph{relative invariant} for 
the  \'etale affine representation $\rho$. Its character
is the function $\chi= \chi(\rho): G \ra \bbR^{>0}$, which is defined as
$$ \chi(g) = \det \Ad(g)^{-1} \,  \det \ell(\rho(g)) $$
\end{definition}

Note that an open orbit $\rho(G) x$ of an affine \'etale representation 
$\rho$ is contained in (and it is actually a connected component of) 
the open semi-algebraic subset 
$$   U_{\delta} = \{ x \in \bbR^n \mid  \delta(x) \neq 0 \} \subset \bbR^n \; . 
$$
By a theorem of Whitney \cite{Whitney}, a real semi-algebraic set $U_{\delta}$ has
only finitely many connected components in the standard Euclidean
topology on $\bbR^n$. Thus we have proved the following result: 

\begin{corollary} Let $\rho: G \ra \Aff(n)$ be an affine \'etale representation.
Then $G$ has  finitely many open orbits on $\bbR^n$. Moreover, either
$G$ is transitive or each orbit 
is a connected component of the complement of a hypersurface of 
degree $n$. 
\end{corollary}

\begin{example}
The list of \'etale affine representations of Lie groups $G$, with $\dim G =2$,
is contained in Example \ref{example:domains}. 
\end{example}

The relative invariant $\delta$ may also be interpreted in terms of growth
of the right-invariant volume on $G$, relative to a parallel volume on $\bbA^n$,
see Example \ref{ex:polvol}.

\subsubsection{The relative invariant of a simply transitive representation is constant}
\label{sect:relisc}

The following theorem is actually a direct consequence of
Proposition \ref{prop:ccriterion2}. For transparency, we provide
a direct proof. (For a different proof see \cite{GH_2, Helms}, 
\cite{Kim} for the case of nilpotent groups.)  

\begin{theorem} \label{thm:relinv_strans}
Let $G$ be a Lie group, and $\rho: G \ra \Aff(n)$ be an affine \'etale representation. Then the following conditions are equivalent:
\begin{enumerate}
\item  $\rho$ is simply transitive.
\item  The relative invariant $\delta(\rho)$ is a non-zero constant function. 
\item  $\chi(\rho) \equiv 1$. 
\end{enumerate}
\end{theorem} 
\begin{proof} Clearly, if $\delta(\rho)$ is non-vanishing, 
every orbit of $G$ is open. Hence, there can be only one orbit, and 
$\rho(G)$ is a simply transitive group of affine transformations. 
Conversely, assume that $\rho(G)$ is simply transitive. Then the
relative invariant $\delta(\rho)$ is a nowhere vanishing 
polynomial function. Let $\bG$ denote the Zariski closure of
$\rho(G)$ in $\Aff(\bbR^n)$. Then $\bG$ is a transitive group 
of affine transformations on $\bbR^n$, and $\delta(\rho)$
remains a relative invariant, for some polynomial homomorphism 
$\chi_{\bG}: \bG \ra \bbR^{>0}$, which extends
$\chi$.  Note that also the unipotent radical $\bU$ 
of $\bG$ acts transitively on $\bbR^n$ (even simply transitively, 
see \cite{Auslander}),
and the polynomial character $\chi_{\bG}$ is trivial
on $\bU$. Therefore, by formula \eqref{eq:character}, evaluated on
$\bU$, $\delta(\rho)$ must be constant. 

The latter fact clearly implies that $\chi(\rho)$ is constant. Conversely, if
$\chi(\rho)$ is constant, by  \eqref{eq:character},  $\delta(\rho)$
is constant and non-zero on an open subset of $\bbR^n$. Since, $\delta(\rho)$
is a polynomial it must be constant. 
\end{proof}

\begin{remark} In the case of an \'etale affine representation, the relative invariant 
$\delta(\rho)$ may be identified with the characteristic map $\Phi(\rho): \bbR^n \ra
\bigwedge^n \lg^*$ (as defined in section \ref{sect:charmap}), using an 
isomorphism $\bigwedge^n \lg^* =  \bbR$. See section \ref{sect:red-stab} for further discussion, and, in particular, Proposition \ref{prop:etreps}.
\end{remark}

We call $\rho$ volume preserving if its linear part  satisfies 
$\ell(\rho(G)) \leq \SL(n, \bbR)$. 

\begin{corollary} Let $G$ be unimodular, and $\rho: G \ra \Aff(n)$ an \'etale affine
representation. Then $G$ is simply transitive on $\bbA^n$ 
if and only if $\rho(G)$ is volume preserving. 
\end{corollary}

Here is an immediate application which shows that 
a large class of Lie groups does not admit  \'etale
affine representations. Consider the connected 
abelian group $$ H^1(G) = G/ [G,G] \,  , $$ 
over which every homomorphism of $G$ to an 
abelian group factorises. 

\begin{corollary} If the group $H^1(G)$ is
compact then $G$ does not admit an \'etale 
affine representation.   
\end{corollary}
\begin{proof} By our assumption, for every \'etale
representation $\rho$ of $G$, its character $\chi(\rho)$
must be trivial. Therefore, by the above theorem, $\rho$
must be simply transitive, and $G$ is
simply connected solvable. This shows that, in fact,
$H^1(G,\bbR) = \{ 1 \}$, and since $G$ is solvable,
$G$  must be trivial. A contradiction. 
\end{proof}

\begin{example} The criterion is satisfied, for example, if 
$$ H^1(\lie{g}, \bbR)= \lie{g}/ [\lie{g}, \lie{g}] = \{ 0\} \; .$$
Hence, \emph{no semisimple Lie group admits an \'etale
affine representation}. Also \emph{no compact group}
admits an  \'etale affine representation.
\end{example}

\begin{corollary} If the group $G$ is
compact then $G$ does not admit an \'etale 
affine representation. 
\end{corollary}

\begin{example} Let $G$ be a reductive Lie group with one-dimensional center.
 Then every \'etale affine representation $\rho$ of $G$ is linear. Since $G$ is
unimodular,  $\delta(\rho)$ is 
a non-zero homogeneous  relative invariant polynomial for $G$
of degree $n$, with non-trivial character $\det = \chi(\rho)$. 
Such representations exist, for example, for $G = \GL(n,\bbR)$. 
There do exist other examples as well.  See \cite{BauesT, Burde_2,  
Kimura} for further reference on this problem. 
\end{example}


\subsubsection{The dual tube representation}
\label{sect:dtetalerep}

The dual tube, as discussed in section \ref{sect:Tdomains},
is a useful tool in the construction of Pseudo-Riemannian and
symplectic affine \'etale representations. This is based on the
following simple observation:

\begin{lemma} Let $\rho: G \ra \Aff(n)$ be an \'etale affine representation,
with open orbit $U$. Let $\ell(\rho): G \ra \GL(n)$ denote its linear part. 
Then the associated semi-direct product 
Lie group $G \ltimes_{\ell(\rho)^*} \bbR^n$ 
has an \'etale affine representation on $\bbA^{2n}$, which
has the dual tube domain $\check  T(U)$ as open orbit. 
Moreover, $G \ltimes_{\ell(\rho)^*} \bbR^n$ preserves the natural flat
Pseudo-Riemannian metric of signature $(n,n)$,  and also the natural
symplectic form on $\check T(U)$. 
\end{lemma}
(Here, $\ell(\rho)^*$ denotes the dual representation for $\ell(\rho)$.)\\

\noindent
The lemma is a particular case of Proposition \ref{prop:dtube_psr} applied
to \'etale affine representations.

\section{Affinely homogeneous domains} \label{sect:AHDs}

Let $U \subseteq  \bbA^n$ denote an open subset of $\bbA^n$. We put
$$ \Aff(U) = \{ A \in \Aff(n) \mid A (U) =   U \} $$ 
for the \emph{affine automorphism group}
of $U$. 

\begin{definition}
Let $U \subseteq  \bbA^n$ be a connected  open subset. Then $U$ is
called an \emph{affinely homogeneous domain} if $\Aff(U)$ acts
transitively on $U$.  
\end{definition}
A domain $U$ in $\bbE^{s}$ is called homogeneous, or a 
\emph{Pseudo-Riemannian affine homogeneous domain} if 
its isometry group $\Isom(U) = \Aff(U) \cap E(s)$ acts 
transitively on $U$. 
Similarly, a domain $U$ in symplectic 
affine space $(\bbA^{2k}, \omega)$ is called a \emph{symplectic affine 
homogeneous domain}, if the group of symplectic affine automorphisms
of $U$, $\Aff(U, \omega) = \Aff(U) \cap \Aff(\omega)$ acts transitively. \\

An affine domain $U$ is called \emph{divisible} if there exists a \emph{discrete} subgroup
$$ \Gamma \leq \Aff(U)$$ such that $U/\Gamma$ is compact. 

\subsection{Prehomogeneous affine representations}
Let $G$ be a Lie group and $\rho: G \ra \Aff(n)$ a homomorphism.
The homomorphism $\rho$ is called an affine representation of $G$.
Composition with the natural homomorphism $\Aff(n) \ra \GL(n)$, defines a linear representation
$$ \ell = \ell_{\rho} : G \ra  \GL(n)$$ which is called the linear part of 
$\rho$. 

\begin{definition}
The affine representation $\rho$ is called a \emph{prehomogeneous affine 
representation}, if $\rho(G)$ acts transitively on an open set
 $U \subseteq \bbA^n$. If $G \leq \Aff(n)$ is a subgroup with an open 
orbit in $\bbA^n$ then $G$ is called a \emph{prehomogeneous subgroup}. 
\end{definition}

We let $U_{\rho} 
\subseteq \bbA^n$ denote the union of all open orbits for $G$.\\

\begin{remark}
If $G$ admits a fixed point on $\bbA^n$, then the homomorphism $\rho$ is conjugate by
a translation to the linear representation $\ell_{\rho}$ on the vector space $V= \bbR^n$. 
Such a representation will be called {\em linear}.  A vector space $V$ with
a linear $G$-action, which has an open orbit, is traditionally
called a {\em prehomogeneous vector space} (see \cite{Kimura}). 
\end{remark}

Quite opposite to a linear representation is the situation that 
the affine representation $\rho$ is \emph{transitive} on $\bbA^n$.
In this case, $U_{\rho} = \bbA^n$. \\

The following notion plays a special role. 

\begin{definition} \label{def:etrep}
A prehomogeneous representation $\rho$ on $\bbA^n$ is called an {\em \'etale affine representation} if  
$\dim G = n$.  
\end{definition}





\subsection{Some examples}
There is a wealth of prehomogeneous affine representations $\rho$, and associated
affinely homogeneous domains. We consider
some simple examples. 

\begin{example} \label{example:first}
\begin{enumerate} 
\item 

The vector group $\bbR^n$ acts simply transitively by translations 
on $\bbA^n$. 

\item 
$G = \GL(n, \bbR)$ acts by left- and also by right-multiplication on the vector
space of matrices $\Mat(n \times n, \bbR)$. Both actions have open 
orbit $U = \GL(n, \bbR) \subseteq \bbR^{n^2} = \Mat(n \times n, \bbR)$, and both 
are simply transitive on $U$.

\item 
$\SL(n,\bbR)$ acts transitively on $U= \bbR^{n} - \{ 0 \}$. Also
$\SL(n,\bbR) \times \SL(n,\bbR)$ has an open orbit $U \times U$
in $\bbR^{2n}$.

\end{enumerate}
\end{example}

Note that all the homogeneous domains in the previous example are divisible by
affine transformations. The first two examples admit transitive  \'etale affine groups,
the first one being the role model for a simply transitive affine representation. 
The second and third examples come from prehomogeneous linear representations
of reductive groups, a topic which is studied extensively in the literature, see
\cite{Kimura} and the references therein.  

\paragraph{Two-dimensional affinely homogeneous domains} 
In dimension two there are up to affine equivalence six homogeneous
affine domains. The following list also exhausts the set of all 
two-dimensional \'etale affine representations up to affine equivalence.
(Compare  \cite{NaganoYagi, DuIh_3}.)
Note that each two-dimensional affinely 
homogeneous domain admits a simply transitive group of
affine transformations.

\begin{example}[\'Etale affine representations in dimension 2] \label{example:domains}
 \hspace{1cm} 
\begin{enumerate} 
\item {\bf  (\emph{The plane $U = \bbA^2$})}
$$ \bbR^2   \text{ and } \;  \cU =
 \left\{ 
\begin{matrix}{ccc}  1 & v & u+ {1 \over 2} v^2  \\  0 & 1 & v \\ 
0 &  0 & 1  
\end{matrix} 
\right\} \; \subseteq \Aff(2)
$$
are simply transitive abelian groups of affine transformations.\\

The group 
$$   \cA_{\lambda} = \left\{ \begin{matrix}{ccc}    \exp(t\lambda)  & 0 & s \\
0 & 1 & t \\ 
0 &  0 & 1  \end{matrix} \right\} \; \subseteq \Aff(2) $$
is a simply transitive solvable, non-abelian, group of affine transformations, for
$\lambda \neq 0$. 

\item {\bf  (\emph{The halfspace $U = {\cal H}_{2}$})} 
Let  ${\cal H}_{2}$ be the
halfspace $x >0$. Then \vspace{1ex} \\
$$ \cB = \left\{ \begin{matrix}{ccc}  \exp(t) & 0 & 0 \\  0 & 1 & v \\ 
0 &  0 & 1  \end{matrix} \right\} \; \subseteq  \; \Aff({\cal H}_{2}) = 
\left\{ \begin{matrix}{ccc} \alpha & 0 & 0 \\  z & \beta & v \\  
0 &  0 & 1  \end{matrix} \Mid \alpha > 0 \right\} $$
is an abelian group of affine transformations. The 
half-spaces $(x,y)$, $x>0$ and $x <0$ are open orbits. \\ 

The groups 
$$   \cC_{\lambda, \tau} = \left\{ \begin{matrix}{ccc}    \exp(t\lambda)  & 0 & s \\
0 &  \exp(- t\lambda)  & (\exp(- t\lambda) -1) \tau \\ 
0 &  0 &  1 \end{matrix} \right\} \; \subseteq \Aff(2) $$
are  solvable, non-abelian, groups of affine transformations, with 
two half\-spaces as open orbits, $\lambda \neq 0$. \\

The group of upper triangular matrices
$$ \cC =  \left\{ 
\begin{matrix}{cc} 
  \exp(t\lambda) & b \\  0 &  \exp(- t \lambda)  
\end{matrix} \right\}  \; \subseteq \GL(2)
$$ 
is a solvable, non-abelian,  \emph{linear}  group of transformations
with two half\-spaces as open orbits. \\

\item {\bf  (\emph{The sector $U = {\cal Q}_{2}$})}
$$   \Aff(\mathcal Q_{2})^0 =  \left\{ 
\begin{matrix}{cc} 
   a  & 0 \\  0 &  b  
\end{matrix} \Mid a>0, b> 0  \right\}  \; \subseteq \GL(2)
$$ 
is an abelian,  \emph{linear}  group of transformations,
which has the four open sectors as open  orbits. 

\item {\bf  (\emph{The parabolic domains $U = \mathcal P_{2}^+, \mathcal P_{2}^-$})}
Let $\mathcal P_{2}^+ = \{ (y, x) \mid y > {1 \over 2}  x^2 \} $ be the convex domain 
enclosed by a parabola.
$$   \Aff(\mathcal P^\pm_{2})^0  = 
\left\{ 
\begin{matrix}{ccc}   \lambda^2  & 0 & 0 \\
0 &  \lambda & 0 \\ 
0 &  0 &  1 \end{matrix}
\begin{matrix}{ccc}  1 & v &  {1 \over 2} v^2  \\  0 & 1 & v \\ 
0 &  0 & 1  
\end{matrix} \mid \lambda >0 
\right\} \; \subseteq \Aff(2)
 $$
is solvable, non-abelian, with open orbits $\mathcal P_{2}^+$
and the exterior $\mathcal P_{2}^-$ of a parabola.    

\item  {\bf  (\emph{The plane with a point removed,  $U = \bbA^2 -{0}$}) }
$$   \mathcal E  = \left\{   \exp(t) \begin{matrix}{cc}    \cos \phi & \sin \phi \\
- \sin \phi  &  \cos \phi  \\ 
\end{matrix} \right\} \; \subseteq \GL(2) = \Aut(U) $$
is an abelian \emph{linear} group, with open orbit $\bbR^2- \{0\}$.
\end{enumerate}
\end{example}

We remark:
\begin{corollary} Every two-dimensional affinely homogeneous domain 
$U \subseteq \bbA^2$ 
admits a simple transitive \'etale group of affine automorphisms in $\Aff(U)$.  
\end{corollary}

\begin{corollary} Every two-dimensional affinely homogeneous domain is divisible, 
except for the parabolic domains $\mathcal P^\pm$. 
\end{corollary}
\begin{proof} If $U \subseteq \bbA^2$
admits a simply transitive abelian group $G$ of affine transformations, then 
$U$ is divisible by a lattice $\Gamma \leq G$, $\Gamma \cong \bbZ^2$.
This is the case for all domains, but $\mathcal P^\pm$.  The automorphism
group $G=\Aff(\mathcal P^\pm)$ is a solvable non-unimodular \'etale group of affine 
transformations. Therefore, every $\Gamma$ which divides $\mathcal P^\pm$ is 
a lattice in $G$. However, the existence of a lattice implies that $G$ is 
unimodular (see Appendix C);  a contradiction.
\end{proof}

\begin{corollary}[{\cite{NaganoYagi}}] \label{cor:dev2t}
Let $U \subseteq \bbA^2$ be the development image of a compact affine two-torus. Then $U$ is an affinely homogeneous 
domain which admits an abelian simply transitive group of affine 
transformations. 
\end{corollary}

This also implies: 

\begin{corollary}
Every two-dimensional divisible domain  $U \subseteq \bbA^2$ is convex 
and homogeneous. 
\end{corollary}

There are attempts to classify \'etale affine representations in low dimensions, 
see \cite{NaganoYagi,FurArr, FriedGoldman,Benoist}, for some results. 
In general, it is a difficult problem to decide
which Lie groups $G$ admit an \'etale affine representation.
See section \ref{sect:etreps}, for further discussion.  


\subsection{Characteristic map for a prehomogeneous representation} 
\label{sect:charmap}

Let  $\rho: G \ra \Aff(n)$ be an affine representation, and let $x \in \bbA^n$.
We define
the \emph{orbit map} at $x$ as $$ o_{x}: \, G \ra \bbR^n,  \;  \;  \, 
g  \mapsto \rho(g)  x \;. $$ 
Let $\lg$ be the Lie algebra of $G$, and 
let  $$ \tau_{x} = \tau_{x}(\rho):\;  \lie{g} \ra \bbR^n$$ denote the 
differential of $o_{x}$ at the identity of $G$. We remark:

\begin{enumerate}
\item The kernel of
the linear map  $\tau_{x}$ is the Lie algebra $\lie{h} = \lie{g}_{x}$
of the stabiliser $H= G_{x} \leq G$ of $x$ under the action $\rho$.

\item The orbit $\rho(G)x$ is open 
if and only if  $o_{x}: G \ra U = \rho(G) x$
is a submersion of $G$ onto the open subset $U$ of $\bbA^n$. 
Moreover, \emph{$G$ has an open orbit at $x$
if and only if $\tau_{x}$ has maximal rank, that is, if $\tau_{x}$ is onto}.
\end{enumerate}

Let us fix a parallel volume form $\nu$ on $\bbA^n$. 
We construct the \emph{characteristic map for $\rho$}
$$  \Phi = \Phi_{\rho}:   \; 
\bbA^n  \lra   \bigwedge^n \lie{g}^*  \; ,  \; \, \Phi ( x )  = \tau_{x}^* \nu \; ,  $$
by taking the pull back of $\nu$ along $\tau_{x}$. Its image 
$$ \Phi(\bbA^n) \subseteq \bigwedge^n \lie{g}^* $$ 
is called the \emph{characteristic image} of $\rho$. \\

The map $\Phi$ and
the characteristic image $\Phi(\bbA^n)$ carry fundamental information
about $\rho$ and the associated homogeneous domain $U_{\rho}$.
This theme will be further developed in section \ref{sect:funddiag}.

Note, in particular, that the map  $\Phi$ is a  \emph{polynomial map} in the natural coordinates of 
$\bbA^n$. 
Moreover,  by the above considerations, 
\emph{$G$ has an open orbit at $x$, if and only if the form $\Phi_{x} = \Phi(x)$ 
does not vanish}. \\

\vspace{1ex}
As a first application, we have two basic remarks:  

\begin{proposition} \label{prop:basic}
 \hspace{1cm}  \begin{enumerate}
\item 
Every affinely homogenous domain $U$ is an open semi-algebraic 
subset of $\bbR^n$.
\item 
Let $\rho: G \ra \Aff(n) $ be an affine representation. 
Then $G$ has only finitely many open orbits. 
\end{enumerate}
\end{proposition}
\begin{proof}  Let $G$ be a connected group which is transitive on $U$.
The set  $U_{\rho} = \{ x \mid \Phi_{x} \neq 0 \}$ is a 
Zariski open subset of $\bbA^n$. Its connected components are precisely the 
open orbits for $G$. By \cite{Whitney}[Theorem 4], there are only finitely many
components in the complement of a real algebraic variety. Therefore, $U_{\rho}$
has only finitely many components. Each component is a semi-algebraic subset
of $\bbR^n$. 
\end{proof}

\subsection{The automorphism group of an affine homogeneous domain}
\label{sect:affautU}

Let $U \subseteq \bbA^n$ be an affine homogeneous domain, and  
$ \Aff(U) $ 
its affine automorphism group. 
The following explains the prominent role which algebraic groups
play in the study of affinely homogeneous domains. (For a review on
linear algebraic groups and basic definitions, consult Appendix A.)

\begin{lemma} \label{lemma:autisa}
Let $U \subseteq \bbA^n$ be an affine homogeneous domain. Then
the automorphism group $\Aff(U)$ of  $U$ 
is a finite index subgroup of a real algebraic group. 
In particular, $\Aff(U)$ has finitely many connected components. 
\end{lemma}
\begin{proof}  Let $\tilde{U}  \subseteq \bbA^n$ be the domain, where the Lie group 
$\Aff(U)$ has open orbits. Then $U$ is a connected component of  $\tilde{U}$.  
The complement $W= \bbA^n -  \tilde{U}$ is an algebraic
subset of $\bbA^n$,  by Proposition \ref{prop:basic}.
Since $\Aff( \tilde{U}) = \Aff(W)$, we conclude that
$\Aff( \tilde{U})$ is a real algebraic subgroup of $\Aff(n)$. Since $\tilde{U}$
has only finitely many components (cf.\ Proposition \ref{prop:basic}), 
$\Aff(U)$ has finite index in $\Aff(\tilde{U})$. 
\end{proof}

Let $G \leq \Aff(n)$ be a prehomogeneous Lie subgroup. The real Zariski 
closure $\ac{G}^\bbR \leq \Aff(n)$ of $G$,  is clearly also a prehomogeneous 
real linear algebraic group.

\begin{corollary}  \label{cor:closU}
Let $\rho: G \ra \Aff(n)$ be a prehomogeneous representation. Then 
the Zariski closure $\ac{G}^\bbR \leq \Aff(n)$ %
preserves the maximal domain $U_{\rho}$.  
\end{corollary}

\subsubsection{Automorphisms of complex affine domains}
A connected domain $\cU \subseteq \bbC^n$
is called a \emph{complex affine homogeneous domain}, if its group of
complex affine transformations $$ 
\Aff_{\bbC}(\cU) = \{ A \in \Aff(\bbC^n) \mid A \,  \cU  = \cU \}$$
acts transitively on $\cU$.
We mention:

\begin{lemma} Let $\cU \subseteq \bbC^n$ be a complex affine homogeneous 
domain. Then $\cU$ is a Zariski-open, hence irreducible and Zariski-dense 
connected subset of $\bbC^n$. The automorphism group $\Aff_{\bbC}(\cU) \leq \Aff(\bbC^n)$ 
is a complex linear algebraic group. 
\end{lemma}
\begin{proof} 
The proof follows along the lines of the proof of Lemma \ref{lemma:autisa}.
Note that any open and connected orbit $\cU= \bG^0 x$, where $\bG$ is 
an algebraic subgroup of $\Aff(\bbC^n)$, is Zariski-dense in $\bbC^n$. 
By Appendix \ref{sect:oclosure}, $\bG x$ is also dense in the Euclidean topology of $\bbC^n$. 
Therefore, there is only one such orbit. In particular, $\tilde{\cU} = \cU$ is 
a connected domain, and $\Aff_{\bbC}(\cU) = \Aff_{\bbC}(\tilde \cU)$.  
\end{proof} 

\subsubsection{Centralisers of prehomogeneous representations}

Let $\rho: G \ra \Aff(n)$ be a prehomogeneous representation. We let 
$$ \Z_{\Aff(n)}(G) = \{ A \in \Aff(n) \mid A \,  \rho(g) = \rho(g) A  \text{ , for all $g \in G$ }\}  \,
 \leq \Aff(n) $$ 
denote the centraliser of $G$ in $\Aff(n)$. Since the elements of 
$\Z_{\Aff(n)}(G)$ permute the open orbits
of $G$,  $\Z_{\Aff(n)}(G)$ preserves the maximal set $U_{\rho}$. 
In particular, it follows that  $$ \Z_{\Aff(n)}(G)   \leq \Aff(U_{\rho}) \; .  $$

Note that $\Z_{\Aff(n)}(G)$ 
is an algebraic subgroup of $\Aff(n)$.
Moreover, we have  the following fact:

\begin{lemma}   \label{lemma:centisass} 
The Lie algebra $\lie{z}_{\rho}(G)$ of $\Z_{\Aff(n)}(G)$ forms an
associative subalgebra of $\aff(n)$.  
\end{lemma} 
\begin{proof} Let $\varphi \in \aff(n)$.
Clearly,  $\varphi \in \lie{z}_{\rho}(G)$ if and only if $A \varphi A^{-1} = \varphi$,
for all $A \in \rho(G)$. Thus, if $\varphi, \psi \in \lie{z}_{\rho}(G)$ then
$A (\varphi  \psi) A^{-1} = \varphi  \psi \in \lie{z}_{\rho}(G)$. 
\end{proof}

The following lemma is useful in this context: 

\begin{lemma} \label{lemma:nilopsalg}
Let $\cA \leq \Mat(n \times n)$ be an associative subalgebra
of linear operators, such that $\trace \phi = 0$, for all $\phi \in \cA$. Then
every element $\phi \in \cA$ is nilpotent. Moreover, the Lie 
algebra which belongs to $\cA$ (by taking commutators in $\cA$) 
is nilpotent. 
\end{lemma}

\begin{corollary} \label{cor:centisu}
Let $\rho:G \ra \Aff(n)$ be a prehomogeneous representation such that  
every element of the  centraliser algebra $\lie{z}_{\rho}(G)$ contains 
only elements of trace zero. Then the centraliser $\Z_{\Aff(n)}(G)$ 
contains only unipotent elements, and, in particular, it is a connected 
nilpotent Lie group. 
\end{corollary}
\begin{proof}   
Let $A \in \Z_{\Aff(n)}(G)$. Viewing  $A$ as an element of  $\GL(n+1, \bbR)$, 
we obtain  $A-E_{n+1} \in \lie{z}_{\rho}(G)$. Thus, $A-E_{n+1}$ is nilpotent,
by Lemma \ref{lemma:nilopsalg}.
This implies that the real algebraic group $\Z_{\Aff(n)}(G)$ has
only unipotent elements.
\end{proof}

This situation occurs, for example, if $\rho$ is transitive:  

\begin{proposition} \label{prop:centisu}
Let $\rho$ be a transitive affine representation. 
Then the centraliser of $\rho(G)$ in $\Aff(n)$ is a connected 
unipotent group of dimension $\leq n$.  
\end{proposition}
\begin{proof} If $G \leq \Aff(n)$ is a transitive subgroup, $\Z_{\Aff(n)}(G)$ 
acts without fixed points on $\bbA^n$. The centraliser 
$\Z_{\Aff(n)}(G)$ is an algebraic subgroup of $\Aff(n)$. 
Its reductive part $T$ has a fixed point, by Lemma \ref{lemma:fixedpoint}. 
Therefore, $T$ must be trivial.  
\end{proof}

We shall show below (see Proposition \ref{prop:cstep2} and Proposition 
\ref{prop:scisunip}) that  the centralisers of Pseudo-Riemannian 
and symplectic prehomogeneous groups in the group of isometries of $\bbE^{s}$,
respectively symplectic affine transformations of $(\bbA^n,\omega)$ are
unipotent groups. Such results have immediate geometric implications.
As the following example shows, the latter phenomena
depend on the orthogonal and symplectic structure:  

\begin{example} \label{example:cnonunipotent}
$G= \SL(n,\bbR) \times \SL(n,\bbR)$ is a volume preserving prehomogeneous subgroup of $\Aff(\bbR^n \oplus \bbR^n)$, as in 2.\ Example
\ref{example:first}. For $\lambda \in \bbR$, define $t_{\lambda} (u, v) =
(\lambda u, \lambda^{-1} v)$. Then $T = \{ t_{\lambda} \mid \lambda >0 \} \leq \SL(\bbR^n 
\oplus \bbR^n)$ is an abelian group of volume preserving linear maps, 
consisting of semisimple elements, and $T$
is contained in the centraliser of $G$.
\end{example}

Note also that the volume preserving group 
$T$ in the example leaves invariant the natural symplectic form $\omega$ 
on $\bbR^n \oplus \bbR^n$.

\subsection{Tube like domains} \label{sect:Tdomains}
Affinely homogeneous domains play an important role in the theory of Hermitian 
symmetric spaces, by the \emph{tube construction}, which itself is an important source for the construction of complex homogeneous domains from affine homogeneous
domains. 
An analogous dual construction plays a role for the theory of  Pseudo-\-Rie\-mannian
and symplectically affine homogeneous domains. We develop its basic 
properties and draw some consequences concerning the classification of
Pseudo-Riemannian and symplectically affine homogenous domains. 

\subsubsection{Complex tube domains}
The classical construction of \emph{tube domains} goes as follows. 
Let $U \subseteq \bbR^n$ be a connected open set. Then 
$$  T(U) = U + i \bbR^n  \, \subseteq \bbC^n $$ is called the
\emph{complex tube domain} associated to $U$. 
Every $A \in \Aff(\bbR^n)$ extends to a complex affine transformation 
$A_\bbC \in \Aff(\bbC^n)$,  where $$ A_\bbC (x+iy ) =  A x + i \,  \ell(A) y \; . $$ 
If $A \in \Aff(U)$ is an affine automorphism of $U$ then $A_\bbC \in \Aff_{\bbC}(T(U))$. 
Also the group of purely imaginary translations
$i \,  \bbR^n \subset \bbC^n$ acts on $T(U)$, 
and this group is normalised by the transformations $A_{\bbC}$. 
In particular, if $U$ is a homogeneous affine domain then 
$T(U)$ admits a transitive group of complex affine 
transformations. That is, $T(U)$ 
is a \emph{complex affine homogeneous domain}.

\begin{example}
Riemannian Hermitian symmetric spaces (cf.\ \cite{Helgason}[Ch.\ VIII]) 
are obtained by using the 
Bergman metric on a complex homogeneous domain $T(U)$, 
where $U \subset \bbR^n$ is a convex non-degenerate homogeneous
selfdual cone. 
\end{example}
Complex domains  $T(U)$, 
where $U$ is a convex non-degenerate 
cone are also called \emph{Siegel-domains} of the first kind.
Here a convex cone is called non-degenerate if it does not contain
any straight line. See \cite{Koszul_1, Koszul_2} and 
Matsushimas expository paper \cite{Matsushima} for 
results on homogeneous tube domains over convex non-degenerate cones, 
and their relation with \emph{bounded domains} in $\bbC^n$. 
See  \cite{FaGi} for a generalisation in the context of 
Pseudo-Hermitian symmetric spaces.

\subsubsection{Dual tube domains} \label{sect:dtc}
A process which is dual to the construction of the tube domain 
$T(U) \subseteq \bbC^n$, $A \mapsto A_\bbC$, allows to construct 
Pseudo-Riemannian 
and symplectically affine homogeneous 
domains from arbitrary homogeneous affine domains.
Perhaps surprisingly, this construction suggests that 
\emph{the determination
of homogeneous Pseudo-Riemannian and symplectic affine
domains is at least as complicated as the classification of all
affinely homogeneous domains.} 

\paragraph{The dual tube domain $\check{T}(U)$}
Let $V= \bbR^n$ denote a real vector space, and $V^*$ the
dual space of $V$. If $\varphi$ is a linear map 
of $V$, we let $\varphi^*$ denote the dual (transposed)
map of $V^*$. If  $U \subseteq V$ is a connected open subset
we call 
$$ \check{T}(U) = U \times V^* \subseteq V \oplus V^* $$ 
the \emph{dual tube domain over $U$}. \\ 

Let $A \in \Aff(\bbR^n)$. We extend $A$ to an affine map
$\check{A}  \in \Aff(V \oplus V^*)$, by declaring
$$ \check{A} \, (u + \lambda ) =  A u + ( \ell(A)^{-1})^* \lambda \;  . $$ 
If $A$ is an affine automorphism of $U$ then $\check{A}$
preserves $\check{T}(U)$. Note also that the group of translations 
$V^* \subseteq \Aff(\check{T}(U))$ is normalised by the transformation 
$\check{A}$. 
In particular, if $U$ is a homogeneous affine domain then 
$\check{T}(U)$ is a homogeneous affine domain, as well. 

\paragraph{Natural Pseudo-Riemannnian metric on $\check{T}(U)$}
The vector space $V \oplus V^{*}$ admits a naturally defined 
scalar product $\check{\g}$ of signature $(n,n)$. This is defined by 
 $$ \check{ < u + \lambda ,  u' + \lambda' > }  =  \lambda(u') + \lambda'(u) \; . $$
The scalar product $\check\g$ is evidently preserved by $\check A$, $A \in \Aff(n)$. 
Thus: 
\begin{proposition} The map $A \mapsto\check{A}$ defines a faithful 
homomorphism $$ \Aff(U) \ra \Isom(\check{T}(U), \check{\g}) \; . $$ 
\end{proposition}
 
In particular, we remark:     

\begin{proposition}   \label{prop:dtube_psr}
Let $U \subseteq \bbA^n$ be an affine homogeneous domain. 
Then the tube domain $(\check{T}(U), \check \g)$ is a Pseudo-Riemannian 
homogeneous domain of signature $(n,n)$.   
\end{proposition}
\begin{proof} In fact let, $G \leq \Aff(n)$ be a subgroup which
acts transitively on $U$. Then $V^* \rtimes \check{G} \leq E(n,n)$
acts transitively on $\check{T}(U) = U \times V^*$.
\end{proof}

Note that $V^*$ is a maximal 
isotropic subspace in $\bbR^{n,n} = V \oplus V^*$. 
The following geometric characterisation of Pseudo-Riemannian domains of 
tube type is straightforward: 

\begin{lemma} Let  ${\cal D} \subseteq \bbE^{n,n}$ be 
a  homogeneous Pseudo-\-Rie\-mann\-ian domain. If 
$\Isom(\cD)$ contains a maximal isotropic subgroup  
of translations then 
$$ {\cal D} = \check{T}(U) \; ,  $$
for some affinely homogeneous domain $U \subseteq \bbA^n$.
\end{lemma}

The following notion is related: 

\begin{definition} A domain $U \subseteq \bbE^s$ is called \emph{translationally 
isotropic} if the group of translations $t(U)$ in $\Isom(U)$ satisfies $t(U)^\perp \subseteq t(U)$. 
\end{definition} 

Every tube domain $(\check{T}(U), \check \g) \subseteq \bbE^{n,n}$ is 
\emph{translationally isotropic}, since $V^* \subseteq t(U)$ is a maximal 
isotropic subspace, and thus $t(U)^\perp \subseteq V^*$.  

\begin{proposition}[{\cite{DuIh_2}[Theorem 5.3]}]  \label{prop:tisotropic}
Every translationally isotropic homogeneous Pseudo-\-Rie\-mann\-ian
domain ${\cal D} \subseteq \bbE^{m,l}$ is of the form 
$$ {\cal D} = \check{T}(U) \times \bbE^{m-k, l-k} \; , $$
for some affinely homogeneous domain  $U \subseteq \bbA^k$, 
$l,m  \geq k$.
\end{proposition} 

\begin{question} By \cite{DuIh_1}, the Lorentzian homogeneous 
domains are $\bbE^{n-1,1}$, and the translationally isotropic 
domains $\check T(\bbR^{>0}) \times \bbE^{n-2}$. 
Is every homogeneous Pseudo-\-Rie\-mann\-ian domain 
translationally isotropic? 
\end{question} 

The question suggests, in particular, that homogeneous
Pseudo-Riemannian domains admit large groups of
translations.  

\paragraph{Natural symplectic geometry on $\check{T}(U)$}
On the vector space $V \oplus V^{*}$ there is also a naturally defined 
symplectic form $\check\omega$, which is given by 
$$  \check\omega ( \, u + \lambda ,  u' + \lambda' \, )  =  \lambda(u') - \lambda'(u) \; . $$
The symplectic form $\check \omega$ is evidently preserved by 
$\check A$, $A \in \Aff(n)$. Most considerations for the Pseudo-Riemannian
case carry over analogously to the symplectic case. We just 
remark: 

\begin{proposition} The map $A \mapsto\check{A}$ defines a faithful 
homomorphism $$ \Aff(U) \ra \Aff(\check{T}(U), \check{\omega}) \; . $$ 
\end{proposition}

Thus:

\begin{corollary}   \label{cor:dtube_s}
Let $U$ be an affine homogeneous domain. 
Then $(\check{T}(U), \check \omega)$ is a symplectic 
affine homogeneous domain.  
\end{corollary}

Moreover, we mention:     

\begin{lemma} Let  ${\cal D} \subseteq \bbA^{2n}$ be 
a  symplectic affine homogeneous domain. If 
$\Aff(\check{T}(U), \check{\omega})$ contains a 
Lagrangian subgroup  
of translations then $$ {\cal D} = \check{T}(U) \; ,  $$
for some affinely homogeneous domain $U \subseteq \bbA^n$.
\end{lemma} 

As in the Pseudo-Riemannian case, one might wonder which
role  translationally isotropic symplectically affine homogeneous 
domains play in the classification of all symplectically affine 
homogeneous domains.  However, here we have

\begin{example} \label{example:snontube}
 $(U = \bbR^2- \{ 0 \}, \; \omega = dx \wedge dy)$ is a symplectically 
affine homogeneous domain, and $\Aff(U,\omega) = \SL(2,\bbR)$.
In particular,  $U$ is not translationally isotropic. 
\end{example}

\paragraph{Para-K\"ahler geometry of $\check{T}(U)$}
We remark that the natural geometry on $\check{T}(U)$, 
which is preserved by the transformations $\check A$, is actually
determined by the symplectic form $\check{\omega}$ together
with the metric $\check \g$. 
It may be expressed in a manner which
is analogous to classical K\"ahler geometry as follows.\\

We define a linear operator $\check{J} \in \GL(V \oplus V^*)$, which is
skew with respect to $\check{\omega}$ and satisfies $\check{J}^2 = \id$
by the relation 
$$ \check{g} ( \cdot , \cdot ) = \check{\omega}(\check{J} \cdot , \cdot ) \; .$$
It is computed by the formula
$$ \check{J} (u+\lambda) = ( -u + \lambda) \; . $$ 
Note, in particular, that the transformations $\check{A}$ are $\check{J}$-linear, that is, 
$$ \ell( \check{A})  \check{J}  = \check{J} \ell (\check{A})  \; . $$ 
A geometric structure $(\check{\omega}, \check{J})$ of this kind, 
is thus a sort of analogue of complex K\"ahler geometry. Such structures have 
attracted recent interest of mathematicians and also physicists under
the names of \emph{Para-K\"ahler} or \emph{Bi-Lagrangian}
geometry. See, for example, \cite[\S 5.2]{Bryant} and, in particular,
the work of Kaneyuki et al.\ on Para-K\"ahler symmetric spaces 
\cite{Ka1,Ka2,KaKo}.

%


%




%



\subsection{Pseudo-Riemannian affine homogeneous domains} \label{sect:psrdomains}

As follows from Proposition \ref{prop:dtube_psr}, there are as many 
Pseudo-Riemannian affine homogeneous domains as there are
homogeneous affine domains. In particular, the classification
of  Pseudo-Riemannian affine homogeneous domains might be 
a quite untractable problem. However, for small index $s$, 
the possible types of homogeneous domains in $\bbE^s$
seem rather restricted. \\

Recall that a homogeneous Riemannian manifold is geodesically complete, since 
it is complete with respect to the Riemannian distance. Thus: 

\begin{example}
There is only one Euclidean homogeneous domain $U \subseteq \bbE^n$, 
namely,  Euclidean space $\bbE^n$ itself. 
\end{example}

As follows from Example \ref{example:domains}, we note for the case $n=2$:

\begin{example} The two-dimensional Pseudo-Riemannian
homogeneous domains are $\bbE^2$, $\bbE^{1,1}$, and the half-space 
${\cal H}_{2} = T(\bbR^{>0})  \subset \bbE^{1,1}$.
\end{example}

In the Lorentzian case,  we have the following result:

\begin{proposition}[\cite{DuIh_1}] 
Let $U \subseteq \bbE^{n-1,1}$ be a homogeneous 
Lorentzian domain. Then $U = \bbE^{n-1,1}$ or $U = {\cal H}_{2}  \times 
\bbE^{n-2} \subset \bbE^{n-1,1}$. 
\end{proposition}
In particular, every Lorentzian affine homogeneous domain
is translationally isotropic. \\

We call a Pseudo-Riemannian domain \emph{irreducible} if it
does not admit a decomposition as in Proposition \ref{prop:tisotropic}. 

\begin{example} The (translationally isotropic) irreducible 
domains in $\bbE^{2,2}$ are the dual tube 
domains $(\check{T}(U), \check{\g})$, where 
$U = \bbA^2 - \{0 \}$, ${\cal H}_{2}$, ${\cal Q}_{2}$,
and $\mathcal P^{\pm}_{2}$. (Compare Example \ref{example:domains}.)
\end{example}

It seems unknown, whether (up to products) the latter list exhausts all
homogeneous domains in  $\bbE^{n-2,2}$.

\subsubsection{Centralisers of prehomogeneous groups of isometries}
Let $\rho: G \ra E({s})$ be a prehomogeneous group of isometries.
We consider the subgroup of isometries in $E({s})$, which 
centralise $\rho(G)$.  

\begin{corollary} 
\label{cor:psrcisunip}
Let $U \subseteq \bbE^s$ be a Pseudo-Riemannian affine homogeneous 
domain. Then the centraliser of any prehomogeneous subgroup 
$G \leq \Isom(U)$ in the group $E(s)$
is a connected unipotent group. Moreover, it is a nilpotent group 
of nilpotency class at  most two. 
\end{corollary} 
\begin{proof} By Corollary \ref{cor:Gcentralg}, the Lie algebra of 
$\Z_{\Aff(n)}(G) \cap E({s})$ consists of nilpotent elements 
of  $\aff(n)$, and it is a nilpotent Lie algebra 
of nilpotentency class at most two. 
Therefore, the identity
component $(\Z_{\Aff(n)}(G) \cap E({s}))^0$ is a unipotent group, and
it is nilpotent of nilpotency class at most two.  
We remark next that $\Z_{\Aff(n)}(G) \cap E({s})$ is connected.
(To prove that  $\Z_{\Aff(n)}(G) \cap E({s})$ is connected, we can 
argue as in the proof of Proposition \ref{prop:scisunip}.)
\end{proof} 

The following observation is due to Wolf \cite{Wolf} (for transitive groups), 
see also \cite{DuIh_2} for the general case:  

\begin{proposition} \label{prop:cstep2}
The elements in $g \in \Z_{\Aff(n)}(G) \cap E({s})$ 
are unipotent, and satisfy $(g-E_{n+1})^2 = 0$. 
\end{proposition}
\begin{proof} 
By the above, $\Z_{\Aff(n)}(G) \cap E({s})$ is a
unipotent linear algebraic group. Using the exponential representation
of $g \in  \Z_{\Aff(n)}(G) \cap E({s})$, Proposition \ref{prop:Gcentrvf}
implies that $g-E_{n+1}$ represents a Killing
vector field in $\aff(n)$, and $(g-E_{n+1})^2 = 0$. 
\end{proof}

The previous two results play a central role in the determination of all flat
Pseudo-Riemannian homogeneous manifolds, which are complete.
(See section \ref{sect:homhol_psr} for further discussion.) Another application 
was observed in \cite{DuIh_2}:

\begin{corollary} \label{cor:niltransp}
Every nilpotent prehomogeneous group $N$ of isometries 
of $\bbE^{s}$ is transitive on $\bbE^{s}$.
\end{corollary}
\begin{proof} The set of semisimple elements in a nilpotent algebraic group
forms a central subgroup, see \cite{Borel}. Thus, by Proposition \ref{prop:cstep2}, 
$N$ must 
be a unipotent group. By Corollary, \ref{cor:uistrans}, $N$ acts 
transitively on $\bbA^n$. 
\end{proof}

The latter result generalises to nilpotent groups of 
volume preserving affine transformations, see  Corollary \ref{cor:niltransitive}. 
For the proof, other methods are required.

\subsection{Symplectic affine homogeneous domains} \label{sect:sympdomains}
As for Pseudo-Riemannian domains, the tube construction 
produces a symplectic affinely homogeneous domain 
$(\check T(U), \check{\omega}) \subseteq \bbA^{2k}$
from each affinely homogeneous domain $U \subseteq \bbA^k$.
However, we already noted in Example \ref{example:snontube}
that this construction does not exhaust
the set of  symplectic affinely homogeneous domains:

\begin{example} The two-dimensional symplectic affine homogeneous domains
of tube type are $\bbA^2$ and the halfspace $\mathcal H_{2} = \check T(\bbR^{>0})$.
\end{example}

On the other hand, from Example \ref{example:domains}, we note:

\begin{example} The two-dimensional symplectic affinely homogeneous
domains are $\bbA^2- \{ 0 \}$, and the tube type domains $\bbA^2$, $\mathcal H_{2}$.
\end{example}

Here is another similiarity to the Pseudo-Riemannian
situation (compare Corollary \ref{cor:psrcisunip}):  

\begin{proposition} \label{prop:scisunip}
Let $U \subseteq \bbA^{2k}$ be a symplectic affine homogeneous 
domain. Then the centraliser of any prehomogeneous subgroup 
$G \leq \Aff(U,\omega_{k})$ in the group of
symplectic affine transformations of $\bbA^{2k}$ is a connected
unipotent group (and it is also abelian). 
\end{proposition}
(The connected component of the symplectic centraliser of a prehomogeneous group of symplectic transformations is abelian, 
by Proposition \ref{prop:symplecticentraliser}.)\\

Hence, we get:  

\begin{corollary}  \label{cor:niltranss}
Every nilpotent prehomogeneous subgroup of  
$\Aff(\omega_{k})$ is transitive on $\bbA^{2k}$.
\end{corollary}

For the proof of Proposition \ref{prop:scisunip},
let us first formulate a lemma:

\begin{lemma} \label{lemma:ssdecomp}
Let $T$ be an abelian subgroup of semisimple
elements in $\SP(V,\omega)$, and $G \leq \SP(V,\omega)$ a subgroup
which centralises $T$.  Let $V^T$ denote the 
subspace of invariants for $T$. Then there exists a
decomposition $V = V^T \oplus (W_{1} \oplus W_{2})$, which is invariant by 
$G$ and $T$, such that
$W_{1}$ and $W_{2}$ are isotropic subspaces, and the 
restriction of $\omega$ to the subspaces $V^T$ and $W_{1} \oplus W_{2}$ 
is non-degenerate. Moreover, $V^T$ and $W_{1} \oplus W_{2}$ are
orthogonal with respect to $\omega$. 
\end{lemma}

The lemma is a direct consequence of the decomposition of $V$
in eigenspaces with respect to $T$. \\

\begin{prf}{Proof of Proposition \ref{prop:scisunip}.}  Let $G \leq \Aff(U,\omega_{k})$ be a prehomogeneous subgroup with orbit $U$. For $g \in \Aff(n)$, $t(g) \in V= \bbR^n$ denotes the translation part of $g$. Let $T$ denote a subgroup of semisimple elements in 
$\Z_{\Aff(\omega_{k})}(G)$. After a change of origin, we may assume that $T$ is linear,
and, in particular, $T \leq \SP(\omega_{k})$. Since $T$ centralises $G$, $t(g)\in V^T$, for all $g \in G$,
and the subspace $V^T$ is invariant by  $\ell(G)$.
Let $x = x_{0} + x_{1} \in U$, where $x_{0} \in V^T$ and $x_{1} \in W_{1} \oplus W_{2}$,
as in Lemma \ref{lemma:ssdecomp}. Let $\pi_{2}: V \ra W_{1} \oplus W_{2}$ be the
corresponding projection operator. We conclude that $\pi_{2}(U) = 
\pi_{2} (Gx ) = \pi_{2}(G x_{1}) =
\pi_{2} (\ell(G) x_{1}) = \ell(G) x_{1}$. Since $\ell(G)$ also preserves a Lagrange decomposition of  $W_{1} \oplus W_{2}$, it cannot have an open orbit in 
$W_{1} \oplus W_{2}$. (In fact, the dual pairing of $W_{1}$ and $W_{2}$ with
respect to $\omega$ produces a non-trivial $\ell(G)$-invariant polynomial.)
Unless $T= \{1\}$, this is a contradiction, since $\pi_{2}(U)$ is open in 
$W_{1} \oplus W_{2}$. 
\end{prf} 



\section{A criterion for transitivity of prehomogeneous representations} 
\label{sect:crittrans}

A distinguished  class of prehomogeneous affine representations
is formed by those representations, which are  transitive on affine
space $\bbA^n$. Our main concern in this section 
will be to formulate a criterion, which ensures that a 
given prehomogeneous affine representation $\rho: G \ra \Aff(n)$ is 
transitive on $\bbA^n$. \\

In section \ref{sect:trans1}, we recall how
transitivity depends on the unipotent radical of the Zariski closure. 
A closer analysis of the characteristic map $\Phi(\rho)$, for
a prehomogeneous representation $\rho$,  is given in section \ref{sect:funddiag}.
This allows to interpret transitivity for $\rho$ in terms of the characteristic image 
$$  \Phi(\rho)(\bbA^n) \, \subseteq \, \bigwedge^n \lie{g}^*   \; .  $$ 
An immediate application (see Corollary \ref{cor:niltransitive}) 
is the following result: 

\begin{theorem} Let $U$ be an affine homogeneous domain, which admits a nilpotent transitive subgroup of volume preserving affine transformations. Then $U = \bbA^n$.  
\end{theorem} 

In particular, for nilpotent prehomogeneous groups, transitivity is equivalent to
volume preservation. Various special cases of this result have been treated in the literature before. 
For the case of \'etale affine representations see \cite{Kim,FGH}, and for 
Pseudo-Riemannian prehomogeneous affine representations of nilpotent groups
see \cite[Corollary \ref{cor:niltransp}]{DuIh_2}. 

\subsection{Transitivity for prehomogeneous groups} \label{sect:trans1}
As a first remark, we note that transitivity may be recognised by looking
at the real Zariski closure:  

\begin{corollary}  \label{cor:clostrans}
Let $G \leq \Aff(n)$ be a prehomogeneous subgroup. Then  
$G$ is transitive on $\bbA^n$ if and only if its Zariski closure 
$\ac{G}^\bbR \leq \Aff(n)$ is transitive.  
\end{corollary}

The corollary is a direct consequence of Corollary \ref{cor:closU}. \\

Transitivity for prehomogeneous algebraic groups is determined by 
the unipotent radical, as follows:
\begin{proposition} \label{prop:unitrans}
Let $G \leq \Aff(n)$ be a prehomogeneous real algebraic group. Then the
following are equivalent:
\begin{enumerate}
\item $G$ is transitive on $\bbA^n$. 
\item The unipotent radical $U(G)$ is transitive on $\bbA^n$.
\end{enumerate}
\end{proposition}
\begin{proof}  By section \ref{sect:alggroups}, $G = U(G) H$, where $H$ 
is reductive. Moreover, $H$ has a fixed point on $\bbA^n$, by Lemma 
\ref{lemma:fixedpoint}. Thus, if $G$ is transitive, $U(G)$ must be transitive as well. 
\end{proof} 

Let $G \leq \Aff(n)$ be a prehomogeneous subgroup.
Let  $\bG = \ac{G}^Z  \leq \Aff(\bbC^n)$ denote the complex 
Zariski closure of $G$. Then $\bG$ acts prehomogeneously on 
complex affine space 
$\bbC^n$.  

\begin{corollary} Let $G \leq \Aff(n)$ be prehomogeneous. 
Then  $G$ is transitive on $\bbA^n$ if and only if 
its complex Zariski closure $\bG$ is transitive on $\bbC^n$. 
\end{corollary} 
\begin{proof} If $G$ is transitive, so is the unipotent radical $U$
of its real closure  $\ac{G}^\bbR$. The Zariski closure of $U$, $\bU \leq \bG$ is prehomogeneous on $\bbC^n$ and unipotent. By the closed orbit property, 
$\bU$ must be transitive. Therefore,  $\bG$ is transitive.

Conversely, assume $\bG$ is transitive. Note that Proposition \ref{prop:unitrans} holds analogously for complex algebraic actions on $\bbC^n$. We deduce that 
$\bU$ acts transitively on $\bbC^n$. 
 Note that $U=  \bU_{\bbR}$ is
Zariski-dense in $\bU$ (see \cite{Borel}). In particular, 
$\dim_{\bbR} U = \dim \bU$. Thus $U$, and, in particular, also $\ac{G}^\bbR$
acts transitively on $\bbA^n$. By Corollary \ref{cor:clostrans}, $G$
is transitive,  as well. 
\end{proof} 

\paragraph{Transitivity of nilpotent groups}
We remark that unipotent prehomogeneous groups
are always transitive: 

\begin{corollary} \label{cor:uistrans} 
Let $N \leq \Aff(n)$ be a prehomogeneous group which
is nilpotent. Then $N$ is transitive on $\bbA^n$ if and only if
$N$ is unipotent. 
\end{corollary}
\begin{proof} In fact, if $N$ is unipotent every orbit on $\bbA^n$ is closed (cf.\
Appendix A, Proposition \ref{prop:unip_closed}). 
The converse is a consequence of Proposition \ref{prop:centisu}.
\end{proof}

\subsection{The fundamental diagram} \label{sect:funddiag}
We let $N= N_{\Aff(n)}(\rho(G)) \leq \Aff(n)$ 
denote the normaliser of $G$ in $\Aff(n)$. Then $N$ acts by
conjugation on the image $\rho(G)$. Assuming that $\rho$ is 
faithful, this gives also rise to an action $C: N \ra \Aut(G)$, where, for $A \in N$,
$$ C(A): G \ra G \; \;  $$ is defined by the relation
$$ \rho \left( C(A) (g) \right)  \; =  \; A \, \rho(g) \, A^{-1} . $$  
Furthermore, we let $c: N \ra \Aut(\lie{g})$ denote the induced representation  on 
the Lie algebra $\lie{g}$. \\

We pick up the notation used in section  \ref{sect:charmap}. 
Recall, in particular, that $\tau_{x}: \lg \ra \bbR^n$ denotes the differential of
the orbit map $o_{x}: G \ra \bbA^n$ at $x \in \bbA^n$. 
For all $x \in \bbA^n$, $A \in N$, the following
diagram is commutative: 
\begin{equation} \label{eq:fundiag}
\begin{CD} 
\lie{g}  @> \tau_{x} >>  \bbR^n \\
 @V {c(A)}VV     @VV {\ell(A)}V  \\
\lie{g}  @>{ \tau_{A x}}>>  \bbR^n
\end{CD}
 \; \;  \; \;  \; \;  .
\end{equation}
Now consider the 
characteristic map 
$$ \Phi= \Phi({\rho}): \bbA^n \ra \bigwedge^n \lg^* \;  ,\; \, \Phi ( x )  = \tau_{x}^* \nu \;  $$ 
for $\rho$. The  commutative diagram \eqref{eq:fundiag} 
implies  the following relation, which holds for all $A \in \N_{\Aff(n)}(\rho(G))$: 
\begin{equation}  \label{eq:fundform}
\Phi_{A x}  \; = \;  \det \ell(A) \,   \left( (c(A)^{-1})^* \Phi_{ x}  \right)   \; . 
\end{equation}

By section \ref{sect:charmap}, the following is evident: 

\begin{proposition} \label{prop:ccriterion1}
Let $\rho: G \ra \Aff(n)$ be a prehomogeneous
affine representation. Then the following conditions are equivalent:
\begin{enumerate}
\item[i)] $G$ acts transitively on $\bbA^n$. 
\item[ii)] $0$ is not contained in the characteristic image $\Phi(\bbA^n)$.
\end{enumerate}
\end{proposition}

For an algebraic  representation, we note
the following refinement:

\begin{proposition}  \label{prop:ccriterion2}
Let $\rho: G \ra \Aff(G)$ be a prehomogeneous algebraic representation
with maximal domain $U_{\rho}$. 
Then the following are equivalent: 
\begin{enumerate}
\item[i)] $G$ acts transitively on $\bbA^n$. 
\item[ii)] $0$ is not contained in the closure of the characteristic image $\Phi(\bbA^n)$.
\item[iii)]  
The image $\Phi(U_{\rho}) \subseteq \bigwedge^n \lie{g}^*$  is Zariski-closed.
\item[iv)]  The image $\Phi(U_{\rho}) \subseteq \bigwedge^n \lie{g}^*$ is closed in the Euclidean
topology.  
\end{enumerate}
\end{proposition}
\begin{proof} 
We show first, if $G$ acts transitively on $\bbA^n$ then $\Phi(\bbA^n)$
is Zariski closed in $\bigwedge^n \lie{g}^*$. In fact, if $G$ acts transitively
then its unipotent radical $U(G)$ acts transitively as well. Thus, using
relation \eqref{eq:fundform}, we deduce that 
$$ \Phi(\bbA^n) = \Phi(U(G) x) = c(U(G)) \Phi_{x} $$   is the 
orbit of a unipotent linear algebraic group. Therefore, 
$\Phi(\bbA^n)$ is Zariski closed.  Hence, i) implies ii) 
and also iii). In particular, by Proposition \ref{prop:ccriterion1}, 
i) is equivalent to ii). 

Now, assume condition iii). Then $\Phi(U_{\rho})$ is Zariski-closed. 
Note that the open set $U_{\rho}$ is
Zariski-dense in $\bbA^n$. Since $\Phi$ is a polynomial map (a morphism of varieties),
it preserves the Zariski closure.
Thus, $$ \Phi(\bbA^n) = \Phi( \ac{U_{\rho}}^\bbR ) \subseteq \;
\ac{\Phi(U_{\rho})}^\bbR = \Phi(U_{\rho}) \; .  $$ 
Since $\Phi$ is non-zero on $U_{\rho}$, $0$ is not contained in $\Phi(\bbA^n)$. 
This shows that $U_{\rho} = \bbA^n$,  and, hence, $G$ acts transitively. 
We proved that i) is equivalent to iii). 

Since $\Phi(U_{\rho})$ is a finite union of orbits,  $\Phi(U_{\rho})$ is Zariski closed 
if and only if it is closed in the Euclidean topology. (See the remarks in
section \ref{sect:alggroups}.) Hence, iii) is equivalent to iv). 
\end{proof}

\begin{remark} It might well happen that the characteristic image $\Phi(\bbA^n)$ is Zariski-closed in  $\bigwedge^n \lie{g}^*$  and also contains $0$.
\end{remark}

For \'etale affine representations the above criterion implies:
\begin{example} Let $n = \dim G$. 
Then $G$ acts simply transitively on $\bbA^n$ if and only if $\Phi$ is  constant, and non-zero. In particular, $\Phi(\bbA^n) = \{ \Phi_{0} \}$, 
where $\Phi_{0} \neq 0$. 
\end{example}

\subsection{Transitivity of nilpotent prehomogeneous groups} \label{sect:niltransitive}

As a direct application of Proposition \ref{prop:ccriterion2}, we
can deduce that, for nilpotent groups, transitivity is equivalent with
volume preservation: 

\begin{corollary} \label{cor:niltransitive}
Let  $\rho: G \ra \Aff(n)$ be a prehomogeneous 
affine representation,  where $G$ is a connected nilpotent Lie group.
Then the following are equivalent: 
\begin{enumerate}
\item[i)] 
$\rho$ is transitive.
\item[ii)]  $\rho(G)$ is unipotent. 
\item[iii)] 
$\rho$ is volume 
preserving (that is, $\ell(\rho(G)) \subset \SL(n)$). 
\end{enumerate}
\end{corollary}
\begin{proof} 
Suppose $G$ is acting transitively on $\bbA^n$. Without loss of generality, we may  assume that $G$ is an algebraic subgroup of $\Aff(n)$. A connected nilpotent linear algebraic group 
is a direct product $G = U(G) T$, 
where $U(G)$ is unipotent and $T$ is a connected abelian group of
diagonalisable matrices, which is contained in the center of $G$, cf.\ \cite{Borel}. 
Since, $U(G)$ acts transitively, its centraliser is unipotent.
Hence, $T = \{ 1 \}$ and $G = U(G)$ is unipotent. Thus, i) implies ii). 

Since every algebraic character of a unipotent group is trivial,  ii) implies iii). 

Suppose now that $\rho(G)$ is volume preserving. Let $x_{0} \in U_{\rho}$. Then the characteristic orbit 
$\Phi( G \, x_{0})  = c(G)^* \Phi_{x_{0}}$ is (Zariski-) closed, since the adjoint action of $G$ on $\lg$ is unipotent.  Then
$\Phi(\bbA^n) = \Phi( \ac{ G \, x_{0} }^Z)  \subseteq \ac{ \Phi( G \, x_{0}) }^Z = c(G)^* \Phi_{x_{0}}$ does not
contain $0$.  Hence, $G$ acts transitively. Thus, iii) implies i). 
\end{proof} 



\section{Characteristic cohomology classes associated with affine representations}
\label{sect:Ccc}

Let $\rho: G \ra \Aff(n)$ be an affine representation. 
The characteristic map $\Phi_{\rho}$  gives
rise to certain characteristic cohomology classes, which 
carry information about $\Phi_{\rho}$, and, in particular, about 
the transitivity properties of $\rho$.

\subsection{Construction of the characteristic classes} \label{sect:charclasses}
Let $\bar{\rho}: \lie{g} \ra \aff(n)$ denote the differential of $\rho$, and 
$\bar{\ell}: \lie{g} \ra \lie{gl}(n)$ the corresponding linear part. 
The representation $\bar{\ell}$  turns $\bbR^n$ into
a $\lie{g}$-module, which we denote by $\bbR^n_{\bar{\ell}}$. 
Furthermore, we let $\bbR_{\bar{\ell}}$ denote the one-dimensional $\lie{g}$-module, 
which is induced by the trace of $\bell$. We shall also consider the Lie algebra 
cohomology groups $H^n(\lg, \bbR_{\bell})$,  $H^n(\lg, \lh, \bbR_{\bell})$,
where $n = \dim \lg / \lh$.  (See Appendix \ref{sect:Lieco} for definitions
and notation on Lie algebra cohomology.) 

\begin{lemma} \label{lemma:relclass}
Let $x \in \bbA^n$ and $\lh = \lg_{x}$ the
stabiliser of $x$ under the action $\bar{\rho}$. Let $\Phi: 
\bbA^n \ra \bigwedge^{n} \lg^*$ be the characteristic map for $\rho$. 
Then $\Phi_{x}$ is an element of
$C^n(\lg/\lh,  \bbR_{\bell})^\lh 
\subseteq C^n(\lg, \bbR_{\bell})$. In particular,  $$ d_{\bell} \,  \Phi_{x} = 0 \; . $$ 
Moreover, if $x \in U_{\rho}$ then $\Phi_{x} \neq 0$, and $ \lh$ acts 
trivially on $C^n(\lg/\lh,  \bbR_{\bell})$.  
\end{lemma}
\begin{proof}  Clearly, $\Phi_{x}$ is in $C^n(\lg/\lh,\bbR_{\bell})$. By the
relation \eqref{eq:fundform}, the form $\Phi_{x}$ is stabilised by
the twisted adjoint action of $H= G_{x}$. Taking derivatives shows 
that $\Phi_{x} \in  C^n(\lg/\lh,\bbR_{\bell})^\lh$ is an invariant for $\lh$. 
Since  $C^*(\lg/\lh,\bbR_{\bell})^\lh$ is a
subcomplex for $d_{\bell}$ and $C^k(\lg/\lh,M) =\{ 0 \}$, $k>n$, (for any module
$M$),
we deduce that $d_{\bell} \Phi_{x} = 0$. Since  $C^n(\lg/\lh,\bbR_{\bell})$
is one-dimensional and spanned by $\Phi_{x}$, $\lh$  acts trivially.
\end{proof}

\paragraph{The absolute class}

As a consequence of  relation \eqref{eq:fundform}, we deduce:

\begin{proposition} \label{prop:absclass} 
The forms $\Phi_{x}$, $x \in \bbA^n$, are elements of the cocycle 
vector space  $Z^n(\lg, \bbR_{\bell})$, and represent a unique cohomology class 
$ \bar c = [ \Phi_{x} ] \in H^n(\lg, \bbR_{\bell})$, which does not depend on $x$. 
\end{proposition}
\begin{proof} 
By Lemma \ref{lemma:relclass}, $\Phi_{x}$
is a cocycle.  
We show now that the associated cohomology classes $[\Phi_{x}]
\in H^n(\lg, \bbR_{\bell})$ do not depend on $x$. We argue as follows. 
Recall that the representation $\ell$ defines a twisted adjoint action 
of $G$ on  the cohomology groups $H^k(\lie{g}, \bbR_{\bar{\ell}})$, whose
derivative is induced by the operators $L_{X}$, $X \in \lg$, cf.\ Appendix \ref{sect:Lieco}. 
Assuming $G$ is connected, this inner action of $G$ on the cohomology 
is trivial, cf.\  \cite{EMSII}. If $\rho$ is not prehomogeneous then $\Phi_{x} = 0$, 
for all $x \in \bbA^n$. Hence, $[ \Phi_{x} ] = 0$, and the claim is proved. 
Otherwise, let $x_{0} \in U_{\rho}$. The transformation formula  \eqref{eq:fundform} shows
that the forms $\Phi_{x}$, $x \in G \, x_{0}$, form  a twisted adjoint orbit in 
$Z^n(\lg, \bbR_{\bell})$. By the above remark, they represent the same cohomology 
class $c \in H^n(\lg, \bbR_{\bell})$.
Therefore, the polynomial map $\bbA^n \ra H^n(\lg, \bbR_{\bell})$, $x \mapsto [\Phi_{x}]$ 
is constant on the  Zariski dense subset $G \, x_{0}$. Hence, it is constant.
\end{proof}

\begin{definition} We call $\bar c = \bar c(\rho) \in H^n(\lg, \bbR_{\bar{\ell}})$
the (absolute) characteristic class of the affine representation $\rho$. 
\end{definition}

The geometric importance of the absolute class $\bar c(\rho)$ stems from: 

\begin{proposition} Let $\rho: G \ra \Aff(n)$ be an affine representation. If 
the absolute characteristic class $\bar c(\rho)$ is different from
$0$ then $\rho$ is transitive. 
\end{proposition}
\begin{proof} 
Assume that $\rho$ is not transitive. Then, by Proposition \ref{prop:ccriterion1}, 
$\Phi(\bbA^n)$ 
contains $0$. By Proposition \ref{prop:absclass}, $\bar c(\rho) = 0$.
\end{proof} 

However, the class $\bar c$ may vanish for a transitive 
affine representation. The following 
example of a nilpotent transitive group is taken from \cite{GH_2}: 
\begin{example} \label{ex:GHex} 
Consider 
$$ \lg = \left\{ X(s,t,u,v) =  \begin{matrix}{cccc} 
0 & t & v & s \\ 0 & 0 & u & t \\ 0 & 0 & 0 & u \\
0 & 0 & 0 & 0 
\end{matrix} \mid s,t,u,v  \in \bbR \right\} \; \subset \lie{aff}(3)
$$
Then $\lg$ is the tangent algebra to a transitive
action of a four-dimensional unipotent group $G \leq \Aff(n)$.
The tangent algebra of the stabiliser at $0 \in \bbA^3$ is 
$\lh = \{ X(0,0,0,v) \mid v \in \bbR \}$. We put $\sigma, \tau, \omega, \nu$, for 
the left invariant one-forms
dual to $S= X(1,0,0,0)$, $T= X(0,1,0,0)$, $U= X(0,0,1,0)$, $V = X(0,0,0,1)
\in \lie{g}$. Then  $d\tau= d \omega = 0$, and $d\nu = - \tau \wedge \omega$,
$d\sigma = \omega \wedge \nu$.  Thus $\Phi_{0} = \sigma \wedge \tau \wedge \omega = 
- d(\sigma \wedge \nu)$. But note that $H^3(\lg,\lh, \bbR) \, = \,  <\! [ \Phi_{0} ] \! > \,  \neq \{ 0 \}$. 
\end{example}

The example motivates the following construction of the relative classes.   

\paragraph{The relative classes} Let $x \in \bbA^n$. 
We put $\lh = \lg_{x}$ for the
stabiliser of $x$ under the action $\bar{\rho}$.
%
By Lemma \ref{lemma:relclass}, the form $\Phi_{x}$ represents a 
cohomology class  $\bar{c}(\rho,x) \in H^n(\lg,\lh, \bbR_{\bell})$.

\begin{definition} We call $\bar{c}(\rho,x) = [ \Phi_{x} ] \in H^n(\lg, \lh,  \bbR_{\bell})$
the \emph{relative} characteristic class of the affine 
representation $\rho$, at $x \in \bbA^n$. 
\end{definition}

Note that the forms  $\Phi(\rho,x)$ vanish precisely outside the union of
all open orbits of $G$. 
Hence:

\begin{lemma} If $\bar{c}(\rho,x) \neq 0$, for some $x \in \bbA^n$, then
$\rho$ is prehomogeneous.
\end{lemma}

In fact, we also have:
\begin{lemma} Let $x \in U_{\rho}$. Then $\bar{c}(\rho,x) \neq 0$ if and only 
if $H^n(\lg, \lh,  \bbR_{\bell}) \neq 0$.  
\end{lemma}

As example \ref{ex:GHex} shows, there exist transitive representations
$\rho$, where $G$ is nilpotent, and the absolute class $\bar c(\rho)$ 
is zero, while all relative versions $\bar{c}(\rho,x)$ are non-zero. 
In fact, in  in section \ref{sect:significancerc} below we show that 
$\bar{c}(\rho,x)$ is non-zero, whenever $\rho$ is transitive. 
We will also show that, \emph{for nilpotent groups $G$}, the 
non-vanishing of the class $\bar{c}(\rho,x)$, for some $x \in \bbA^n$, 
implies transitivity of $\rho$.

\paragraph{Naturality properties 
} 
\begin{enumerate}
\item By construction,  the relative classes $\bar{c}(\rho,x)$ map to the absolute class
$\bar{c}(\rho)$ under the natural homomorphism
$$   H^n(\lg,\lh, \bbR_{\bell}) \ra H^n(\lg, \bbR_{\bell}) \; .$$
\item Let $L \leq G$ be a closed subgroup, and let $\rho_{L}: L \ra \Aff(n)$
denote the restriction of $\rho$ to $L$.  Then in the following commutative diagram 
(vertical arrows denoting restriction homomorphisms)
$$ 
\begin{CD} 
H^n(\lg,\lh, \bbR_{\bell}) @>>> H^n(\lg, \bbR_{\bell}) \\
 @VVV     @VVV  \\
  H^n(\lie{l},  \lie{l} \cap \lh, \bbR_{\bell})  @>>> H^n(\lie{l}, \bbR_{\bell})
\end{CD}
$$
the classes $\bar{c}(\rho,x)$, $c(\rho)$,  map to 
$\bar{c}(\rho_{L},x)$,  $c(\rho_{L})$ respectively. 
\item Let $A \in \N_{\Aff(n)}(\rho(G))$, and $c(A): \lg \ra \lg$ 
the induced automorphism of $\lg$. If $\lh = \lg_{x}$ then 
$c(A) \lh = \lg_{Ax}$. Consider the map
of cohomology groups
$$ c(A)^*: H^n( \lg, \lg_{Ax}, \bbR_{\bar \ell}) \ra H^n( \lg, \lg_{x}, \bbR_{\bar \ell}) . $$
Then,  by \eqref{eq:fundform},
 $$ c(A)^* \bar{c}(\rho,A x) = \det \ell(A) \bar{c}(\rho,x) \; . $$
\end{enumerate}


\subsubsection{Alternative construction of the characteristic classes} \label{sect:radobs}
Let $\Gamma$ be a group. To any affine representation $\rho: \Gamma \ra \Aff(n)$ the 
translation part $t: \Gamma \ra \bbR^n$ defines a cohomology class 
$[t]  \in H^1(\Gamma, \bbR^n_{\ell})$, \emph{which vanishes if and only if $\Gamma$
has a fixed point on $\bbA^n$}.  (See \cite {Hirsch,GH_1} for some
applications.) \\

Analogously, for a Lie algebra 
representation $\bar{\rho}: \lie{g} \ra \aff(n)$, the translation part 
$$ t: \lg \ra \bbR^n $$ defines
a class $u \in H^1(\lie{g}, \bbR^n_{\bar{\ell}})$, which is called the 
{\em radiance obstruction}. 
Goldman and Hirsch \cite{GH_1, GH_2} considered the 
exterior powers $$ u^{i} = \Lambda^{i} u \in  H^{i}(\lg,\bigwedge^{i} \bbR^n_{\bar{\ell}}) \; . $$ 
By \cite{GH_2}, if $u^{i} \neq 0$ then every orbit of $G$ has at least dimension $i$.
Note that the fixed $n$-form $\nu$ defines an 
isomorphism of $\lie{g}$-modules $$\bigwedge^n \bbR^n_{\bell} \cong 
\bbR_{\bell} \; . $$ Under this  isomorphism, the $n$-th exterior power of $u$ 
identifies with $ \bar c(\rho)$: 
$$        \bar c(\rho) \,  = \,  \Lambda^{n} u  \; . $$

The representation of $ \bar c(\rho)$ as an exterior product of a class 
in degree one, leads to an important application for compact affine manifolds, 
see section \ref{sect:GHvol}.



\subsection{Significance of the relative classes}  \label{sect:significancerc}
Let $\rho: G \ra \Aff(n)$ be a prehomogeneous representation.
We show below that the cohomology groups $H^n(\lg,\lh,\bbR_{\bell})$
are computed by a characteristic character $\chi_{G/H}(\rho)$, 
which is associated to $\rho$ and $G$. We then explain how  the
group $H^n(\lg,\lh,\bbR_{\bell})$ is linked to the geometry of the semi-invariant
measure, which is induced by $\rho$ on the homogeneous space $X= G/H$. 
As a first application, we deduce that,  
\emph{for a transitive affine representation $\rho$},
all relative cohomology classes $\bar{c}({\rho,x})$ are non-zero
(contrary to what may happen for the absolute
class $\bar{c}({\rho})$). 



\subsubsection{The character $\chi_{G/H}(\rho)$}
Let $x \in U_{\rho}$.
We put $N_G(H)$ for the normaliser of 
$H= G_{x}$ in $G$.  The adjoint representation $\Ad: N(G) = N_{\Aff(n)}(G) \ra \GL(\lg)$ induces a
quotient representation
$$\Ad_{G/H}:  N_G(H) \ra \GL(\lg/\lh) \; . $$
We then define the \emph{unimodular character for $G/H$ and $\rho$} as 
$$  \chi_{G/H}(\rho) \,  =  \,  \det \ell(\rho) \, \det \Ad_{G/H}^{-1} : N_G(H) \ra \bbR^{>0}  
\; .   $$ 
The group $N(G)$ acts 
on $C^n(\lg, \bbR_{\bell})$ by the twisted adjoint action as
in equation \eqref{eq:fundform}. Since $N_G(H)$ normalises $H$ and $\lh$,  
this action preserves the subspace 
$$ C^n(\lg/\lh,\bbR_{\bell}) \subseteq C^n(\lg, \bbR_{\bell}) \; . $$ 
Then
$N_G(H)$ acts on the one-dimensional space $$ C^n(\lg/\lh,\bbR_{\bell}) \cong \Hom(\Lambda^n \lg/\lh, \bbR_{\bell})$$ with the character $\chi_{G/H}(\rho)$.  
For all $A \in N_G(H)$, equation \eqref{eq:fundform}
turns to
\begin{equation} \label{eq:fundformngh} 
\Phi_{A x}  \; = \;     \chi_{ G/H}(\rho) (A) \;  \Phi_{x}    \; \; 
 . 
\end{equation} 
In particular, this implies that $\chi_{G/H}(\rho)$ factorises over $H$, that is,
$$ \chi_{G/H}(\rho)_{| H} \equiv 1   \; . $$

\paragraph{Associated Lie algebra cohomology group}

By Proposition \ref{prop:biinvmeasures}, the character $ \chi_{ G/H}(\rho)$ computes the
cohomology group $H^n(\lg,\lh,\bbR_{\bar{\ell}})$:

\begin{proposition} \label{prop:semiinvmX}
The following conditions are equivalent:
\begin{enumerate}
\item $ \chi_{ G/H}(\rho): (N_G(H)/H)^0 \ra \bbR^{>0}   \; \, \equiv 1$.
\item  $H^n(\lg,\lh,\bbR_{\bar{\ell}}) \neq  \{ 0 \}$. 
\end{enumerate}
\end{proposition} 

\begin{example} If $N_G(H)^0 = H^0$ then $H^n(\lg,\lh,\bbR_{\bar{\ell}}) \neq  \{ 0 \}$. 
For example, the automorphism 
group $G= \Aut({\cal H}_{2})$ of a halfspace has
$H^2(\lg,\lh,\bbR_{\bar{\ell}} )= \bbR$.  
\end{example}

\paragraph{Induced semi-invariant measure on X= G/H}
We put $X= G/H$ for the associated homogeneous space, where $H = G_{x}$. 
The pull back of the parallel volume form $\nu$ on $\bbA^n$ defines
a semi-invariant measure on $X$. Its properties strongly interact with the 
geometry of the representation $\rho$. \\ 

The following geometric 
interpretation of the characteristic character $\chi_{G/H}(\rho)$
is a consequence of  Lemma \ref{lemma:biinvmeasures}:  

\begin{proposition} \label{prop:semiinvmX2}
Let $\rho: G \ra \Aff(n)$ be a prehomogeneous representation, $x \in U_{\rho}$ 
and $X= G/H$ the associated homogeneous space.
Then the following conditions are equivalent:
\begin{enumerate}
\item $ |\chi_{ G/H}(\rho)|: N_G(H)  \ra \bbR^{>0}   \; \, \equiv 1$.
\item The semi-invariant measure on $X$, which is induced from a parallel volume form $\nu$ on $\bbA^n$, is invariant by the right-action of $N_G(H)$ on $X$.
\end{enumerate}
\end{proposition} 


\paragraph{Homogeneous domains with non-vanishing relative class}
We say that the affinely homogeneous domain $U$ has a
\emph{non-vanishing relative class}, if the relative 
cohomology groups for the group $G= \Aff(U)$ are non-vanishing.

\begin{corollary} \label{cor:cnvcentisu2}
Let $U$ be a homogeneous affine domain, $G= \Aff(U)$, and
$H = G_{x}$, for  $x \in U_{\rho}$. If $H^n(\lg,\lh,\bbR_{\bar{\ell}}) \neq  \{ 0 \}$ 
then the center $\Z(G)$ of $\Aff(U)$ is a connected unipotent group. 
\end{corollary}
\begin{proof} Clearly, $Z(\Aff(n),G)^0 \leq N_G(H)$. 
Thus, by our assumption, the character $ \chi_{ G/H}(\rho)$ is trivial on
$Z(\Aff(n), G)^0$.  Since 
$Z(\Aff(n),G)^0$ is in the kernel of the representation $\Ad_{G/H}$, this implies 
that the elements of  $Z(\Aff(n),G)^0$
are volume preserving on $\bbA^n$.  By Corollary \ref{cor:centisu}, 
$Z(\Aff(n),G)$ is a connected  unipotent group. Thus, 
$\Z(G) = Z(\Aff(n), G)$.
\end{proof}

\begin{example} Among the two-dimensional homogeneous domains, 
precisely  $U =\bbR^2 - \{0\}$, and the sector $U= {\cal Q}_{2}$ 
have vanishing relative class.
\end{example}

\subsubsection{Transitive representations}
As remarked before, $\chi_{G/H}(\rho)$ is trivial on $H$, for any
prehomogeneous affine representation $\rho$. If $\rho$ is transitive 
then $\chi_{G/H}(\rho)$ is trivial on all of $N_G(H)$: 

\begin{corollary} \label{cor:chingh}
If $\rho$ is a transitive affine representation then
$| \chi_{G/H}(\rho) | \equiv 1$ on $N_G(H)$. 
\end{corollary}
\begin{proof} By Proposition \ref{prop:ccriterion2}, the orbit 
$\Phi_{A x} $, $A \in N_G(H)$ must be bounded away from $0$. 
By formula \eqref{eq:fundformngh}, this can only be if 
$\chi_{G/H} (A) \in \{ 1, -1\}$. 
\end{proof}
 
We have the following consequence:

\begin{theorem} \label{thm:transitive}
Let $\rho: G \ra \Aff(n)$ be a transitive affine 
representation. For $x \in \bbA^n$, put $\lh = \lg_{x}$. 
Then, for all $x \in \bbA^n$, the cohomology class 
$\bar{c}(\rho,x) \in H^n(\lg,\lh,\bbR_{\bell})$ is non-zero. 
\end{theorem} 
\begin{proof}  
Since $G$ is transitive, the character 
$ \chi_{G/H}$ is trivial, by Corollary \ref{cor:chingh}.
This implies, 
$H^n(\lg,\lh,\bbR_{\bell}) \cong \bbR$, by Proposition \ref{prop:semiinvmX}.
More specifically, the natural map
induces an  isomorphism of the group of relative 
$n$-cocycles $C^n(\lg,\lh,\bbR_{\bell})$ with $H^n(\lg,\lh,\bbR_{\bell}) \cong \bbR$.
Since $\Phi_{x} \neq 0$, $[ \Phi_{x} ] = \bar c(\rho,x) \in H^n(\lg,\lh,\bbR_{\bell})$
is a non-zero generator of $H^n(\lg,\lh,\bbR_{\bell})$. 
\end{proof}

Remark that, conversely, if,  for all $x \in \bbA^n$, $\bar c(\rho,x) \neq 0$ then
$G$ acts transitively. 

\begin{corollary} Let $\rho: G \ra \Aff(n)$ be a transitive affine 
representation. Then  $H^n(\lg,\lh,\bbR_{\bell}) \neq \{ 0 \}$.
\end{corollary}

We thus note (compare Proposition \ref{prop:semiinvmX2}): 
 
\begin{corollary} Let $\rho: G \ra \Aff(n)$ be a transitive affine 
representation.  Then the semi-invariant measure on 
$X=G/H$ with character $\det \ell(\rho)$ is right-invariant by $N_G(H)$.
\end{corollary}

For a transitive \emph{algebraic} representation, we add the 
following observation:

\begin{proposition} \label{prop:transvol}
Let $\rho: G \ra \Aff(n)$ be a transitive algebraic affine 
representation. Then the following are equivalent: 
\begin{enumerate}
\item $\rho$ is volume preserving.
\item $\det \Ad_{\lg/\lh} \equiv 1$ on $H$.
\item $\det \Ad_{\lg/\lh} \equiv 1$ on $N_G(H)^0$.
\end{enumerate}  
\end{proposition}
\begin{proof} Put $\lambda = \det \ell(\rho)$. If $\rho$ is
volume preserving then $H \leq \ker \lambda$. Moreover,
$\chi_{G/H}(\rho) \equiv 1$ on $H$, as for any prehomogeneous
representation. Thus, $\det \Ad_{\lg/\lh} \equiv 1$ on $H$. 
Conversely, assume $\det \Ad_{\lg/\lh} \equiv 1$ on $H$. 
Then $H \leq \ker \lambda$. 
Since the unipotent radical of $\rho(G)$ acts transitively, $\ker \lambda$ acts
transitively. Thus $\rho(G) 
= \rho( ( \ker \lambda) \,  H) = \rho( \ker \lambda) $. Therefore,
 $\lambda \equiv 1$.
It follows the equivalence of 1.\ and 2.  

By Corollary \ref{cor:chingh}, the analogous argument may be used
to show the equivalence of 1.\ and 3. 
\end{proof}

\begin{corollary} \label{cor:transvol}
Let $\rho: G \ra \Aff(n)$ be a transitive affine 
representation.  Then $\rho$ is volume preserving if and
only if $X=G/H$ admits a $G$-invariant measure.  
\end{corollary}

\subsection{Affine representations of nilpotent groups}
In the case of nilpotent groups, we can 
summarise  as follows: 

\begin{theorem} \label{thm:niltransitive2}
Let $\rho: G \ra \Aff(n)$ be an affine prehomogeneous representation, 
where $G$ is  a connected nilpotent Lie group. Let 
$x \in U_{\rho}$, and put $H = G_{x}$.  Then the following are equivalent: 
\begin{enumerate}
\item $\rho$ is transitive on $\bbA^n$.
\item $\chi_{G/H}(\rho) \equiv 1$ on $N_G(H)^0$.
\item  $ H^n(\lg,\lh,\bbR_{\bell}) \neq 0$.
\item The semi-invariant measure induced by a parallel volume on 
$X= G/H$ is 
invariant by the right-action of $N_G(H)^0$.
\item The induced semi-invariant measure on $X$ is a $G$-invariant measure. 
\item $\rho(G)$ is volume preserving. 
\item $\rho(G)$ is unipotent. 
\end{enumerate}
\end{theorem}
\begin{proof} If $\rho$ is transitive then  $ H^n(\lg,\lh,\bbR_{\bell}) \neq 0$, 
by Theorem \ref{thm:transitive}. By Proposition \ref{prop:semiinvmX} and
Proposition \ref{prop:semiinvmX2}, this implies that 
$\chi_{G/H}(\rho) \equiv 1$ on $N_G(H)^0$, and also that the induced
semi-invariant measure is invariant by $N_G(H)^0$.  Since $G$ is nilpotent this implies 
(cf.\ Example \ref{example:nilinvm}) that the character $\det \ell ( \rho) $
is trivial on $N_G(H)^0$. Moreover $\rho(G) \leq \ac{\rho(G)}^\bbR =TU$, 
where $T$ is a central subgroup of diagonalisable elements, 
and $U$ is unipotent. By functoriality we may assume now that 
$\rho(G) = \ac{\rho(G)}^\bbR$. It follows that $\det \ell (\rho) =1$  on
the central subgroup $T$. Since $\det = 1$ on $U$, $\rho(G)$ is
volume preserving, and, in particular, the  
induced measure is an invariant measure. 
By Corollary \ref{cor:niltransitive},  $\rho(G)$ is volume preserving
implies that $\rho$ is transitive.    
\end{proof}

Moreover, transitivity of a nilpotent affine representation is 
determined by the relative class
$\bar{c}(\rho_{L},x)$  at an arbitrary point $x \in \bbA^n$:  

\begin{corollary} \label{cor:niltransitive2}
Let $\rho: G \ra \Aff(n)$ be an affine representation, 
where $G$ is  a connected nilpotent Lie group, and let 
$x \in \bbA^n$. 
Then $\rho$ is a transitive representation 
 if and only if $$ 0 \,  \neq \bar{c}(\rho,x) \, \in H^n(\lg,\lh, \bbR_{\bell}) \; . $$
\end{corollary}
\begin{proof}  Assume that $\bar{c}(\rho,x) \in H^n(\lg,\lh, \bbR_{\bell})$ 
is different from zero. In particular,  this implies that 
$0 \neq \Phi_{x} \in C^n(\lg,\lh,\bbR_{\bell})$.
Thus, the representation $\rho$ is prehomogeneous at $x$. 
By Theorem \ref{thm:niltransitive2}, $\rho$ is transitive. 
\end{proof}

Note, that in the situation of Theorem \ref{thm:niltransitive2}, $ \bbR_{\bell} = \bbR$ is
the trivial $\lg$-module. \\

Theorem \ref{thm:niltransitive2}, together with Corollary  
\ref{cor:niltransitive2} characterise transitivity
for nilpotent affine actions. In particular, the results summarise various
special cases, which have been obtained in the literature before: see
\cite{Kim} for \' etale affine representations of nilpotent
groups, \cite{GH_2} for relation with the absolute classe $\bar{c}(\rho)$, 
and \cite{DuIh_2} for nilpotent representations with invariant scalar
product (compare also sections \ref{sect:psrdomains}, \ref{sect:sympdomains}).  

\subsubsection{Minimal classes for algebraic group actions}
Let $\rho$ be an algebraic representation.
\begin{definition}
A relative class $\bar c(x,\rho)$ will be called \emph{minimal} if 
$N_G(H)$ contains a maximal $\bbR$-split torus in $G$,
where $H = G_{x}$. 
\end{definition}
For example, if $G$ is nilpotent
then $\bar c(x,\rho)$ is minimal, for all $x \in \bbA^n$.
 Every reductive subgroup of 
$G$ fixes a point $x \in \bbA^n$. Therefore, every 
$\rho$ has minimal relative classes $\bar c(x,\rho)$. \\

The following result
generalises Corollary \ref{cor:niltransitive2}:

\begin{theorem} \label{thm:transitive_a}
Let $\rho: G \ra \Aff(n)$ be an affine representation, where $G$ is algebraic.
Let $x \in \bbA^n$ such  
that $\bar{c}(\rho,x)$ is minimal, and put $\lh = \lg_{x}$.  
Then $\rho$ is transitive if and only if 
$\bar{c}(\rho,x) \in H^n(\lg,\lh, \bbR_{\bell})$ is
different from zero. 
\end{theorem}
\begin{proof}  Assume $\bar{c}(\rho,x) \neq 0$. Then, in particular,
$\rho$ is prehomogeneous at $x$. Since $H^n(\lg,\lh, \bbR_{\bell})$
is different from zero, $\chi_{G/H}(\rho) \equiv1$, on $N_G(H)^0$.
In particular, by \eqref{eq:fundformngh}, $N_G(H)^0$ is contained
in the stability group of $\Phi_{x}$. 
By Proposition \ref{prop:minpoints}, the orbit $c(G)^* \Phi_{x}$ 
is closed in $C^n(\lg,\lh,\bbR_{\bell})$. Hence, Proposition \ref{prop:ccriterion2}
implies that $\rho$ is transitive.
\end{proof}

\subsection{Prehomogeneous representations with reductive stabiliser} 
\label{sect:red-stab}

A homogeneous space $G/H$ is called \emph{reductive} 
if the adjoint action of $H$ on $\lg$ is reductive. In particular, the 
Lie algebra $\lh$ acts reductively on $\lg$, and there is
a $\lh$-invariant direct decomposition $\lg = \lh \oplus \lie{p}$,
where $\lie{p}$ is a vector subspace isomorphic to $\lg/\lh$, 
and $\lh$ acts reductively on $\lg/\lh$.  Let $\rho: G \ra \Aff(n)$ be a 
prehomogeneous affine representation. We call $\rho$  
\emph{reductive} if, for some $x \in U_{\rho}$, the isotropy algebra
$\lh$ for the homogeneous space $G/H$, $H = G_{x}$ is reductive 
in $\lg$.\footnote{A particular important class of examples appears in the work of
Sato and Kimura  \cite{SK} on regular reductive
prehomogeneous vector spaces, where both $G$ and $H$
are assumed reductive. See also \cite{Kimura}, and Example \ref{RPVs}.}
In particular $H^0$ is a reductive subgroup in $\Aff(n)$.
Also $\rho$ is reductive if $ \lh = \{ 0 \}$. This covers
the important special case of \'etale affine representations, as well.

\paragraph{The character $\chi$} 
Let $\rho: G \ra \Aff(n)$ be a prehomogeneous representation, where
$G$ is connected. We define the \emph{characteristic character} 
$$ \chi= \chi(\rho): \; G\,  \ra \, \bbR^{>0}$$ by putting 
$$   \chi(\rho)(g) \, = \, \det \ell(g)  \det \Ad(g)^{-1} \; . $$
Let  $\Ad_{\lh}$
denote the adjoint representation of $N_G(H)$ restricted to $\lie{\lh}$.
On  $N_G(H)$ we have
$  \chi(\rho) = (\det \Ad_{\lh})^{-1} \, \chi_{G/H} (\rho)$.

\begin{lemma} \label{lemma:chired}
Let $\rho$ be a reductive representation, where
$G$ is connected. Then $\chi(\rho) = 
 \chi_{G/H} (\rho)$ on $N_G(H)^0$. Moreover,  $\chi(\rho) \equiv 1$
 if and only if $\chi_{G/H} (\rho) \equiv 1$ on $N_G(H)^0$.
\end{lemma}
\begin{proof} 
Since $\lh$ is a reductive Lie Algebra $\det \Ad_{\lh}(h) = 1$, for 
all $h \in H^0$. Let $\lie{p}^\lh$ be the elements of $\lie{p}$, centralised
by $\lh$. Since $\lh$ acts reductively, 
\begin{equation} \label{eq:ndec} 
\lie{n}(\lg, \lh) = \lh +  \lie{p}^\lh 
\end{equation}
and 
\begin{equation} \label{eq:ldec}
\lg = \lie{n}(\lg, \lh) + [ \lie{p}, \lh ] \; .    
\end{equation}
It follows from \eqref{eq:ndec} that $\det\Ad_{\lh}(n) = 1$, for all $n \in N_G(H)^0$. 
Therefore,  $\chi(\rho) = \chi_{G/H} (\rho)$ on $N_G(H)^0$. 
In particular,  $\chi(\rho) \equiv 1$, implies $\chi_{G/H} (\rho)
\equiv 1$ on $N_G(H)^0$. 

For the converse, assume that 
$\chi_{G/H} (\rho) \equiv 1$ on $N_G(H)^0$. Then the 
corresponding infinitesimal character $\bar{\chi}_{G/H}: \lie{n}(\lg,\lh) \ra \bbR$
vanishes. Since $\bar{\chi}_{G/H}$ is the restriction of the infinitesimal
character $\bar{\chi}: \lg \ra \bbR$, which belongs to $\chi({\rho})$, the
decomposition \eqref{eq:ldec} implies that $\bar{\chi} \equiv 0$. 
Hence, the character $\chi({\rho})$ is constant. 
\end{proof}

Note, in particular, that $\chi(\rho) \equiv 1$ on $H$. 

 \begin{corollary} 
 Let $\rho: G \ra \Aff(n)$ be a reductive prehomogeneous 
 representation, which is transitive.  Then $ |\chi(\rho)| \equiv 1$.
\end{corollary}

Recall (cf.\  \cite{Koszul, HS}) that $\lh$ is called 
\emph{not homologous to zero in $\lg$} if the
restriction homomorphism $H^q(\lg,\bbR) \ra H^q(\lh, \bbR)$
is surjective, for all $q$. For example, if $\lg = \lie{j} \oplus \lh$,
where $\lie{j}$ is an ideal in $\lg$ then $\lh$ is not homologous to zero.
If $\lh$ is reductive in $\lg$, then surjectivity of $H^q(\lg,\bbR) \ra H^q(\lh, \bbR)$,
$q = \dim \lh$, is sufficient for $\lh$ being not homologous to zero.
If $\lh$ is not homologous to zero then (cf.\  \cite[Theorem 12]{HS})
it follows, in particular, that the natural map on cohomology $H^\bullet(\lg,\lh,\bbR_{\lambda}) \ra 
H^\bullet(\lg, \bbR_{\lambda})$ is injective, for any one-dimensional module $\bbR_{\lambda}$.  \\

We summarise:

 \begin{corollary} \label{cor:redcharchar}
 Let $\rho: G \ra \Aff(n)$ be a reductive prehomogeneous 
 representation.  If\/  $\lh = \lie{g}_{x}$ is not homologous to zero in $\lg$, then the following
are equivalent: 
\begin{enumerate}
\item $\rho$ is transitive. 
\item $| \chi(\rho) | \equiv 1$. 
\item $H^n(\lg,\lh,\bbR_{\bell}) \neq \{0 \}$.
\item $H^n(\lg,\bbR_{\bell}) \ni \bar{c}(\rho) \neq 0$. 
\end{enumerate}
\end{corollary}

\begin{remark} If $\lh$ is reductive in $\lg$, we could (in lieu of Lemma \ref{lemma:chired})
have applied  Poincar\'e duality for relative
cohomology (see \cite{Knapp}). Thereby, 
$$  H^k(\lg,\lh, \bbR_{\lambda}) \cong H^{n-k}(\lg, \lh, \bbR_{ad} \tensor {} {\bbR_{\lambda}^*})^*  \; , $$
where $n= \dim \lg/ \lh$. The module $ \bbR_{ad} \tensor {}  \bbR_{\bell}^*$ 
is determined by the character $-\bar \chi$.  
Moreover, $$ H^0(\lg, \lh,  \bbR_{ad} \tensor {}  \bbR_{\bell}^* ) = 
H^0(\lg, \lh,  \bbR_{\bar \chi}) = \bbR_{\bar{\chi}}^\lg$$ is non-zero 
if and only if  $\bar \chi  = 0$. 
\end{remark}

\subsubsection{Prehomogeneous domains of reductive algebraic representations}
\begin{corollary} Let $\rho: G \ra \Aff(n)$ be a reductive prehomogeneous 
representation, which is algebraic. Then $\rho$ is transitive if and only
if the absolute class $\bar{c}(\rho) \in H^n(\lie{g}, \bbR_{\ell})$ is non-vanishing.
\end{corollary}
\begin{proof} If $\rho$ is a transitive algebraic representation then 
$\rho(G) = U T$, where $T$ is maximal reductive. In particular, $T = G_{x}$ fixes a
point $x \in \bbA^n$. By Theorem \ref{thm:transitive}, $0 \neq \bar c(\rho,x) \in H^n(\lg, \lie{t}, \bbR_{\bell})$. Since $\lie{t}$ is not homologous to zero in $\lg$, its image 
$\bar{c}(\rho) \in H^n(\lie{g}, \bbR_{\ell})$ is non-zero. 
\end{proof} 

\begin{example}[Regular prehomogeneous vector spaces] \label{RPVs}
Let $\rho$ be a linear prehomogeneous algebraic representation. If $H$
is reductive  then the set $U_{\rho}$ is the
complement of a hypersurface, and it is also the set of real points of 
an affine algebraic variety.  (See, for example, \cite[p.41ff]{Kimura}, for reference.) 
In particular,  any open orbit $\rho (G) x$ 
is (a connected component) of the set of real points of an 
affine algebraic variety.
If furthermore $\rho(G)$ is reductive then the representation space
$(\rho(G), \bbR^n)$ is called a
\emph{regular} prehomogeneous vector space. By \cite[Proposition 2.26]{Kimura}
there exists a \emph{non-degenerate} relative invariant for $\rho(G)$ 
which has character $(\det \rho)^2$. In particular, $\chi_{\rho} = \det \rho \neq 1$.
\emph{It follows, in particular, 
that $H^n(\lg, \lie{h},\bbR_{\ell}) = \{ 0 \}$, for a regular prehomogeneous vector
space.}
\end{example}

Using arguments as in  \cite[loc.\ cit.]{Kimura}, we note:
\begin{proposition} Let $\rho$ be a prehomogeneous algebraic 
representation. If $H$ is reductive then the affine homogeneous domain  
$U_{\rho}$ is the complement of a hypersurface, 
and it is also the set of real points of 
an affine algebraic variety.
\end{proposition}
In particular, in the situation of the proposition, there exists a relative 
polynomial invariant 
$$ \delta: \bbA^n \ra \bbR  \, ,  $$
satisfying $ \delta \ \rho(g)  \, = \, \tau(g) \, \delta$,  
for some character $\tau: G \ra \bbR^{\neq 0}$, and
$$ U_{\rho} = U _{\delta} = \{ x \in \bbA^n \mid \delta(x) \neq 0 \} \; . $$
One can show that 
there is a one-to-one correspondence between algebraic 
characters of  $G$ which factor over $H$, and the irreducible 
components of the complement of $U_{\rho}$. Hence, the group
of characters, which factor over $H$ is of rank one if $U_{\rho} \neq \bbA^n$. 
We conclude that $\chi(\rho)$ must be a power of $\tau$.

 \begin{corollary} \label{cor:redalgcharchar}
 Let $\rho: G \ra \Aff(n)$ be an algebraic  reductive prehomogeneous 
 representation.  Then the following
are equivalent: 
\begin{enumerate}
\item $\rho$ is transitive. 
\item $| \chi(\rho) | \equiv 1$. 
\item $H^n(\lg,\lh,\bbR_{\bell}) \neq \{0 \}$.
\item $H^n(\lg,\bbR_{\bell}) \ni \bar{c}(\rho) \neq 0$. 
\end{enumerate}
\end{corollary}
\begin{proof} 
By the previous corollary,  we have the equivalence of
1.,  3.\ and 4. Assume that 1.\ holds. By the above remarks,
this implies that the character $|\chi(\rho)|$ is constant. Conversely,
if $|\chi(\rho)|$ is constant, then $ \chi_{G/H} (\rho)$ is constant on
$N_{G}(H)^0$. By Proposition \ref{prop:semiinvmX}, this implies 3.
\end{proof}

By Corollary \ref{cor:niltransitive}, the equivalences 1.-3.\ of Corollary \ref{cor:redcharchar} 
and Corollary \ref{cor:redalgcharchar} also hold if $G$ is nilpotent. However, in general, $\chi(\rho) \equiv 1$
does not imply transitivity:

\begin{example} Consider the 2-dimensional linear prehomogeneous representation 
of the group $\SL(2,\bbR)$ on $\bbR^2$. It has open orbit  $\bbR^2 - \{ 0 \}$. 
The stabiliser 
$H= G_{x}$, $x \neq 0$, is a unipotent subgroup. 
Moreover, $\chi_{\SL(2,\bbR)} \equiv 1$,
but $\chi_{\SL(2,\bbR)/H} \neq 1$. In particular, 
$H^2(\lie{sl}(2,\bbR),\lh, \bbR)  = \{ 0 \}$. 
\end{example}

Also, the character $\chi$ may not be trivial for a transitive representation,
as the following example shows: 

\begin{example} Let $G = \bbR^2 \rtimes H \leq \Aff(2)$, where $H \leq \SL(2,\bbR)$
denotes the 2-dim.\ solvable group of upper triangular matrices, and $\bbR^2$ is the
group of translations. Then $G$ is transitive and $\chi_{G} \neq 1$.
\end{example}

\subsubsection{The case of \'etale affine representations}
If the affine representation $\rho$ is \'etale then $\lh = \{ 0 \}$, and there are no 
relative versions of the characteristic class $c(\rho)$. As a special case of 
Corollary \ref{cor:redcharchar}, we obtain: 

\begin{proposition} \label{prop:etreps}
 Let $G$ be a Lie group, and $\rho: G \ra 
\Aff(n)$ be an affine \'etale representation. Then the following
are equivalent
\begin{enumerate}
\item $\rho$ is transitive. 
\item $\chi(\rho) \equiv 1$.
\item $H^n(\lg,\bbR_{\bell}) \neq \{0 \}$.
\end{enumerate} 
\end{proposition} 

This result (in particular, the equivalence of 1.\ and 2.) appears 
in the work of Helmstetter \cite{Helms}, Kim \cite{Kim}, whereas the relation
with the cohomology group $H^n(\lg, \bbR_{\bell})$ is introduced
in \cite{GH_2}. It also plays a role in the theory of 
prehomogeneous vector spaces, see \cite{Kimura}. 
As shown in section \ref{sect:etreps}, if $\rho$ is an \'etale affine representation,
the character $\chi$ relates to a relative invariant $\delta$ of $\rho$. 
This may also be interpreted as follows: 

\begin{example}  \label{ex:rightvol}
Any right-invariant volume
form on $G$ induces a volume form $\eta$ on $U_{\rho}$.
By comparison with the parallel volume $\nu$ on $\bbA^n$, 
we can write $$ \eta = f \nu \; , $$ 
for some function $f$ on $U_{\rho}$.
It follows that $\delta = f^{-1}$ is a relative invariant for $G$ with
character $\chi$, satisfying $\delta \, \eta = \nu$.
 \end{example}
Choose  linear independent affine vector fields $\bar X_{1}, 
\ldots \bar X_{n}$,  whose flows generate $\rho(G)$, and such that
$\eta(\bar X_{1},  \ldots, \bar X_{n}) =1$. Then 
$$ \delta(x) = 
\nu(\bar X_{1}(x),  \ldots, \bar X_{n}(x)) \; . $$
It follows that the function $\delta$ is a polynomial on $\bbR^n$, and 
the zeros of $\delta$ describe the boundary of $U_{\rho}$.  

\section{Compact affine manifolds and prehomogeneous algebraic groups} 
\label{sect:parvol}
 
Here, we develop a theme, which is contained in  a series of 
papers \cite{Hirsch, GH_1,GH_2} by Goldman and Hirsch. Their ideas show a
strong link of the parallel volume conjecture 
with the theory of prehomogeneous affine representations, as
presented in the previous two sections.

\subsection{Holonomy of compact complete affine manifolds}

A basic remark concerning compact complete affine 
manifolds is: 
 
\begin{theorem} \label{thm:comptrans}
Let $M$ be a compact complete affine manifold, and $\Gamma = h(\pi_{1}(M))$ its
affine holonomy group. Then the Zariski closure  
$A(\Gamma)$ of\/ $\Gamma$ acts transitively on $\bbA^n$. 
\end{theorem}
\begin{proof} Let $U$ be the unipotent radical of $A(\Gamma)$. 
Every orbit $Ux$ is a contractible closed submanifold of $\bbA^n$
which is preserved by $\Gamma$. Since $\Gamma$ acts
properly discontinuously and freely on $Ux$, the quotient $\Gamma \backslash Ux$
is a compact manifold. Considering the cohomological dimension
$cd(\Gamma)$
of the group $\Gamma$ over $\bbZ$ (see \cite{Brown}), 
we have $n = \dim \bbA^n = \dim M = cd(\Gamma) = 
\dim \Gamma \backslash Ux$. Thus $U x = \bbA^n$. 
\end{proof}

\paragraph{Polynomial volume forms}
A volume form 
on $M$ is called a \emph{polynomial volume form}
if it is expressed by a polynomial function in the
affine coordinates of $M$. The following observation 
is thus an immediate consequence of Theorem \ref{thm:comptrans} 
and Corollary \ref{cor:transvol}: 

\begin{proposition}
Let $M$ be a compact complete affine manifold. Put $G = A(\Gamma)$, and $H= G_{x}$, for
some $x \in \bbA^n$. Then the following 
conditions are equivalent:
\begin{enumerate}
\item $M$ has a parallel volume form.
\item $M$ has a polynomial volume form.
\item The homogeneous space $G / H$ admits a $G$-invariant volume.
\end{enumerate}
\end{proposition}

In particular, \emph{every polynomial volume form on a compact
complete affine manifold is parallel}. \\

It seems natural to conjecture that on a compact 
affine manifold $M$ a polynomial volume form must be parallel.
In \cite{FGH}, this is proved under the assumption
that the holonomy group $\Gamma$ of $M$ is nilpotent. Not
much seems to be known in the general case.  

\begin{question} Let $M$ be a compact affine manifold.
Is every polynomial volume form on $M$ parallel?
\end{question}

\begin{example}  \label{ex:polvol}
Let $\rho: G \ra \Aff(n)$ be an \'etale affine
representation. Assume that $G$ admits a uniform lattice $\Gamma$. 
Every polynomial volume form on $$ M \, = \, \Gamma \, \backslash G $$ lifts
to a $\Gamma$-left-invariant volume form $\mu$ on $G$, which is polynomial.
We may integrate $\mu$ with respect to a finite $G$-invariant volume $\lambda$ 
on $ \Gamma \, \backslash G$, to obtain a polynomial volume form $\eta$, which is
left-invariant on $G$. Explicitly, if $X_{i}$ are vector fields on $G$,
$$ \eta(X_{1}, \ldots ,X_{n}) =   \int_{\Gamma \, \backslash G} l_{g}^* \mu \, (X_{1}, \ldots, 
X_{n}) \, d\lambda \; . $$
Since $\Gamma$ is uniform, $G$ is unimodular and, therefore, $\eta$ is 
also right-invariant on $G$. 
As in Example \ref{ex:rightvol}, 
we can write $$\eta = f \nu \; , $$ where $\delta = f^{-1}$ is a polynomial. Since
$\eta$ is polynomial, also $f$ is a polynomial. It follows that
$f$ must be constant. Hence, the parallel volume on $G$ is left-invariant.
This implies that $G$ and $M$ are complete (see Corollary \ref{cor:Gcomplete_um}).
By the above, since $M$ is complete, the polynomial volume form $\mu$ is
parallel. 
\end{example}




\subsection{Holonomy of compact volume preserving affine manifolds}
\label{sect:GHvol}

Let $M$ be a compact affine manifold, and $\Gamma = hol(\pi_{1}(M))$ 
its affine holonomy group. We  put $A(\Gamma)$ for the Zariski closure
of $\Gamma$ in $\Aff(n)$. 
The following striking  observation is the main result of  \cite{GH_2}. 

\begin{theorem}[Goldman-Hirsch, \cite{GH_2}] \label{thm:GHvol}
Let $M$ be a compact affine manifold with parallel volume.  Then $A(\Gamma)$ acts
transitively on $\bbA^n$. 
\end{theorem}

We put $G = A(\Gamma)$, and $H= G_{x}$, for
some $x \in \bbA^n$. 
The proof of Theorem \ref{thm:GHvol},  which 
we will explain below, then shows:  

\begin{corollary}  \label{cor:GHvol}
Let $M$ be a compact affine manifold with parallel volume.  
Then the absolute characteristic class  $\bar{c}(G)  \in H^n(\lg, \bbR)$
is non-zero. 
\end{corollary}

We also obtain the following strong restriction on the
homogeneous spaces $G/H$, which may appear as Zariski closures 
of the holonomy groups of compact volume preserving affine manifolds.

\begin{corollary} 
Let $M$ be a compact affine manifold with
parallel volume.  Then $$ \dim H^n(\lg,\lh, \bbR) = 1 \; , $$ and the
natural homomorphism $$H^n(\lg,\lh, \bbR) \ra 
H^n(\lg,\bbR)$$ is injective. 
\end{corollary}
\begin{proof}  By Theorem \ref{thm:GHvol}, $A(\Gamma)$ acts transitively, 
and by Corollary  \ref{cor:GHvol}, the absolute class $\bar{c}(G) \in H^n(\lg, \bbR)$
is non-zero.  Since, as described in section \ref{sect:charclasses}, $\bar{c}(G)$ is the image
of a relative class $\bar{c}(G,x) \in H^n(\lg,\lh, \bbR)$, the claim follows.
\end{proof}

Here is an important consequence of Theorem \ref{thm:GHvol}
concerning the \emph{divisibility of homogeneous domains}. 

\begin{corollary} Let $U$ be a homogenous domain which is divisible
by a properly discontinuous group $\Gamma$ of affine transformations.
Then either $U = \bbA^n$ or there exists an element $\gamma \in \Gamma$  
with $\det \gamma >1$. 
\end{corollary}
\begin{proof} Follows immediately, since $A(\Gamma)^0 \leq \Aff(U)$ 
preserves $U$. (See section \ref{sect:affautU}.)
\end{proof}
More generally, Theorem \ref{thm:GHvol} implies that the
development image of a compact volume preserving affine manifold
$M$ is not contained in a proper semi-algebraic
subset of $\bbA^n$. (See  \cite{GH_2}).  \\

One can also use Theorem \ref{thm:GHvol} to deduce: 

\begin{corollary}[\cite{GH_2}]
A compact homogeneous affine manifold $M$ with 
parallel volume  is
complete.
\end{corollary}
\begin{proof} In fact, the universal cover $X$ of a homogeneous affine
manifold $M$ develops onto an affine homogeneous domain $U$, and
the development map is a covering map. If $M$ is volume 
preserving $U$ is divisible by a a volume preserving properly 
discontinuous group $\Gamma$ of affine 
transformations. Hence, $U = \bbA^n$. It follows that the development
map of $M$ is a covering map onto $\bbA^n$. Therefore, it is 
a diffeomorphism. Thus, $M$ is complete.
\end{proof}

See section \ref{sect:comphom}, for an independent proof of the latter fact.

For an application of Theorem \ref{thm:GHvol} concerning the
structure of the groups of isometries of compact flat Pseudo-Riemannian
manifolds, see section \ref{sect:compvp}. 

\settocdepth{none}
\subsubsection{Proof of Theorem \ref{thm:GHvol}} 
\settocdepth{subsubsection}

We outline the proof of Theorem \ref{thm:GHvol} using
the main ideas of \cite{GH_2}, which relate nicely to the 
concepts discussed in section \ref{sect:charclasses}. 
The basic new tools which we require stem from the cohomology theory of discrete groups (as documented for example in \cite{Brown}), 
and the cohomology of algebraic linear groups $G$, as developed 
in \cite{Hochschild}. 

\paragraph{Cohomology of algebraic linear groups} 
Let $G= U T$ be an algebraic linear group, 
where $U$ is the unipotent radical and  $T$ is a maximal reductive
subgroup.  We let $\lg$, $\lu$, $\lie t$ denote their corresponding Lie algebras, 
which are subalgebras of $\aff(n)$. Let $V$ be a rational $G$-module. The
algebraic cohomology groups $H^{i}_{alg}(G,V)$ are defined in \cite{Hochschild}. 
For any subgroup $\Gamma$ of $G$, 
there exists then a natural restriction homomorphism $$ r_{\Gamma}: H^{i}_{alg}(G,V)
\ra H^{i}(\Gamma,V) $$ into the ordinary cohomology groups of $\Gamma$. 
Moreover, the groups $H^{i}_{alg}(G,V)$ may be computed by the Hochschild 
isomorphism $$ hs: H^i_{alg}(G,V) \ra H^i(\lu,V)^T \; , $$ which identifies $H^i_{alg}(G,V)$
with the $T$-invariants of the Lie algebra cohomology of the unipotent radical.
(See \cite{Hochschild} for reference.)

\paragraph{Discrete cohomology of $\Gamma$, and cohomology of $M$}
Let $M$ be a manifold, $\Gamma =\pi_{1}(M)$. There exists a
natural homomorphism $$ q: H^{i}(\Gamma,V) \ra H^i(M, {\cal V}) \, $$ 
where $\cal V$ is the local coefficient system on $M$ induced by $V$. If $M$ has contractible universal cover $X$ then $q$ is an isomorphism. Moreover $q$ is an
isomorphism on $H^1$. (See for example \cite[VII, \S7]{Brown} or
\cite[IV, \S 11]{MacLane}.) \\

The proof of Theorem  \ref{thm:GHvol} then builds on the following two 
remarks: 

\paragraph{Functoriality of the characteristic class $u^n$}

Let $\Gamma = hol(\pi_{1}(M))$ be the holonomy of a compact volume
preserving affine manifold, and $G= A(\Gamma) =U T$ its Zariski closure in
$\Aff(n)$.  As remarked in section \ref{sect:radobs}, 
the translation part $t$, defines a characteristic class $u^n$, within the
$n$-th cohomology groups with coefficients in $\bbR = \bigwedge^n \bbR^n$.   
The characteristic class $u$ (and therefore $u^{n}$) is naturally 
defined in the discrete group-, algebraic group- and Lie algebra 
cohomology theory, and, as observed in \cite{GH_2}, it is compatible
with the restriction homomorphisms and the Hochschild isomorphism
$hs$. 
 
\paragraph{Representation of $u^n$ in de Rham cohomology}
It is proved in $\cite{GH_1}$ that $q(u^n) \in H^n(M,\bbR)$ 
is represented by the parallel volume form in the de Rham cohomology
group $H^n_{dR}(M,\bbR)$ of $M$ .\\

\begin{prf}{Proof of Theorem \ref{thm:GHvol}}
Since $M$ is compact, the de Rham cohomology class $q(u^n)$
of the parallel volume form on $M$ does not vanish.
Therefore, in particular $u^n \neq 0$. By the Hochschild
isomorphism this also holds for  $u^n \in H^n(\lu,\bbR)$. 
As remarked in section \ref{sect:radobs},  $u^n$ corresponds
to the class $\bar{c}(\lu)$, which is therefore is non-vanishing.
Therefore, $U$ acts transitively
on $\bbA^n$. 
\end{prf}

%

%







\newpage 
\begin{appendix}


\section{Linear algebraic groups} \label{sect:alggroups}
We provide some background material 
and basic facts from the theory of linear algebraic groups and algebraic 
group actions.  

\subsection{Definition of linear algebraic groups}
A subgroup $\bG \leq \GL(n,\bbC)$ is a 
\emph{linear algebraic group} if it is the zero
locus of polynomial equations in its matrix entries. 
A subgroup $G \leq \GL(n,\bbR)$ is said to be a \emph{real
linear algebraic group} if it is closed with respect
to the Zariski topology on $\GL(n,\bbR)$.  A real linear 
algebraic group $G$ is thus defined as the real zero locus 
of a set of real polynomials in its matrix entries. In particular, 
$G$ is of the form $G = \bG(\bbR)$, for some linear algebraic
group $\bG$, which is defined over the reals. 
A homomorphism $\rho: \bG \ra \bL$ between linear algebraic groups is
called a \emph{morphism of algebraic groups} or algebraic homomorphism if it may be
expressed by polynomials in the matrix entries of $\bG$. 


If $G \leq \GL(n,\bbR)$ is any subgroup, we let $\ac{G}^\bbR \leq \GL(n,\bbR)$ 
denote its real Zariski closure, which is the smallest real algebraic subgroup
containing $G$. Correspondingly, let $G \leq \GL(n,\bbC)$. Then  its Zariski closure
$\ac{G}^Z \leq \GL(n,\bbC)$ is the smallest linear algebraic group containing $G$.¬

\subsection{Structure theory} 
\paragraph{Unipotent groups}
A linear algebraic group is called 
\emph{unipotent} if it is conjugate to a subgroup  of the 
group of all upper triangular matrices which 
have only $1$ as  eigenvalue. All its elements are 
unipotent matrices. Every connected subgroup of unipotent 
matrices in $\GL(n,\bbR)$ is a unipotent real linear
algebraic group. Also every unipotent  linear algebraic group
is connected. (See  \cite{Borel} or \cite{On_Vin} for more details.)

\paragraph{Reductive groups}
A subgroup $G$ of $\GL(n,\bbR)$ is called \emph{reductive}
if every $G$-invariant subspace in $\bbR^n$ has a $G$-invariant complement in 
$\bbR^n$. Main examples are compact subgroups, semisimple groups 
and (complex-) diagonalisable groups. A linear algebraic group $\bG$ is
reductive if and only if its unipotent radical $\bU$ is trivial. A connected
abelian group of semisimple elements is called an algebraic torus. 

\paragraph{Levi-splitting} 
Every linear algebraic group $\bG \leq \GL(n,\bbC)$
splits as a semidirect product of linear algebraic 
groups $\bG = \bH \bU$, where $\bH$ is reductive
and $\bU$ is the maximal unipotent normal subgroup
of $\bG$. The group $\bU$ is called the \emph{unipotent radical}
of $\bG$. The splitting induces a corresponding splitting 
of the groups of real points. In particular, every real linear 
algebraic group $G \leq \GL(n,\bbR)$
splits as a semidirect product of linear algebraic 
groups $G= H  U(G)$, where $H$ is reductive
and $U(G)$ is the unipotent radical
of $G$. See \cite{Borel, On_Vin}.

\subsection{Orbit closure} \label{sect:oclosure}

Every orbit $\bG x$ of a linear algebraic group $\bG$ on a vector space 
(or algebraic variety) contains a closed orbit in its 
Zariski closure $\ac{\bG x}^Z$, and the orbit $\bG x$
is open in its closure. This is called the \emph{closed orbit lemma}, 
see \cite[1.8]{Borel}. By the latter fact and \cite[I.10]{Mumford}, 
\emph{orbit closure in the Zariski-topology and 
in the the Euclidean (Hausdorff) topology coincide}. 
Concerning the topology of real algebraic group 
actions we have:

\begin{proposition}  \label{prop:realtop}
If  $\bG$ is a linear algebraic group defined over $\bbR$, then 
the real points $\bX_{\bbR}$ of the orbit $\bX = \bG x$ form a finite 
union of orbits of $\bG_{\bbR}$. Moreover, if $\bX$ is closed, the orbits of  
$\bG_{\bbR}$ in the real algebraic variety $\bX_{\bbR}$ are closed in 
the Euclidean topology. 
\end{proposition}
For this result, see \cite[Proposition 2.3]{BH}.
If $\bG$ is reductive
also the converse holds, cf.\  \cite{Birkes}.

\begin{proposition}[See \cite{Birkes, Borel, Rosenl}] \label{prop:unip_closed}
Let $\bU$ be a unipotent linear algebraic group, which acts on 
a vector space (or on an affine algebraic
variety). Then every orbit of $\bU$ is closed.
\end{proposition}

The analogous result holds for actions of unipotent real linear 
algebraic groups. That is,  every orbit of a unipotent real algebraic 
group is closed in the real Zariski-topology. \\


For algebraic actions, which are defined over $\bbR$, 
the preceding result generalises as follows. 
Let $\bG \leq \GL(n,\bbC)$ be a 
linear algebraic group which is defined over $\bbR$.
Recall that a torus $\bT \leq \bG$ is called $\bbR$-split if it can 
be diagonalised with respect to a real basis. 

\begin{proposition}[See \cite{Birkes}]   \label{prop:minpoints}
If $\bG$ is defined over $\bbR$, and
 the stabiliser $\bG_{x}$, of $x 
\in \bbR^n$, contains a maximal $\bbR$-split torus of $\bG$ then  
the orbit $\bG x$ is Zariski closed. 
\end{proposition} 

A corresponding result holds for real algebraic group actions. Namely, if $G_{x}$
contains the maximal connected diagonalisable subgroup $A$ of $G$, then $Gx$
is closed in the Euclidean topology. (This can be proved directly as 
a consequence of Proposition \ref{prop:unip_closed} 
and the Iwasawa decomposition (see \cite{On_Vin}) 
$G= K A N = K N A$, where $K$ is compact, $N$ is unipotent.)

\section{Lie algebra cohomology} \label{sect:Lieco}
Here we introduce Lie algebra cohomology, and relative 
Lie algebra cohomology. We compute the top cohomology group 
in relative Lie algebra cohomology with one-dimensional coefficients. 

\subsection{Definition of Lie algebra cohomology}
Let $\lie{g}$ be a Lie algebra over the reals. (We
shall use the real numbers as ground field, for convenience.)
Let $V= V_{\lambda}$ be the $\lg$-module, corresponding to an action 
$\lambda: \lg \ra \lie{gl}(V)$. 
Let $$C^k(\lg,V) = (\bigwedge^k \lg)^* \tensor {} V$$
be the module of $k$-forms on $\lg$ with values in $V$.
The Lie algebra $\lg$ acts on $C^k(\lg, V)$ by the adjoint action on 
$\bigwedge^k \lg^*$ twisted 
with $\lambda$. For $X \in \lg$, the corresponding
operator $L_{X}: C^k(\lg, V) \ra C^k(\lg, V)$ is called the Lie derivative. 
For $\omega \in C^k(\lg, V)$, we put $\iota_{X} \omega \in 
C^{k-1}(\lg, V)$ to denote the contraction with $X$. Then
the following commutation formula holds: 
\begin{equation} \label{eq:prule}
 \iota_{Y}\,  L_{X}  =  L_{X}  \, \iota_{Y}  - \iota_{[X,Y]}  \; . 
\end{equation} 

The boundary operator $d_{V}: C^k(\lg, V) \ra  C^{k+1}(\lg, V)$ 
is a differential of degree one, which commutes with the operators 
$L_{X}$. It is defined inductively by the relation   
\begin{equation} \label{eq:ibop}
 \iota_{X} \, d_{V}  + d_{V} \, \iota_{X}  = L_{X} \; . 
\end{equation} 
Recall, that  $d_{V}$ may be
explicitly computed as 
\begin{equation*} \label{eq:bop}
\begin{split} 
d_{V} \omega\, (Y_{1} \wedge \cdots \wedge Y_{{k}} ) &= \sum_{l =1} ^k (-1)^{l+1} \lambda(Y_{l})
\left( \omega(Y_{1} \wedge \cdots \wedge \hat{Y_{l}} \wedge \cdots \wedge Y_{k}) \right) \\
 &\quad 
 + \sum_{r < s} (-1)^{r+s} 
 \omega([ Y_{r}, Y_{s}] \wedge Y_{1} \wedge \cdots \wedge \hat{Y_{r}} \wedge \cdots \wedge \hat{Y_{s}} \wedge \cdots \wedge Y_{k})  \; \; . 
\end{split} 
\end{equation*}
As usual, $Z^k(\lg, V) = \{
\omega \in C^k(\lg,V)  \mid d_{V} \omega = 0 \}$ denotes the group
of $k$-cocycles, and $ B^k(\lg,V) = \{  d_{V} \eta \mid \eta \in C^{k-1}(\lg,V) \} 
\subset Z^k(\lg, V)$ the group of coboundaries. The
complex $( C^\bullet(\lg, V), d_{V})$ is called the 
\emph{Koszul-complex}. Its cohomology vector spaces 
$H^k(\lg, V_{\lambda}) =   Z^k(\lg, V) / B^k(\lg,V)$  are called
the cohomology groups of $\lg$ with coefficients in $V$.

\subsection{Relative Lie algebra cohomology} Let $\lh \subset \lg$ be
a subalgebra. We shall also consider the relative cohomology groups 
$H^k(\lg, \lh, V)$. Let $C^k(\lg/\lh,V)$ denote the subspace of cochains in
$C^k(\lg,V)$ which vanish if one argument is contained in $\lh$. 
Note that $C^k(\lg,\lh, V) \cong \Lambda^{k} (\lg/\lh^* , V)$.
As follows from \eqref{eq:ibop}, the $\lh$-invariants 
 $$ C^\bullet(\lg/\lh,V)^\lh \cong \Hom_{\lh}( \Lambda^{\bullet} \lg/\lh , V) \; , $$
are preserved by $d_{V}$, and thus form a subcomplex $C^\bullet(\lg,\lh,V)$ of the Koszul complex $(C^\bullet(\lg, V), d_{V})$.  This subcomplex is called the complex of relative cochains. Its cohomology groups are the \emph{relative cohomology groups}
$H^k(\lg, \lh, V)$.
Note, in particular, that the inclusion of cochain complexes,  induces  natural homomorphisms 
$$    H^k(\lg, \lh, V) \ra H^k(\lg, V) \; . $$

See \cite{HS,Knapp, Koszul} for detailed 
reference on Lie algebra (co-) homology.

\subsubsection{Top cohomology group with one-dimensional coefficients}
Note that $H^k(\lg,\lh,V) = \{  0\} $, $k >n =  \dim \lg/ \lh$ .
Let $\lambda: \lie{g} \ra \bbR$ be a one
dimensional representation. We shall compute the top
cohomology group $H^n(\lg, \lh, \bbR_{\lambda})$.\\

Let $\lie{n}(\lg, \lh)$ denote the normaliser of $\lh$ in $\lg$. 
Let $\ad: \lie{g} \ra \lie{gl}(\lg)$ denote the adjoint representation, 
$ \ad_{\lg/\lh}: \lie{n}(\lg, \lh) \ra \lie{gl}(\lg/\lh)$ the quotient
representation. Note, in particular, that (induced by the adjoint 
action) $\lie{n}(\lg, \lh)$ acts on the one dimensional
module $C^n(\lg/\lh,\bbR)$. Since $C^n(\lg/\lh,\bbR) \cong \Lambda^n (\lg/\lh)^*$, 
$\lie{n}(\lg, \lh)$ acts with the character $- \trace \ad_{\lg/\lh}$. 
The Lie derivatives $L_{X}$,  for 
$X \in \lie{n}(\lg, \lh)$, preserve $C^n(\lg/\lh,\bbR_{\lambda}) \subseteq  
C^n(\lg,\bbR_{\lambda})$,
and, thereby, $\lie{n}(\lg, \lh)$ acts with the character 
\begin{equation} \label{eq:chil}
 \bar{\chi}_{\lambda} =  \lambda - \trace \ad_{\lg/\lh}  
\end{equation}
on  $C^n(\lg/\lh, \bbR_{\lambda})$, that is, for all $\tau \in C^n(\lg/\lh,\bbR_{\lambda})$,
\begin{equation} \label{eq:chil2}
 L_{X} \tau =  \bar{\chi}_{\lambda} (X) \, \tau  \; . 
\end{equation}
Now the following holds: 


\begin{proposition} \label{prop:relco}
Let $n =  \dim \lg/ \lh$. Then the group
$H^n(\lg,\lh,\bbR_{\lambda})$ is non-zero if and only if 
$\bar{\chi}_{\lambda} \equiv 0$.
\end{proposition}

\begin{proof}  Assume that $H^n(\lg,\lh,\bbR_{\lambda}) \neq \{0 \}$. 
Then, in particular, there exists a non-zero generator $\tau$ of the
module $C^{n}(\lg,\lh, \lambda) =  \left(\Lambda^{n} (\lg/\lh)^*\right)^\lh$.
We compute the boundary operator 
$d_{\lambda}: C^{n-1}(\lg,\lh, \lambda) \ra C^{n}(\lg,\lh, \lambda)$: 
Let $\omega \in \left(\Lambda^{n-1} (\lg/\lh)^*\right)^\lh$ be a relative
$n-1$ cochain. By duality in $C^\bullet(\lg/\lh, \bbR)$, there exists 
$X \in \lg$, such that $\iota_{X} \tau = \omega$. Moreover, by \eqref{eq:prule}, 
$ 0= L_{H} \omega = L_{H} \iota_{X} \tau = \iota_{[H,X]} \tau $, for all $H \in \lh$.
In fact, this implies that $X \in \lie{n}(\lg, \lh)$.
Using \eqref{eq:ibop}, we compute 
\begin{equation} \label{eq:bound}
 d_{\lambda} \omega = d_{\lambda} \iota_{X} \tau = L_{X} \tau = 
\bar{\chi}_{\lambda}(X) \tau  \; .  
\end{equation}
Therefore, if $\bar{\chi}_{\lambda} \neq 0$, 
$H^n(\lg,\lh,\bbR_{\lambda}) = \{0 \}$.

For the converse,  assume $H^n(\lg,\lh,\bbR_{\lambda}) = \{0 \}$.
Then either $\left(\Lambda^{n} \lg/\lh^*\right)^\lh = \{ 0 \}$, or 
$d_{\lambda} \neq 0$. In the first case, we must have
$\bar{\chi}_{\lambda} (H) \neq 0$, for some $H \in \lh$.
In the second case, we may use the  computation of
the boundary \eqref{eq:bound}, to conclude $\bar{\chi}_{\lambda} \neq 0$.
\end{proof}

The special case,  $H^n(\lg) \neq 0$, $n = \dim \lg$, 
if and only if $\lg$ is unimodular,  is due to Koszul \cite{Koszul}.

\begin{example} If $\lg$ is nilpotent then the group
$H^n(\lg,\lh,\bbR_{\lambda})$ is non-zero
if and only if $\lambda \equiv 0$ on $\lie{n}(\lg, \lh)$.
\end{example}

For an interpretation of $H^n(\lg,\lh, \bbR_{\lambda})$ concerning invariant measures
on homogeneous spaces, see Proposition \ref{prop:biinvmeasures}. 


\section{Invariant measures on homogeneous spaces} \label{sect:invmeasures}
We briefly give some standard background material on the existence
of invariant measures on homogeneous spaces. One may
consult \cite{Raghunathan}[Chapter I] or \cite{Weil_1} for
more detailed reference.  
We then proceed to  show that the top relative Lie algebra 
cohomology group with one-dimensional coefficients, which
is associated to a semi-invariant measure on a homogeneous space,
carries information about the measure preserving automorphisms 
of the space. This extends a well known non-vanishing result of Koszul \cite{Koszul}
on the top cohomology of unimodular Lie groups.

\subsection{Semi-invariant and invariant measures} 
Let $G$ be a Lie group,
$H \leq G$ is a closed subgroup. We put $X = G/H$ 
for the associated $G$-homogeneous space, $n = \dim X$. 
Let $L_{g}: X \ra X$ denote left-multiplication with $g \in G$.
A Borel measure $\mu$ on $X$ is called \emph{semi-invariant} with character
$\lambda$ if there exists a continuous homomorphism $\lambda: G \ra \bbR^{>0}$
such that $L_{g}^* \mu = \lambda(g) \mu$, for all $g \in G$.  The measure $\mu$ is
called \emph{invariant} if $\lambda \equiv 1$.  

\paragraph{Haar measure and unimodular character}
Every locally compact group $G$ has a 
(up to scalar multiple) 
unique (left)-invariant measure $\mu =\mu_{G}$, which is called
the \emph{Haar measure} of $G$. 

Let $R_{g}: G \ra G$ denote right-multiplication
with $g \in G$. Then $ R_{g}^* \mu = \Delta_{G}(g) \mu$ is another 
Haar measure for $G$. The homomorphism 
$$ \Delta = \Delta_{G}: G \ra \bbR^{>0} $$
is  called the \emph{unimodular character} of $G$. 
If $\Delta_{G} \equiv 1$ 
the Haar-Measure is also right-invariant. Therefore, 
$G$ has a \emph{biinvariant measure}
if and only if $\Delta_{G} \equiv 1$. In this case, $G$ is called 
\emph{unimodular}. 

Since $G$ is a Lie group, the Haar measure 
can be computed by integration relative to a left-invariant $n$-form $\omega \neq 0$ 
on $G$, $n = \dim G$.  We have $$ \Delta_{G}(g) =  | \det \Ad(g) | \; , $$ 
where $ \det \Ad(g)$ is the determinant of the adjoint representation.  
In particular,  $G$ is unimodular if and only if the determinant of
the adjoint representation has absolute value one.

\paragraph{Existence of semi-invariant measures}
Not every homogeneous space admits an invariant measure.
A precise criterion is as follows:\\

For $h \in H$, define $\Delta_{G/H}(h) = 
\Delta_{G}(h) \Delta_{H}(h)^{-1}$. 

\begin{proposition}[\mbox{\cite{Raghunathan}[Lemma 1.4]}] 
\label{prop:sinvmeasure}
The homogeneous space 
$X= G/H$ admits a (unique up to scalar)
semi-invariant measure with character $\lambda$ if and only if,  
for all $h \in H$, $$ \lambda(h) = \Delta_{G/H}(h) \; . $$
\end{proposition}

In particular, $X$ admits an invariant measure if and only if
 $\Delta_{G/H} \equiv 1$ on $H$. \\
 
Let $\Ad_{G/H}$ denote the adjoint representation of $H$ on 
$\lg /\lh$.  Then,
for all $h \in H$,  $ \Delta_{G/H}(h) = | \det \Ad_{G/H}(h) |$.

\subsection{The unimodular character of $X =G/H$} 
Let $N_G(H)$ denote the normaliser of $H$ in $G$. 
The adjoint representation $G \ra  \GL(\lie{g})$
induces a quotient representation
 $$ \Ad_{G/H}: N_G(H) \ra \GL(\lie{g} /\lie{h}) \; . $$
Note that the restriction of $\Ad_{G/H}$ to $H$ corresponds
to the isotropy representation of $H$ on the tangent space of 
$X$ at $H$.
 
For every $g \in N_G(H)$, define $$ C(g) = L_{g }R_{g}^{-1}: X \ra X \; , $$
where $R_{g}$ denotes right-multiplication on $X$. 
Let $\mu$ be a semi-invariant measure on $X$  with character $\lambda$.
The relation $$ C(g)^*\mu =  \Delta_{G/H}(g) \mu , $$ defines a
unimodular character, $$  \Delta_{G/H}: N_G(H) \ra \bbR^{>0} \; , $$
independently of $\lambda$. In fact,  $\Delta_{G/H} =  | \det \Ad_{G/H} |$. 
Obviously, we have:

\begin{lemma} \label{lemma:biinvmeasures}
The semi-invariant measure $\mu$ with character $\lambda$
is right-invariant by $g \in N_G(H)$ if and only if 
$\lambda(g) \,  \Delta_{G/H}^{-1}(g) =1$.
\end{lemma} 

In addition, we obtain the following: 

\begin{proposition} \label{prop:biinvmeasures}
Let $X= G/H$ be a homogeneous space and $\mu$ a semi-invariant measure
with (smooth) character $\lambda$. Let $\bar{\lambda}: \lie{g} \ra \bbR$ denote the 
derivative of $\lambda$. 
Then the following conditions are equivalent:
\begin{enumerate}
\item $\mu$ is invariant by the right action of $N_G(H)^0$ on $X$.
\item $\chi_{\lambda} = \lambda\;  \Delta_{G/H}^{-1}: \; N_G(H)^0 \ra \bbR^{>0}   \; \, \equiv 1$.
\item  $H^n(\lg,\lh,\bbR_{\bar{\lambda}}) \neq  \{ 0 \}$. 
\end{enumerate}
\end{proposition} 
\begin{proof} This a direct consequence of the above lemma 
and Proposition \ref{prop:relco}.  
\end{proof}

\begin{example} \label{example:nilinvm} Let $G$ be nilpotent.
Then $\Delta_{G/H} \equiv 1$. Therefore, $H^n(\lg,\lh,\bbR_{\bar{\lambda}}) \neq  \{ 0 \}$
is equivalent to $\lambda \cong 1$ on $N_G(H)^0$. 
\end{example}

\begin{corollary} \label{cor:biinvmeasures}
Let $X= G/H$ be a homogeneous space which admits an invariant 
measure $\mu$. Then $\mu$ is right-invariant by the action of
$N_G(H)^0$ if and only if $H^n(\lg,\lh,\bbR) \neq  \{ 0 \}$.
\end{corollary}

As a special case of the corollary, we have the well known 
theorem of Koszul \cite{Koszul}, which
states that $H^n(\lie{g}, \bbR) \neq 0$ if and only if $G$
is unimodular.

\end{appendix}





\begin{thebibliography}{} 

\bibitem{Abels} H.\ Abels, {\em Properly discontinuous groups of affine transformations. A survey}, Geometriae Dedicata {\bf 87}  (2001), 309-333.


\bibitem{AS} S.\ Adams, G.\ Stuck, {\em 
The isometry group of a compact Lorentz manifold. I, II. }
Invent. Math. {\bf 129} (1997), no. 2, 239-261, 263-287. 

\bibitem{Aubert-Med}  A.\ Aubert, A.\ Medina, {\em 
Groupes de Lie pseudo-riemanniens plats}, 
Tohoku Math. J. (2) 55 (2003), no. 4, 487-506.


\bibitem{Auslanderc}  L.\ Auslander, {\em The structure of complete locally affine manifolds},
Topology {\bf 3} (1964) suppl.\ 1, 131-139.

\bibitem{Auslander}  L.\ Auslander, {\em Simply transitive groups of
affine motions},  Amer.\ J.\ Math.\  {\bf 99} (1977), no.\ 4, 809-826. 

\bibitem{AusM}  L.\ Auslander, L. Markus, {\em Flat Lorentz 3-manifolds}, 
Mem.\ Amer.\ Soc.\ Math.\  {\bf 30} (1959). 

\bibitem{AusKur} L.\ Auslander, M.\ Kuranishi,
{\em On the holonomy group of locally Euclidean spaces}, 
Ann.\ of Math.\ (2) {\bf 65} (1957), 411-415.


\bibitem{BauesT} O.\ Baues, {\em Left-symmetric algebras for ${\mathfrak {gl}}_n$}, 
Transactions of the AMS, {\bf 351}  (1999), no.\  7.

\bibitem{BauesG} O.\ Baues, {\em Gluing affine two-manifolds with polygons},
Geom.\ Dedicata {\bf 75}  (1999), no.\ 1, 33-56

\bibitem{BauesV} O.\ Baues, {\em Varieties of discontinuous groups}, in
{\em Crystallographic groups and their generalizations},  
Contemporary Mathematics {\bf 262} (2000), 147-158.


\bibitem{BauesF} O.\ Baues, 
{\em Finite extensions and unipotent shadows of affine
crystallographic groups},
C.\ R.\ Math.\ Acad.\ Sci.\ Paris
{\bf 335}  (2002),  no.\ 10, 785-788.

\bibitem{BC_1} O.\ Baues, V.\ Cort\' es, 
{\em Simply transitive abelian groups of symplectic type}, 
Ann.\ Inst.\ Fourier (Grenoble) {\bf 52} (2002), no. 6, 1729-1751.

\bibitem{BG} O.\ Baues, W.M.\ Goldman, {\em Is the deformation space of complete affine structures on the 2-torus smooth?}, ``Geometry and dynamics'', 69-89, Contemp.\ Math., {\bf 389}, Amer.\ Math.\ Soc., Providence, RI, 2005.

\bibitem{Bieberbach_1} L.\ Bieberbach, {\"Uber die Bewegungsgruppen
der Euklidischen R\"aume}, Math.\ Ann.\ {70} (1911), 297-336.

\bibitem{Bieberbach_2} L.\ Bieberbach, {\"Uber die Bewegungsgruppen
der Euklidischen R\"aume II}, Math.\ Ann.\ {72} (1912), 400-412.

\bibitem{BS} 
 A.\ Borel,  J.-P.\  Serre, 
{\em Th\' eor\` emes de finitude en cohomologie galoisienne},
Comment.\ Math.\ Helv.\ {\bf  39} (1964),  111-164.

\bibitem{BH}  A.\ Borel,  Harish-Chandra, {\em 
Arithmetic subgroups of algebraic groups}, 
Ann.\ of Math.\ (2) {\bf 75} (1962), 485-535.

\bibitem{Benoist} Y.\ Benoist, 
{\em Une nilvari\'et\'e non affine}, 
J.\ Diff.\ Geom.\ {\bf 41} (1995), 21-52.

\bibitem{Benzecri} J.-P.\ Benz\'ecri, {\em 
Sur les vari\'et\'es localement affines et localement projectives}, 
Bull.\ Soc.\ Math.\ France {\bf 88} (1960), 229-332. 

\bibitem{Birkes} D.\ Birkes, {\em 
Orbits of linear algebraic groups}, 
Ann.\ of Math.\ (2) {\bf 93} (1971), 459-475. 

\bibitem{Borel} A. Borel, {\em Linear algebraic groups},  
Second edition,  Graduate Texts in Mathematics {\bf 126},  
Springer-Verlag, (1991).

\bibitem{Brown} K.S.\ Brown, {\em Cohomology of groups},
Graduate Texts in Mathematics, {\bf 87}. Springer-Verlag, 
New York-Berlin, (1982).

\bibitem{Bryant} R.L.\ Bryant, {\em 
Bochner-K\"ahler metrics}, J. A.M.S. {\bf 14}  (2001), no. 3, 623-715.

\bibitem{Burde}
 D.\ Burde, 
{\em Affine structures on nilmanifolds},  Internat.\ J.\ Math.\ {\bf  7} (1996), 599-616.

\bibitem{Burde_2}
 D.\ Burde, 
{\em Left-symmetric algebras, or pre-Lie algebras in geometry and physics}, Cent. Eur.\ J.\ Math.\ {\bf 4} (2006), no. 3, 323-357.

\bibitem{BuKa} P.\ Buser, H.\ Karcher, {\em Gromov's almost flat manifolds},
 Ast\'erisque, {\bf 81},  Soci\'et\'e Math\'ematique de France, Paris, {\bf 1981}

\bibitem{Carriere} Y.\ Carri\` ere,
{\em Autour de la conjecture de L.\ Markus sur les vari\' et\' es affines},
Invent.\ Math.\ {\bf 95} (1989), 743-753.

\bibitem{CS} V.\ Cort\' es, L.\  Sch\"afer, {\em Geometric structures on Lie groups with flat biinvariant metric}, preprint March 2008

\bibitem{DSS} K.\ Dekimpe, M.\ Sadowski, A.\ Szczepa\'nski, 
{\em Spin structures on flat manifolds},  Monatsh. Math. {\bf 148} (2006), no. 4, 283-296.

\bibitem{DR}
P.\ Doyle, J.P. Rossetti, \emph{
Tetra and Didi, the cosmic spectral twins},
Geom.\ Topol.\  {\bf 8} (2004), 1227-1242. 

\bibitem{DG}
T. A.\ Drumm, W.\ Goldman, {\em Complete flat Lorentz 3-manifolds with free fundamental group}, Internat. J. Math. {\bf 1} (1990), 149-161.

\bibitem{DuIh} D.\ Duncan, E.\ Ihrig, 
{\em Incomplete flat homogeneous geometries},
Differential geometry: geometry in mathematical physics and related topics (Los Angeles, CA, 1990), 197-202, Proc. Sympos. Pure Math., 54, Part 2, Amer. Math. Soc., Providence, RI, 1993. 

\bibitem{DuIh_3} D.\ Duncan, E.\ Ihrig, 
{\em  Translationally isotropic flat homogeneous manifolds with metric signature 
$(n,2)$},   Ann. Global Anal. Geom. {\bf 11} (1993), no. 1, 3-24.

\bibitem{DuIh_2}  D.\ Duncan, E.\ Ihrig, 
{\em Flat pseudo-Riemannian manifolds with a nilpotent transitive group of isometries}, Ann.\ Global Anal.\ Geom.\ {\bf 10} (1992), no. 1, 87-101. 

\bibitem{DuIh_1}  D.\ Duncan, E.\ Ihrig, {\em Homogeneous spacetimes of zero curvature},
Proc.\  Amer.\ Math.\ Soc.\ {\bf 107}  (1989), no. 3, 785-795.

\bibitem{EMSII} {\em Lie groups and Lie algebras. II},  
Encyclopaedia of Mathematical Sciences, {\bf 21},  
Springer-Verlag, Berlin, (2000). 


\bibitem{Epstein} R.D.\ Canary, D.B.A.\ Epstein, P.\ Green, {\em
Notes on notes of Thurston}, Analytical and Geometric Aspects of Hyperbolic Space,
London Mathematical Society Lecture Note Series  {\bf  111}, 
Cambridge University Press, (1984). 

\bibitem{FaGi} J.\ Faraut, S.\ Gindikin, {\em Pseudo-Hermitian symmetric spaces of tube type}, Topics in geometry, 123-154, Progr. Nonlinear Differential Equations Appl., 20, BirkhŠuser Boston, Boston, MA, 1996.

\bibitem{Fried} D.\ Fried, {\em Flat spacetimes}, J.\ Differential Geom. {\bf 26} (1987), 385-396.

\bibitem{FriedGoldman} D.\ Fried,  W.M.\ Goldman,
 {\em Three-dimensional affine crystallographic groups},  Adv. in Math. {\bf 47}
(1983), no. 1, 1-49.

\bibitem{FGH} D.\ Fried, W.M.\ Goldman, M.W.\  Hirsch, 
{\em Affine manifolds with nilpotent holonomy},  
Comment. Math. Helv. {\bf 56}  (1981), no. 4, 487-523.

\bibitem{FurArr} P.\ Furness, D.\ Arrowsmith, 
{\em Locally symmetric spaces},
J.\ London Math.\ Soc.\ (2) {\bf 10} (1972), 487-499  

\bibitem{Goldman} W.M.\ Goldman, 
{\em Geometric structures on manifolds and varieties of representations}, 
Geometry of group representations, (Boulder, Colorado, 1987), 
Contemporary Mathematics {\bf 74} (1988), 169-198.

\bibitem{Goldman1}  W.M.\ Goldman, {\em Projective geometry on manifolds}, 
Lecture notes v.0.3 (1988).  

\bibitem{GH_2} W.M.\ Goldman,  M.W.\  Hirsch,
\emph{Affine manifolds and orbits of algebraic groups},
 Trans. Amer. Math. Soc. {\bf 295} (1986), no. 1, 175-198.

\bibitem{GH_1} W.M.\ Goldman,  M.W.\  Hirsch,
\emph{The radiance obstruction and parallel forms on affine manifolds},
Trans. Amer. Math. Soc. {\bf 286} (1984), no. 2, 629-649.
 
\bibitem{GoldmanKamishima} W.M.\ Goldman, Y.\  Kamishima, 
{\em Topological rigidity of developing maps with applications 
to conformally flat structures},
Geometry of group representations, (Boulder, Colorado, 1987), 
Contemporary Mathematics {\bf 74} (1988), 199-203 

\bibitem{GoldmanKamishima_2} W.M.\ Goldman, Y.\  Kamishima, 
{\em The fundamental group of a compact flat Lorentz space form is virtually polycyclic}, J.\ Differential Geom.\ {\bf 19} (1984), no. 1, 233-240.


\bibitem{GM} F.\ Grunewald, G.\ Margulis, 
{\em Transitive and quasitransitive actions of affine groups preserving a generalized Lorentz-structure},
J.\ Geom.\ Phys.\  {\bf 5}  (1988), no.\ 4, 493-531.


\bibitem{GS}  F.\ Grunewald, D.\ Segal, 
{\em On affine crystallographic groups},
J.\  Differential Geom.\ {\bf  40} (1994), no.\ 3, 563-594.

\bibitem{Grom_1} M.\ Gromov, \emph{Almost flat manifolds},
 J. Differential Geom. {\bf 13} (1978), no. 2, 231-241. 
 
\bibitem{Gued} M.\ Guediri, \emph{Compact flat spacetimes}, Differential Geom. Appl. {\bf 21} (2004), no. 3, 283-295.  

\bibitem{Guan} D.\ Guan, {\em On compact symplectic manifolds with Lie group symmetries},
Trans.\ Amer.\ Math.\ Soc.\ {\bf  357} (2005), no. 8, 3359-3373. 

\bibitem{Helgason} S.\ Helgason, {\em Differential geometry, Lie groups, 
and symmetric spaces}, Pure and Applied Mathematics {\bf 80},
Academic Press, Inc., 1978 

\bibitem{Helms} J.\  Helmstetter, 
{\em Radical d'une alg\`ebre sym\'etrique \`a gauche}, 
Ann. Inst. Fourier (Grenoble) {\bf 29} (1979), no. 4, viii, 17-35.

\bibitem{Hermann} R.\ Hermann, {\em 
An incomplete compact homogeneous Lorentz metric},  
J.\ Math.\ Mech.\ {\bf 13} (1964),  497-501. 

\bibitem{Hermann_2} R.\ Hermann, {\em 
Geodesics and classical mechanics on Lie groups},
J.\ Mathematical Phys.\ {\bf 13} (1972), 460-464.

\bibitem{Hirsch} M.\ W.\  Hirsch, {\em Flat manifolds and the cohomology of groups},  Algebraic and geometric topology (Proc. Sympos., Univ. California, Santa Barbara, Calif., 1977), pp. 94-103, Lecture Notes in Math., {\bf 664}, Springer, Berlin, 1978. 

\bibitem{Hochschild} G.\ Hochschild, \emph{Cohomology of algebraic linear groups},
 Illinois J.\ Math.\ {\bf 5} (1961), 492-519.

\bibitem{HS} G.\ Hochschild,  J.-P.\ Serre, \emph{Cohomology of Lie algebras},
 Ann. of Math. (2) {\bf 57}  (1953), 591-603.


\bibitem{Jo_1} K.\ Jo, \emph{Differentiability of quasi-homogeneous convex affine domains}, J.\ Korean Math.\ Soc.\ {\bf 42} (2005), no. 3, 485-498.

\bibitem{JoKiM} K.\ Jo, I.\ Kim, \emph{Convex affine domains and Markus conjecture}, Math.\ Z.\  {\bf 248} (2004), no. 1, 173-182.
 
\bibitem{Jo_2} K.\ Jo, \emph{Quasi-homogeneous domains and convex affine manifolds}, Topology Appl.\ {\bf 134} (2003), no. 2, 123-146.

\bibitem{Ka1} S.\ Kaneyuki, {\em On classification of parahermitian symmetric spaces}, Tokyo Journal of Mathematics {\bf 8} (1985), 473-482.
 
\bibitem{Ka2} S.\ Kaneyuki, {\em Compactification of parahermitian symmetric spaces and its applications}, I: Tube type realizations, "Proceedings of the III International Workshop, Lie Theory and Its Applications in Physics" (H. D. Doebner, V. K. Dobrev, and J. Hilgert, eds.), World Scientific Publishers, 2000, pp. 63-74. 

\bibitem{KaKo} S.\ Kaneyuki,  M. Kozai, {\em Paracomplex structures and affine symmetric spaces}, Tokyo Journal of Mathematics {\bf 8} (1985), 81-98. 

\bibitem{Kim} H.\ Kim, \emph{Complete left-invariant affine structures on nilpotent Lie groups}, J.\ Diff.\ Geom.\  {\bf  24} (1986), pp. 373-394.

\bibitem{Kimura} T.\ Kimura, {\em
Introduction to prehomogeneous vector spaces},
Translations of Mathematical Monographs, {\bf 215}, 
American Mathematical Society, Providence, RI, 2003.

\bibitem{Knapp} A.W.\  Knapp, {\em Lie groups, Lie algebras, and cohomology}, Mathematical Notes, {\bf 34}. Princeton University Press, Princeton, NJ, 1988. 

\bibitem{Kobayashi} S.\ Kobayashi,  ``Transformation Groups in Differential Geometry'', Springer-Verlag, New York-Heidelberg, 1972.

\bibitem{KN}   S.\ Kobayashi,  K.\ Nomizu,   ``Foundations of Differential Geometry, I'', John Wiley \& Sons 1969.     

\bibitem{KS} B.\ Kostant, D.\ Sullivan, \emph{
The Euler characteristic of an affine space form is zero},
Bull.\ Amer.\ Math.\ Soc.\ {\bf 81} (1975), no. 5, 937-938.      

\bibitem{Koszul}  J.\ L.\ Koszul, \emph{Homologie et cohomologie des alg\`ebres de Lie},
Bull.\ Soc.\ Math.\ France {\bf 78} (1950), 66-127.                                  
   
\bibitem{Koszul_1}  J.\ L.\ Koszul, \emph{Domaines born\'es homog\`enes et orbites des groupes de transformations affines},
Bull.\ Soc.\ Math.\ France {\bf 89} (1961), 515-533.                    

\bibitem{Koszul_2}  J.\ L.\ Koszul, \emph{Vari\'et\'es localement plates et convexit\'e},
Osaka J.\ Math.\  {\bf 2} (1965), 285-290.   
                     
\bibitem{Kuiper} N.H.\ Kuiper,
{\em Sur les surfaces localement affines},
G\'eom\'etrie diff\'erentielle. Colloques Internationaux du Centre National
de la Recherche Scientifique, Strasbourg 1953, 79-87.





\bibitem{MacLane} S.\  Maclane, {\em Homology},
Springer-Verlag Berlin, G\"ottingen, Heidelberg, Grundlehren der
Mathematischen Wissenschaften {\bf 114}, (1963).

\bibitem{Malcev} A. I. Malcev, {\em On a class of homogeneous spaces},
 Amer. Math. Soc. Translation {\bf 39} (1951). 
 
\bibitem{Mar}  G. A.\  Margulis, {\em Complete affine locally flat manifolds with a free fundamental group}, J.\ Soviet Math.\ {\bf 134} (1987), 129-134. 

\bibitem{Markus}  L.\ Markus, {\em Cosmological models in differential geometry}, mimeographed notes, Univ.\ of Minnesota, 1962, p.\ 58. 
 
\bibitem{Matsushima} 
Y.\ Matsushima, {\em 
On tube domains},  Symmetric spaces (Short Courses, Washington Univ., St. Louis, Mo., 1969-1970), pp. 255-270. Pure and Appl. Math., Vol. 8, Dekker, New York, 1972. 


\bibitem{Marsden} J.\ Marsden, {\em 
On completeness of homogeneous pseudo-riemannian manifolds},
Indiana Univ.\ J.\ {\bf 22} (1972/73), 1065-1066.

\bibitem{MiPo} 
R.\ Miatello, R.\ Podest\'a, \emph{Spectral properties of four-dimensional compact flat manifolds},  Ann.\ Global Anal.\ Geom.\ {\bf 29} (2006), no.\ 1, 17-50.

\bibitem{Milnor} J.\ Milnor, {\em 
On fundamental groups of complete affinely flat manifolds},
Advances in Math.\  {\bf 25} (1977), no.\ 2, 178-187

\bibitem{Milnor_2}  J.\ Milnor, {\em Curvatures of left invariant metrics on Lie groups}, Advances in Math.\  {\bf 21} (1976), no. 3, 293-329 . 

\bibitem{Milnor_3} J.\ Milnor,
{\em On the existence of a connection with curvature zero\/},
Comment.\ Math.\ Helv.\ {\bf 32} (1958), 215-223.

\bibitem{Mostow} G.D.\ Mostow,
{\em
Homogeneous spaces with finite invariant measure},
Ann.\ of Math.\ (2) {\bf 75} (1962), 17-37.

\bibitem{Mumford} D.\ Mumford, {\em 
The red book of varieties and schemes}, 
Lecture Notes in Mathematics, {\bf 1358}, 
Springer-Verlag, Berlin, (1988).  

\bibitem{NaganoYagi} T.\ Nagano, K.\ Yagi,
{\em The affine structures on the real two-torus},
Osaka J. Math.\ {\bf 11} (1974), 181-210.



\bibitem{ONeill} B.\ O'Neill, {\it Semi-Riemannian Geometry},
Academic Press, (1983). 

\bibitem{On_Vin} A.L.\ Onishchik,  \'E.B. Vinberg, {\em Lie groups and algebraic groups}, Springer Series in Soviet Mathematics. Springer-Verlag, Berlin, (1990).


\bibitem{Pf} F. Pf\"affle, {\em The Dirac spectrum of Bieberbach manifolds},
 J. Geom. Phys. {\bf 35} (4) (2000), 367-385. 

\bibitem{Scheune1} J.\ Scheuneman, 
{\em Examples of compact locally affine spaces}, 
{\em  Bull.\ Amer.\ Math.\ Soc.\ } {\bf 77} (1971),  589-592.

\bibitem{Scheune2} J.\ Scheuneman, 
{\em Affine structures on three-step nilpotent Lie algebras},
Proc.\ Amer.\ Math.\ Soc.\  {\bf 46}  (1974), 451-454.
 
\bibitem{Raghunathan} M.\ S.\ Raghunathan, {\em 
Discrete subgroups of Lie groups},  Ergebnisse der Mathematik und ihrer Grenzgebiete {\bf 68},
 Springer-Verlag,  (1972).

\bibitem{Ratcliffe} J.\ G.\ Ratcliffe,
{\em Foundations of Hyperbolic Manifolds}, 
Graduate Texts in Mathematics  {\bf 149}, Springer-Verlag, (1994).

\bibitem{Rosenl} M.\ Rosenlicht, {\em 
On quotient varieties and the affine embedding of certain homogeneous spaces},  Trans.\ Amer.\ Math.\ Soc.\ {\bf 101} (1961),  211-223. 

\bibitem{SK} M.\ Sato, T.\ Kimura, {\em A classification of irreducible prehomogeneous vector spaces and their relative invariants}, Nagoya Math. J. {\bf 65} (1977), 1-155.


\bibitem{Schie} 
A.\ Schiemann, {Ternary positive definite quadratic forms are determined by their theta series},  Math.\ Ann.\ {\bf 308} (1997), no. 3, 507-517. 



\bibitem{SuThu} D.\ Sullivan and W.\ Thurston,
{\em Manifolds with canonical coordinate charts: some examples},
L'Enseignement Math.\ {\bf 29} (1983), 15-25.

\bibitem{Tol} R.\ Tolimieri, {\em
Homogeneous space with finite invariant measure}, 
Geometriae Dedicata {\bf 1} (1972), no.\ 1, 1-5. 

\bibitem{Thurston} W.\ Thurston, {\em Three dimensional geometry and topology},
Vol.1, Princeton University Press, (1997)

\bibitem{Tralle} A.\ Tralle, K.\ Aleksy, J.\ Oprea, {\em
Symplectic manifolds with no K\"ahler structure} 
Lecture Notes in Mathematics, 1661, (1997).


\bibitem{Varadarajan} V.\ S.\ Varadarajan, 
{\em Lie groups, Lie algebras, and their representations}, 
Reprint of the 1974 edition, Graduate Texts in Mathematics {\bf 102}, 
 Springer-Verlag, (1984)

\bibitem{Vinberg} \'E.\ B.\ Vinberg, \emph{The theory of homogeneous convex cones}, Transl.\ Moscow Math.\ Soc.\ {\bf 12} (1963), 340-403. 
 
\bibitem{Wolf} J.A.\ Wolf,
{\em Spaces of constant curvature}, Publish or Perish (1984)

\bibitem{Wolf_1} J.A.\ Wolf, {\em Homogeneous manifolds of zero curvature},
Trans. Amer. Math. Soc. {\bf 104} (1962), 462-469.

\bibitem{Wolf_4} J.A.\ Wolf, {\em  Isotropic manifolds of indefinite metric},
Comment.\ Math.\ Helv.\ {\bf  39} (1964),  21-64.

\bibitem{Wolf_5} J.A.\ Wolf, {\em  On the geometry and classification of absolute parallelisms. I}, J.\ Differential Geometry {\bf 6} (1971/72), 317-342. 

\bibitem{Wolf_6} J.A.\ Wolf, {\em  On the geometry and classification of absolute parallelisms. II},  J.\ Differential Geometry {\bf 7} (1972), 19-44.

\bibitem{Wolf_2} J.A.\ Wolf,
{\em  Flat homogeneous pseudo-Riemannian manifolds}, Geom. Dedicata {\bf 57} (1995), no. 1, 111-120.

\bibitem{Wolf_3} J.A.\ Wolf,
{\em Isoclinic spheres and flat homogeneous pseudo-Riemannian manifolds},  Crystallographic groups and their generalizations (Kortrijk, 1999), 303-310, Contemp. Math., {\bf 262}, Amer.\ Math.\ Soc., Providence, RI, 2000.


\bibitem{Weil_1} A.\ Weil, {\em  L'int\'egration dans les groupes topologiques et ses applications}, Actual. Sci. Ind., no. 869. Hermann et Cie., Paris, 1940.


\bibitem{Whitney} H.\ Whitney, \emph{Elementary structure of real algebraic varieties}, 
Ann.\ of Math. (2) {\bf 66} (1957)  545-556.


\bibitem{Yagi}  K.\ Yagi, 
{\em On compact homogeneous affine manifolds}, 
Osaka J. Math. {\bf 7} (1970),  457-475.

\bibitem{ZB} Ph.\ B.\ Zwart,  W. M.\ Boothby,  {\em On compact, homogeneous symplectic manifolds},  Ann.\ Inst.\  Fourier (Grenoble) {\bf 30}  (1980), no. 1, vi-vii, 129-157.


\bibitem{Zeghib} A.\ Zeghib, 
{The identity component of the isometry group of a compact Lorentz manifold}, 
Duke Math.\ J.\  {\bf 92} (1998), no. 2, 321-333. 


\bibitem{Zimmer} 
R.\ J.\ Zimmer, {\em On the automorphism group of a compact Lorentz manifold and other geometric manifolds}, 
Invent.\ Math.\ {\bf 83}  (1986), no. 3, 411-424.


\end{thebibliography}
\end{document}